\documentclass{article}
\RequirePackage{amsmath}
\RequirePackage{hyperref}
\usepackage{00_arxiv}
\usepackage[utf8]{inputenc} 
\usepackage[T1]{fontenc}    
\usepackage{hyperref}       
\usepackage{url}            
\usepackage{pgfplots}
\pgfplotsset{compat=1.15}
\usepackage{mathrsfs}
\usetikzlibrary{arrows}
\usepackage{booktabs}       
\usepackage{amsfonts}       
\usepackage{nicefrac}       
\usepackage{microtype}      
\usepackage{lipsum}		
\usepackage{amsmath,amsfonts,amssymb}
\usepackage{graphicx}
\usepackage{subfigure}
\usepackage[super]{nth}
\usepackage{mathptmx}     
\usepackage{multirow}
\usepackage{algorithm}
\usepackage{algorithmicx}
\usepackage{algcompatible}
\usepackage{multicol}
\usepackage{algpseudocode}
\usepackage[super]{nth}
\usepackage{enumitem}
\usepackage{xcolor}
\usepackage{longtable}
\usepackage{bbm}
\usepackage{adjustbox}
\usepackage{cleveref}
\usepackage{acro}
\usepackage[section]{placeins}
\setlist{nosep} 

\newtheorem{assm}{Assumption}

\newcommand{\E}{\mathbb{E}}
\newcommand{\Q}{\mathfrak{Q}}

\newcommand{\exclude}[1]{}

\newcommand{\A}{\chi}

\usepackage{tikz}
\usepackage{soul}
\usepackage{mathtools}

\DeclarePairedDelimiter\floor{\lfloor}{\rfloor}

\definecolor{ududff}{rgb}{0.30196078431372547,0.30196078431372547,1}
\definecolor{qqqqff}{rgb}{0,0,1}
\definecolor{ffqqqq}{rgb}{1,0,0}
\definecolor{pathcolor}{rgb}{1,1,0}
\definecolor{mygray}{rgb}{96,96,96}

\usetikzlibrary{arrows,calc,decorations.pathmorphing,decorations.markings,trees,arrows.meta}
\tikzset{
  treenode/.style = {align=center, inner sep=0pt, text centered,
    font=\sffamily},
solid_circle_node/.style = {treenode, circle, line width=0.5mm, black, draw=black, minimum width=0.7cm, minimum height=0.7cm},
plain_circle_node/.style = {treenode, circle, line width=0.5mm, black, minimum width=0.6cm, minimum height=0.6cm},
solid_square_node/.style = {treenode, line width=1mm, black, draw=black, minimum width=2.9cm, minimum height=2.5cm},
plain_square_node/.style = {treenode, line width=1mm, black, minimum width=2cm, minimum height=3.5cm}
}

\title{Multi-stage Stochastic Programming Methods for Adaptive Disaster Relief Logistics Planning}

\author{
Murwan Siddig \\
  Deutsche Post Chair - Optimization of Distribution Networks\\
  RWTH Aachen University\\
  Aachen, Germany, 52072  \\ 
  \href{mailto:siddig@dpo.rwth-aachen.de}{\texttt{siddig@dpo.rwth-aachen.de}}\\
   \And
 Yongjia Song \\
  Department of Industrial Engineering\\
  Clemson University\\
  Clemson, SC, USA, 29631 \\
  \href{mailto:yongjis@clemson.edu}{\texttt{yongjis@clemson.edu}}\\
}
\usepackage{acro}
\usepackage{array}

\acsetup{
	first-style = long-short,
	list/template = longtable,
	list/sort = true}

\DeclareAcronym{2ssp}{
	short = 2SSP,
	long  = two-stage stochastic programming,
	tag = abbrev
}

\DeclareAcronym{apa}{
	short = APA,
	long  = affected potential area,
	tag = abbrev
}

\DeclareAcronym{ci}{
	short = CI,
	long  = confidence interval,
	tag = abbrev
}

\DeclareAcronym{cv}{
	short = CV,
	long  = clairvoyance,
	tag = abbrev
}

\DeclareAcronym{cvd}{
	short = CV-D,
	long  = CV model with a \textit{deterministic} time of landfall,
	tag = abbrev
}

\DeclareAcronym{cvr}{
	short = CV-R,
	long  = CV model with a \textit{random} time of landfall,
	tag = abbrev
}

\DeclareAcronym{dp}{
	short = DP,
	long  = demand point,
	tag = abbrev
}

\DeclareAcronym{dps}{
	short = DPs,
	long  = demand points,
	tag = abbrev
}

\DeclareAcronym{dpe}{
	short = DPE,
	long  = dynamic programming equations,
	tag = abbrev
}

\DeclareAcronym{dm}{
	short = DM,
	long  = decision-maker,
	tag = abbrev
}

\DeclareAcronym{fema}{
	short = FEMA,
	long  = Federal Emergency Management Agency,
	tag = abbrev
}

\DeclareAcronym{famsp}{
	short = FA-MSP,
	long  = fully adaptive MSP,
	tag = abbrev
}

\DeclareAcronym{famspd}{
	short = FA-MSP-D,
	long = FA-MSP model with a \textit{deterministic} time of landfall,
	tag = abbrev
}

\DeclareAcronym{famspr}{
	short = FA-MSP-R,
	long  = fully adaptive MSP model with a \textit{random} time of landfall,
	tag = abbrev
}

\DeclareAcronym{hdrlp}{
	short = HDRLP,
	long  = hurricane disaster relief logistics planning,
	tag = abbrev
}

\DeclareAcronym{lb}{
	short = LB,
	long  = lower bound,
	tag = abbrev
}

\DeclareAcronym{mc}{
	short = MC,
	long  = Markov chain,
	tag = abbrev
}

\DeclareAcronym{mdc}{
	short = MDC,
	long  = major distribution center,
	tag = abbrev
}

\DeclareAcronym{msp}{
	short = MSP,
	long  = multi-stage stochastic programming,
	tag = abbrev
}

\DeclareAcronym{nhc}{
	short = NHC,
	long  = National Hurricane Center,
	tag = abbrev
}

\DeclareAcronym{noaa}{
	short = NOAA,
	long  = National Oceanic and Atmospheric Administration,
	tag = abbrev
}

\DeclareAcronym{rh}{
	short = RH,
	long  = rolling-horizon,
	tag = abbrev
}

\DeclareAcronym{rh2ssp}{
	short = RH-2SSP,
	long  = Rolling-horizon two-stage stochastic programming,
	tag = abbrev
}

\DeclareAcronym{rh2sspd}{
	short = RH-2SSP-D,
	long  = RH-2SSP model with a \textit{deterministic} time of landfall,
	tag = abbrev
}

\DeclareAcronym{rh2sspr}{
	short = RH-2SSP-R,
	long  = RH-2SSP model with a \textit{random} time of landfall,
	tag = abbrev
}

\DeclareAcronym{sshws}{
	short = SSHWS,
	long  = Saffir–Simpson hurricane wind scale,
	tag = abbrev
}

\DeclareAcronym{sp}{
	short = SP,
	long  = supply point,
	tag = abbrev
}

\DeclareAcronym{sps}{
	short = SPs,
	long  = supply points,
	tag = abbrev
}

\DeclareAcronym{s2ssp}{
	short = S-2SSP,
	long  = static two-stage stochastic programming,
	tag = abbrev
}

\DeclareAcronym{s2sspd}{
	short = S-2SSP-D,
	long  = S-2SSP model with a \textit{deterministic} time of landfall,
	tag = abbrev
}

\DeclareAcronym{s2sspr}{
	short = S-2SSP-R,
	long  = S-2SSP model with a \textit{random} time of landfall,
	tag = abbrev
}

\DeclareAcronym{ub}{
	short = UB,
	long  = upper bound,
	tag = abbrev
}

\DeclareAcronym{us}{
	short = US,
	long  = United States,
	tag = abbrev
}

\begin{document}
\maketitle
\begin{abstract}
We consider a logistics planning problem of prepositioning relief items in preparation for an impending hurricane landfall. This problem is modeled as a multiperiod network flow problem where the objective is to minimize the logistics cost of operating the network and the penalty for unsatisfied demand. We assume that the demand for relief items can be derived from the hurricane's predicted intensity and landfall location, which evolves according to a Markov chain. We consider this problem in two different settings, depending on whether the time of landfall is assumed to be deterministic (and known a priori) or random. For the former case, we introduce a fully adaptive \ac{msp} model that allows the decision-maker to adjust the prepositioning decisions, sequentially, over multiple stages, as the hurricane's characteristics become clearer. For the latter case, we extend the \ac{msp} model with a random number of stages introduced in~\cite{guigues2021multistage}, to the case where the underlying stochastic process is assumed to be stage-wise dependent. We benchmark the performance of the \ac{msp} models with other approximation policies such as the static and rolling-horizon two-stage stochastic programming approaches. Our numerical results provide key insight into the value of \ac{msp}, in disaster relief logistics planning.
\end{abstract}


\section{Introduction}
\label{sec:intro}
Over the past few decades, many significant hurricanes have hit the \ac{us}. These hurricanes have resulted in heavy loss of lives, numerous personal injuries, and consequential material damages. According to the \ac{noaa}, a total of 273 hurricanes struck the \ac{us} between 1851 and 2004. Since 2005 – the year in which the costliest hurricane in history, \textit{Katrina}, made landfall in Louisiana – there have been at least 43 major hurricanes. Altogether, these major hurricanes have caused more than 600 billion dollars worth of damage in the \ac{us} and the Caribbeans, making hurricanes to be among the costliest and the most frequently observed natural disasters. 

Despite their tragic impacts, unlike other unpredictable natural disasters such as earthquakes, hurricanes can be detected a few days before they make landfall. Indeed, approximately five days before a hurricane predicted time of landfall, the \ac{nhc} provides information about the hurricane's predicted trajectory, forward speed, intensity, and endangered areas. This information comes in the form of forecast advisories, and they are often leveraged by humanitarian and governmental agencies to prepare and allocate hurricane relief supplies, such as first-aid commodities, food, water, housing, backup power generators, etc. \exclude{For instance, these forecast advisories can be incorporated into hazard analysis tools, such as PSurge~\cite{psurge}, SLOSH~\cite{slosh}, and HAZUS~\cite{hazus}, to estimate the demand for the percentage of displaced households~\cite{esnard2011index}, damage to the infrastructure~\cite{iloglu2020maximal}, and power outages~\cite{guikema2014predicting}.} Nevertheless, having accurate predictions about the hurricane characteristics is only a necessary condition for a successful disaster relief effort. As it happened, notwithstanding the unprecedented number of relief items prepositioned in preparation for Hurricane Katrina, the \ac{fema} was heavily criticized for failing to prepare enough relief items~\cite{menzel2006katrina}, after the \ac{nhc} provided accurate predictions about 48 hours before the storm hit~\cite{knabb2005tropical}. Indeed, transportation and logistics play central roles in disaster relief efforts since they deal directly with the procurement, allocation, and distribution of critical supplies and services to the affected areas~\cite{ergun2010operations,perez2016inventory}, as well as potential large-scale evacuation operations~\cite{murray2013evacuation}. All of the aforementioned assertions, coupled with the dynamic evolution of a hurricane's characteristics, emphasize the need for a reliable proactive, yet flexible, logistics operations plan. 

In this paper, we consider a typical \ac{hdrlp} problem where a \ac{dm} procures and prepositions relief items in anticipation of a hurricane that is predicted to make landfall in a few days. It is assumed that the \ac{dm} can identify an \ac{apa} that is at risk of being struck by the hurricane. Within the \ac{apa}, there is a set of \ac{dps} that represent the locations from which demand for relief items is expected. These relief items can be procured from a \ac{mdc}, and the procured items can be prepositioned at a set of \ac{sps}. Once the hurricane makes landfall, the prepositioned relief items can be shipped from the \ac{sps} to the \ac{dps}. We assume that the demand for relief items is determined based on two of the hurricane's characteristics: its intensity and its location at the time of landfall. If these characteristics are known in advance, this problem can be modeled as a \textit{deterministic} multiperiod network flow problem, defined over a logistics network (see Figure~\ref{fig:network}). In this problem, the goal of the \ac{dm} is to find an optimal solution for the logistics operations in the network, including (i) procuring the relief items from the \ac{mdc}; (ii) prepositioning the procured items at the \ac{sps}; and (iii) delivering the prepositioned items to \ac{dps} according to the incurred demand, which minimizes the overall logistics cost of procuring and shipping the relief items plus the penalty cost for failing to meet the demand for relief supplies (if any).

Since the demand for relief items is not known in advance -- but can be derived from the hurricane's characteristics which evolve according to a stochastic process with a know probability distribution -- we propose a multi-stage \textit{stochastic} programming (\ac{msp}) model for the \ac{hdrlp} problem. \ac{msp} is a class of mathematical programming models for sequential decision-making under uncertainty. In a conventional \ac{msp} problem, the goal of a \ac{dm} is to construct an optimal policy that controls the behavior of a probabilistic system over multiple decision stages as it evolves through a period known as the planning horizon. To that end, in each decision stage, the \ac{dm} observes the current state of the system and chooses an action that is prescribed by the optimal policy. \ac{msp} models have a colorful range of applications in several areas such as energy~\cite{de2017assessing}, finance \cite{dupavcova2009asset}, and transportation and logistics \cite{fhoula2013stochastic}, among others. 

We assume that the hurricane's characteristics evolve according to an \ac{mc}, where we consider not only the uncertainty about the hurricane's intensity and landfall location but also its \textit{time} of landfall. To that end, we consider the \ac{hdrlp} problem in two different settings, depending on whether the timing of the hurricane's landfall is assumed to be deterministic or to follow a probability distribution. For each setting, we propose an \ac{msp} model that allows the \ac{dm} to make adaptive logistics operational decisions over time as the characteristics of the hurricane at landfall become clearer.

Due to the uncertainty in the problem data and the sequential nature of the decision-making structure in \ac{msp} models, it is well-known that \ac{msp} models can be computationally challenging to solve.\exclude{The use of an \ac{msp} model could become handicapped under a restrictive computational budget such as a strict time limit.}To that end, we consider alternative approximate approaches which are less computationally expensive than the \ac{msp} model. In particular, we consider \ac{s2ssp} models, where the entire planning horizon is aggregated into two stages. We also consider a \ac{rh} approach where an \ac{s2ssp} model is constructed and solved in every period, but only the decisions pertaining to the current period are implemented. We evaluate the performances of these approaches in an extensive numerical experiment using synthetic data and perform various sensitivity analysis that leads to key insights into the value of MSP in hurricane relief logistics planning. 

\exclude{Although recent studies have been proposed to alleviate some of these computational burdens (see e.g.,~\cite{siddig2021adaptive,van2019level}), these enhancements often impose assumptions on the evolution of the stochastic process that are not applicable to our model. For instance, in both~\cite{siddig2021adaptive,van2019level}, it is assumed the stochastic process is assumed to be stage-wise independent. That is, any realization of the random process at a point in time $t$ is entirely independent from its history up to time $t-1$, which is not the case in an \ac{mc}.}

In sum, the contributions of our paper are threefold:
\begin{enumerate}
    \item We propose novel \ac{msp} models to the \ac{hdrlp} problem under demand uncertainty characterized by the hurricane's evolving attributes.
    \item We introduce and solve the \ac{hdrlp} problem where the hurricane's landfall time is assumed to be random in addition to its trajectory and intensity. This extends the scope of \ac{msp} modeling to address situations where the number of stages is random.
    \item We provide insights into the value of \ac{msp} models in disaster relief logistics planning compared to alternative decision policies given by the \ac{s2ssp} and the RH approach.
\end{enumerate}

The remainder of this paper is organized as follows. In Section~\ref{sec:literature}, we review the literature on applying optimization models for dealing with hurricane relief logistics problems. In Section~\ref{sec:MSP_models}, we describe the problem setting, the proposed \ac{msp} models, and the alternative approaches. In Section~\ref{sec:case_study}, we discuss our numerical experiment results. We conclude with some final remarks in Section~\ref{sec:conclusion}. 

\section{Literature Review}
\label{sec:literature}
In the literature, there is a good deal of work on applying optimization models to address different aspects of disaster relief logistics planning. A significant part of the literature is related to the last-mile delivery problems~\cite{balcik2008last}, which focus on the scheduling and delivery of relief items to the affected individuals. Other aspects of the disaster relief effort have also been studied, including the social cost of unmet demand~\cite{holguin2013appropriate,perez2016inventory,rawls2010pre,yi2007dynamic}, evacuation and rescue missions~\cite{murray2013evacuation}, disruptions in the power and supply chain systems~\cite{li2017effective,duque2020optimizing}, among others. Our work pays less attention to the details of the logistics models, but focuses on the interplay between the dynamic evolution of the hurricane characteristics and the logistics operations, using a standard multiperiod logistics network flow model. To that end, we give a brief overview of the literature most relevant to these aspects. 

\exclude{
\subsection{Hurricane Evolution}
\label{subsec:hurricane_evolution}
As previously noted, the \ac{nhc} provides data support to aid with the preparations and relief responses. When a tropical cyclone is likely to form in the Atlantic or eastern North Pacific ocean, the \ac{nhc} issues advisory notices, at least every six hours, starting approximately five days before the anticipated time of landfall. On one hand, given any of the forecast advisories, the forecast cone invoked by the most probable path, and the forecast uncertainty estimated from historical data indicate a set of vulnerable locations from which certain demand for relief items is expected. These locations are thereby considered as the \ac{dps}. On the other hand, the wind speed will reflect the level of demand quantities that ought to be expected from these locations. Specifically, depending on the category of the hurricane, as classified according to the \ac{sshws}~\cite{schott2012saffir}, one can estimate (along with the location information) the demand quantities at different \ac{dps}. It is important to note that these forecasts advisories (especially the early ones) are subject to estimation errors and misspecification. Analyzing and constructing statistical models which capture the multi-layered uncertainty about the hurricane evolution can be deployed in providing proactive decision policies~\cite{duque2020optimizing}. However, the vast majority of the literature which gathers and streams the dynamic uncertainty about the hurricane evolution into a decision model, focuses on hurricane evacuation planning~\cite{davidson2020integrated,regnier2008public,regnier2006dynamic,regnier2020six,yang2017scenario}, which is typically a one-time, irreversible decision. Instead, we focus on a fully adaptive disaster relief supply allocation policy, concurrently with the dynamic evolution of the hurricane characteristics. }

\subsection{Two-stage Stochastic Optimization Models for \ac{hdrlp}}
\label{subsec:2S_stochastic_optimization}
\exclude{The \ac{fema} operates a hierarchical nationwide disaster relief supply chain network. At the forefront of this network are the \ac{mdc}s, which are permanent storage facilities, strategically dispersed throughout the country/world. The \ac{mdc}s are followed by pre-staging and staging areas (\ac{sps}), where the relief items are temporarily stationed in anticipation of a disaster. Finally, the network has points of distributions (\ac{dps}), where the aid to the affected individuals is provided~\cite{vanajakumari2016integrated}. An efficient response to a hurricane disaster must, therefore, ensure that sufficient amounts of supplies are available at the right place and time. Deploying the aforementioned forecasting tools within the logistics of operating the network, plays an important role in meeting this objective~\cite{taskin2016bayesian,pacheco2016forecast}. Many previous studies have analyzed the logistics of operating the relief items supply network under uncertainty, at the strategic level, using \ac{s2ssp} and robust optimization models~\cite{SABBAGHTORKAN20201}. }

Most of the previous studies in the literature have focused on the logistics of operating the relief supply logistics network under uncertainty at the strategic level, using \ac{s2ssp} and robust optimization models~\cite{SABBAGHTORKAN20201}. In most of the \ac{s2ssp} models, see, e.g.,~\cite{alem2016stochastic,duran2011pre,morrice2016supporting,paul2019supply,rawls2010pre,salmeron2010stochastic,velasquez2020prepositioning}, the first-stage decisions are typically defined as preparatory actions before the hurricane makes landfall, such as deciding on the locations of facilities for stockpiling relief items, delivery fleet assignment, and inventory levels. On the other hand, the second-stage decisions are typically defined as recourse actions after observing the realized hurricane attributes. One main concern with these \ac{s2ssp} models is that it does not adequately address the estimation errors inherited in the forecast uncertainty, especially during the early stages of the planning horizon. Additionally, a latent policy, where the mobilization of relief items is postponed until the forecast uncertainty is significantly reduced, might be costly (or even infeasible) due to the limited time that remains to accomplish these tasks. This dilemma underlines the need for a policy that can be adapted to the dynamic evolution of the hurricane.

\subsection{Adaptive Decision Policies for \ac{hdrlp}}
\label{subsec:adaptive_logistics_planning}

One way to achieve an adaptive decision policy is to apply the so-called \ac{rh} approach, where, a look-ahead model is constructed and solved in every period in the planning horizon, but only the decisions pertaining to the current period are implemented. This approach naturally gives rise to an adaptive decision policy as the look-ahead model is constructed based on the most up-to-date status of the system in each period~\cite{alden1992rolling,chand2002forecast,sethi1991theory}. In the context of disaster relief logistics, this procedure has mostly been applied by using a deterministic optimization model as the look-ahead model, where the underlying stochastic process is approximated by a point estimator~\cite{vanajakumari2016integrated}. Recently, \ac{2ssp} models have been a popular choice for the look-ahead model in the \ac{rh} approach~\cite{chang2021stochastic,pacheco2016forecast,rivera2016dynamic,duque2020optimizing}. In our numerical study, we consider an \ac{rh} framework where a \ac{2ssp} model is solved in each period, as an alternative approach to obtain approximate decision policies. 

In addition to using stochastic programs as the look-ahead models, there are alternative approaches to achieve adaptive decision policies that are tailored for the application of \ac{hdrlp}. For example,~\cite{lodree2009supply,pacheco2016forecast} are among those that are most closely related to our work. They consider the problem of prepositioning relief items in preparation for a hurricane landfall. As it is assumed that earlier advisories have greater uncertainty and that the logistics cost increases over time (see the discussion following Assumption 5 in Section~\ref{sec:MSP_models}), they view the process of making the pre-disaster planning decisions as having two layers: (i) when to start the prepositioning process; and (ii) how many relief items should be prepositioned at different \ac{dps} once the prepositioning process has started. From this perspective, in the first layer, the goal of the \ac{dm} is to identify the ``best'' period to start the prepositioning process. Specifically, a period that strikes a balance between maximizing the level of confidence in the predicted hurricane characteristics and minimizing the logistics cost.~\cite{lodree2009supply} approach this problem as an optimal stopping problem, and~\cite{pacheco2016forecast} use a combination of decision theory and stochastic programming techniques. We use a similar Markovian structure for modeling the underlying stochastic process associated with the evolution of the hurricane's characteristics. However, we impose the adaptability of sequential decision-making, to the arrival of new information, more explicitly via \ac{msp} models.

\exclude{One major drawback to the model presented in~\cite{lodree2009supply} is that it views the hurricane's forecast advisories as  statistically independent of one another. This, however, is addressed in~\cite{pacheco2016forecast} where the authors model the dependency between sequential forecast advisories using an \ac{mc}. Another prominent feature of the work presented in~\cite{pacheco2016forecast} is that it uses a dynamic approach for prepositioning the relief items. Unlike the static model presented in~\cite{galindo2013prepositioning}, the dynamic approach presented in~\cite{pacheco2016forecast} allows the \ac{dm} to adapt the prepositioning decisions as new information about the hurricane arrives. To construct such dynamic policy, the authors decompose the state of the system throughout the planning horizon into three possible states: (i) an initial state where no prepositioning has been made yet; (ii) an active state where prepositioning decisions can be made; and (iii) a final state where the estimated time until the hurricane makes landfall would not be enough to perform any prepositioning operations. At first (time $t=1$), the system is assumed to be at the initial state. Then, for each subsequent period, the \ac{dm} has to decide whether to start the prepositioning operations or wait until the next forecast advisory arrives. If the \ac{dm} chooses to begin the prepositioning operations, the system evolves into the active state where the \ac{dm} needs to determine the amount and locations of the relief items to be prepositioned by using the available less accurate forecast information but with a lower logistics cost. Otherwise, the \ac{dm} waits until the next time period and faces the same dilemma again -- but now with a more accurate forecast, albeit at a higher logistics cost. Finally, once the estimated remaining time until the hurricane makes landfall becomes smaller than the time required to perform any prepositioning activity, the system evolves into the final state, and no more prepositioning is made. To solve this problem, the authors do the following. First, the author in~\cite{pacheco2016forecast} consider a set of scenarios $N$, where each scenario $k\in N$ for a hurricane occurrence is characterized by three attributes: (i) location, (ii) intensity, and (iii) time of landfall. Second, for a given scenario $k\in N$, to determine if a given time period $t_k\leq T$ is the best time to start prepositioning, the authors evaluate the cost of acting at time $t_k$ compared to waiting for an additional period (i.e., time $t_k+1$). If the cost of waiting is lower than the risk of acting, the decision of initiating the prepositioning is postponed for time $t_k+1$. Otherwise, the prepositioning is initiated at time $t_k$. To measure the difference between the cost of acting now and the risk of waiting, the authors use a decision theory approach that minimizes the expected cost and the maximum regret. In this paper, we consider a similar problem of prepositioning relief supplies like the one presented in~\cite{pacheco2016forecast}. }

\exclude{
{\color{blue} (There are also some recent works that construct scenario trees for \ac{msp} models using the sample-path forecast data. See, e.g., \url{https://onlinelibrary.wiley.com/doi/full/10.1111/risa.12990} and \url{https://www.sciencedirect.com/science/article/pii/S2212420918312457})}
}

\section{Stochastic Programming Models and Methods for Hurricane Relief Logistics Planning}
\label{sec:MSP_models}
This section is organized as follows. First, we give a general description of the \ac{hdrlp} problem, introduce some basic assumptions, and provide a mathematical formulation in a completely deterministic setting. Second, we introduce a generic formulation for \ac{msp} models. Third, we introduce the \ac{msp} model and the alternative approaches for the situation where the hurricane time of landfall is assumed to be deterministic. Finally, we extend the \ac{msp} model and the alternative approaches to the situation where the hurricane time of landfall is assumed to be random.

\subsection{The \ac{hdrlp} problem: description, assumptions, and deterministic formulation.}
\label{subsec:detrminitc_formulation}
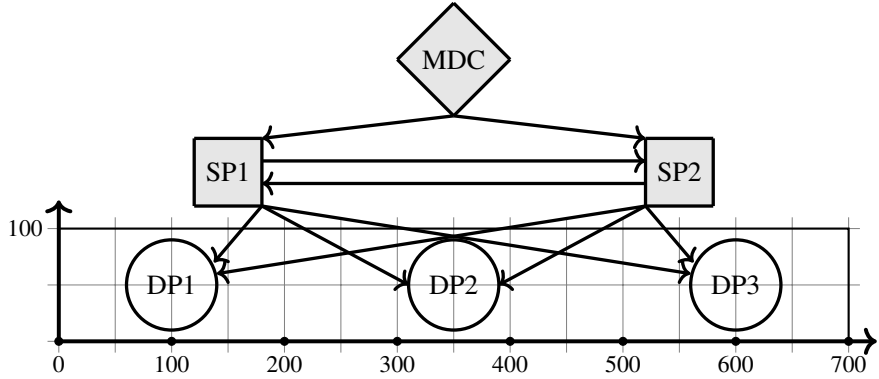
\begin{figure}[htbp]
\begin{center}
\begin{tikzpicture}[scale=1.5, transform shape]
\draw[step=0.5cm,gray,very thin] (-0.1,-0.1) grid (7.1,1.1);
\draw[->,line width=1.5pt] (0,0) -- (7.25,0);
\draw[->,line width=1.5pt] (0,0) -- (0,1.25);
\draw[thick,black]  (0,0) -- (7,0) -- (7,1) -- (0,1) -- (0,0);
\fill[line width=2pt,color=black,fill=black,fill opacity=0.1] (3.50,2.00) -- (4.00,2.50) -- (3.50,3.00) -- (3.00,2.50) -- cycle;
\fill[line width=2pt,color=black,fill=black,fill opacity=0.1] (1.20,1.80) -- (1.80,1.80) -- (1.80,1.20) -- (1.20,1.20) -- cycle;
\fill[line width=2pt,color=black,fill=black,fill opacity=0.1] (5.20,1.80) -- (5.80,1.80) -- (5.80,1.20) -- (5.20,1.20) -- cycle;
\draw [line width=1.25pt,color=black] (1.00,0.50) circle (0.4cm);
\draw [line width=1.25pt,color=black] (3.50,0.50) circle (0.4cm);
\draw [line width=1.25pt,color=black] (6.00,0.50) circle (0.4cm);
\draw [line width=1.25pt,color=black] (3.50,2.00)-- (4.00,2.50);
\draw [line width=1.25pt,color=black] (4.00,2.50)-- (3.50,3.00);
\draw [line width=1.25pt,color=black] (3.50,3.00)-- (3.00,2.50);
\draw [line width=1.25pt,color=black] (3.00,2.50)-- (3.50,2.00);
\draw [line width=1.25pt,color=black] (1.20,1.80)-- (1.80,1.80);
\draw [line width=1.25pt,color=black] (1.80,1.80)-- (1.80,1.20);
\draw [line width=1.25pt,color=black] (1.80,1.20)-- (1.20,1.20);
\draw [line width=1.25pt,color=black] (1.20,1.20)-- (1.20,1.80);
\draw [line width=1.25pt,color=black] (5.20,1.80)-- (5.80,1.80);
\draw [line width=1.25pt,color=black] (5.80,1.80)-- (5.80,1.20);
\draw [line width=1.25pt,color=black] (5.80,1.20)-- (5.20,1.20);
\draw [line width=1.25pt,color=black] (5.20,1.20)-- (5.20,1.80);
\draw [->,line width=1.25pt,color=black] (3.50,2.00) -- (1.80,1.80);
\draw [->,line width=1.25pt,color=black] (3.50,2.00) -- (5.20,1.80);
\draw [->,line width=1.25pt,color=black] (1.80,1.60) -- (5.20,1.6);
\draw [->,line width=1.25pt,color=black] (5.20,1.40) -- (1.80,1.4);
\draw [->,line width=1.25pt,color=black] (1.80,1.20) -- (3.1,0.5);
\draw [->,line width=1.25pt,color=black] (1.80,1.20) -- (5.60,0.6);
\draw [->,line width=1.25pt,color=black] (5.20,1.20) -- (5.625,0.7);
\draw [->,line width=1.25pt,color=black] (5.20,1.20) -- (3.9,0.5);
\draw [->,line width=1.25pt,color=black] (5.20,1.20) -- (1.4,0.6);
\draw [->,line width=1.25pt,color=black] (1.80,1.20) -- (1.375,0.7);
\begin{scriptsize}
\draw [fill=black] (0,0) circle (1pt);
\draw [fill=black] (1.00,0) circle (1pt);
\draw [fill=black] (2.00,0) circle (1pt);
\draw [fill=black] (3.00,0) circle (1pt);
\draw [fill=black] (4.00,0) circle (1pt);
\draw [fill=black] (5.00,0) circle (1pt);
\draw [fill=black] (6.00,0) circle (1pt);
\draw [fill=black] (7.00,0) circle (1pt);
\draw[color=black] (3.5,2.5) node {MDC};

\draw[color=black] (1.5,1.5) node {SP1};

\draw[color=black] (5.5,1.5) node {SP2};
\draw[color=black] (1,0.5) node {DP1};
\draw[color=black] (3.5,0.5) node {DP2};
\draw[color=black] (6,0.5) node {DP3};
\end{scriptsize}
\begin{tiny}
\draw[color=black] (0,-0.2) node {$0$};
\draw[color=black] (1,-0.2) node {$100$};
\draw[color=black] (2,-0.2) node {$200$};
\draw[color=black] (3,-0.2) node {$300$};
\draw[color=black] (4,-0.2) node {$400$};
\draw[color=black] (5,-0.2) node {$500$};
\draw[color=black] (6,-0.2) node {$600$};
\draw[color=black] (7,-0.2) node {$700$};
\draw[color=black] (-0.3,1) node {$100$};
\end{tiny}
\end{tikzpicture}    
\end{center}
\caption{An example hurricane relief logistics network for a $700\times100$ affected potential area.}
\label{fig:network}
\end{figure}

We consider a typical \ac{hdrlp} problem, which is modeled as a \textit{multiperiod  network flow} problem defined on a directed network $G = (V,A)$. In this network, the set of nodes $V = \{0\}\cup I \cup J$ consists of the \ac{mdc} (denoted by node $0$), the \ac{sps} (denoted by the set $I$), and the \ac{dps} (denoted by the set $J$). Assuming that the hurricane is set to make landfall at time $T$, the hurricane relief logistics process unfolds as follows. At any point in time $t=1, \ldots, T$, the \ac{dm} has the opportunity to procure relief items from the \ac{mdc} and preposition them at different \ac{sps}. Prepositioned items at any \ac{sp} $i\in I$ can be rerouted to other \ac{sps} $i'\in I$ at any time $t$. At time $T$, once the hurricane makes landfall, the prepositioned relief items can be shipped from the \ac{sps} to the \ac{dps}, according to the demand incurred at different \ac{dps}. Figure~\ref{fig:network} depicts an example of a logistics network for an \ac{apa} with $|I| = 2$ and $|J|=3$.

In this \ac{hdrlp} problem, the objective is to minimize the total cost, which consists of the logistics cost (including the procurement and transportation costs) and the penalty cost for unsatisfied demand for relief items. The description of the parameters used in our model formulation is provided in Table~\ref{tab:parameters}.
\begin{table}[htbp]
\footnotesize
    \centering
    \begin{tabular}{p{0.125\textwidth} p{0.825\textwidth}}
    \toprule
        Parameter & Description \\ \hline
        ${c^b_{ii',t}}$ & unit cost of transporting relief items from node $i \in \{0\}\cup I$ to node $i' \in I$, in period $=1,\dots,T$.\\
        ${c^h_{i,t}}$ & unit cost of holding items at SP $i\in I$ in period $t$, in period $=1,\dots,T$.\\
        ${h_t}$ & unit cost of procuring relief items from the MDC in period $t$, in period $=1,\dots,T$.\\
        ${c^a_{ij}}$ & unit cost of transporting relief items from node $i \in I$ to node $j\in J$, at the time of landfall $T$.\\
        ${p}$ & unit penalty cost for unsatisfied demand.\\
        ${q}$ & unit salvage value for unused relief items.\\
        ${d_j}$ & demand for relief items at the time of landfall at DP $j \in J$.\\
        ${u_i}$ & inventory capacity of SP $i \in I$.\\
        \bottomrule
    \end{tabular}
    \caption{Description of the parameters in the \ac{hdrlp} multiperiod network flow problem.}
    \label{tab:parameters}
\end{table}

In addition, we make the following assumptions.
\begin{assm}
\label{assm:cap}
The \ac{mdc} has unlimited capacity, whereas each SP $i\in I$ has an inventory capacity of $u_i$.
\end{assm}
\begin{assm}
\label{assm:shipment}
All shipments made at the start of a period $t$ will arrive at their destinations by the start of period $t+1$.
\end{assm}
\begin{assm}
\label{assm:demand}
The demand is only incurred at the time of the hurricane's landfall, i.e., at time $T$.
\end{assm}
\begin{assm}
\label{assm:linearity}
We do not make decisions regarding the selection of the \ac{sps}, nor the timing of their activation. That is, all the \ac{sps} are assumed to be available starting from period $t = 1$.
\end{assm}
\begin{assm}
\label{assm:cost}
The logistics costs $(c^b_{ii',t}, c^a_{ij,t}, c^h_{i,t}, h_t)$ are nondecreasing over time $t = 1,\ldots, T$.
\end{assm}

Assumption~\ref{assm:cap} is customary. Assumptions~\ref{assm:shipment},~\ref{assm:demand}, and~\ref{assm:linearity} are introduced to simplify the underlying logistics optimization model. We note that, although Assumption~\ref{assm:demand} is a simplification, the demand for relief items at different \ac{sps} relies heavily on the hurricane's intensity and its location at the time of landfall. Hence, while it could be interesting to consider the situation where demand occurs over multiple periods, as we shall see in Section~\ref{sec:case_study}, imposing Assumption~\ref{assm:demand} does not undermine the conclusions of this paper. Assumption~\ref{assm:linearity} is employed to avoid introducing binary decision variables to model the SP selection and activation, which would inflate the computational resources needed to solve the resulting multistage stochastic \textit{mixed-integer} programming problem. Assumption~\ref{assm:cost} can be interpreted as the impact caused by the imminence of the hurricane on the day-to-day operations of the logistics network. In particular, the days before the hurricane landfall are typically associated with damage-inducing weather conditions (e.g., high winds, heavy rain, etc.), which may disrupt the logistics network. These disruptions can ultimately lead to a surge in the prices of relief commodities~\cite{beatty2021hurricanes} that are essential to the relief logistics efforts. Moreover, the imminence of a hurricane usually necessitates the need for urgent evacuation efforts by vehicles that can be deployed at short notice (e.g., general aviation aircraft). While such means of transportation have significant benefits, they often come with hefty prices. Nonetheless, as we show in Section~\ref{sec:case_study}, our proposed models and methods will remain valid even when Assumption~\ref{assm:cost} is lifted. The introduction of this assumption is, rather, to bridge the gap between the observations in our numerical experiments and the practical intuitions of \ac{hdrlp}.

\paragraph{Deterministic formulation.} The preamble to the stochastic optimization models presented in this section is the \textit{deterministic} version of the problem. We may view this deterministic version from the perspective of a clairvoyant, i.e., the demand for relief items and the hurricane's landfall time are both known before the logistics planning. We refer to the optimal solution to this deterministic optimization problem as the \ac{cv} solution. To formulate this deterministic multiperiod network flow problem for \ac{hdrlp}, we consider the decision variables shown in Table~\ref{tab:decision_variables}.
\begin{table}[htbp]
\footnotesize
    \centering
    \begin{tabular}{p{0.225\textwidth} p{0.725\textwidth}}
    \toprule
        Decision variable & Description \\ \hline
        $x_t = (x_{i,t})_{i\in I}$ & amount of inventory at SP $i\in I$, at the end of period $t,\, \forall \;t=1, \ldots, T$.\\
        $f_t = (f_{ii',t})_{i\in \{0\} \cup I, i' \in I}$ & amount of relief items shipped from the \ac{mdc} or SP $i\in \{0\} \cup I$ to $i'\in I,\, \forall \;t=1,\ldots, T$.\\
        $y = (y_{ij})_{i\in I, j \in J}$ & amount of relief items shipped from SP $i\in I$ to DP $j\in J$, at the time of landfall $T$.\\
        \bottomrule
    \end{tabular}
    \caption{Description of the decision variables in the \ac{hdrlp} multiperiod network flow problem.}
    \label{tab:decision_variables}
\end{table}

Given the time of landfall $T$, demand levels $(d_j)_{j \in J}$ and initial inventory levels $(x_{i,0})_{i\in I}$ at the \ac{sps}, the corresponding deterministic multiperiod  disaster relief problem is given by: 
\begin{subequations}
    \label{eq:det_form}
    \begin{align}
    & \begin{multlined}
    \displaystyle \min \quad \displaystyle \sum_{t=1}^T \left(\sum_{i\in \{0\}\cup I}\sum_{i'\in I}c^b_{ii',t}f_{ii',t}+\sum_{i\in I}c^h_{i,t}x_{i,t} + h_t\sum_{i\in I}f_{0i,t} \right)\\* 
    \displaystyle \quad + \sum_{i\in I}\sum_{j\in J}c^a_{ij}y_{ij} + \sum_{j\in J} p(d_{j} - \sum_{i\in I }y_{ij}) + \sum_{i\in I} q(x_{i,T}- \sum_{j\in J}y_{ij})
    \end{multlined} \nonumber \\* 
    \text{s.t.} \quad & \displaystyle x_{i,t-1} + \sum_{i'\in \{0\}\cup I, i'\neq i}f_{i'i,t} - \sum_{i'\in I, i'\neq i}f_{ii',t} = x_{i,t} & \forall \;i\in I, \; t = 1,\ldots, T, \label{eq:det_form_cos-1} \\*
     & \displaystyle \sum_{i'\in I, i'\neq i}f_{ii',t} \leq x_{i,t-1} & \forall \;i\in I, \; t = 1,\ldots, T, \label{eq:det_form_cos-2} \\*
     & \displaystyle 0 \leq x_{i,t} \leq u_i & \forall \;i\in I, \; t = 1,\ldots, T, \label{eq:det_form_cos-3} \\*
      & \displaystyle \sum_{j\in J}y_{ij} \leq x_{i,T} & \forall \;i\in I, \label{eq:det_form_cos-4} \\*
      & \displaystyle \sum_{i \in I}y_{ij} \leq d_j & \forall \;j\in J, \label{eq:det_form_cos-5} \\*
      & \displaystyle y_{ij} \geq 0 & \forall \;i\in I, \; j\in J, \label{eq:det_form_cos-6}\\*
      & \displaystyle f_{ii',t} \geq 0 & \forall \;i\in I, i'\in I,\;  t = 1,\ldots, T. \label{eq:det_form_cos-7}
   \end{align}
\end{subequations}
The first four terms in the objective function of formulation~\eqref{eq:det_form} represent the logistics costs. In particular, the first term represents the cost of shipping relief items from the \ac{mdc} to the \ac{sps}, and rerouting the relief items between \ac{sps}; the second term represents the cost of holding the relief items at the \ac{sps}, the third term represents the cost of procuring relief items from the \ac{mdc}, and the fourth term represents the cost for shipping relief items from the \ac{sps} to the \ac{dps} at the time of landfall $T$. Additionally, the fifth term represents the penalty cost for unsatisfied demand at the \ac{dps}, and the last term represents the return from salvaging unused relief items. 

Constraint~\eqref{eq:det_form_cos-1} represents the flow balance constraint, which indicates that the amount of relief items at SP $i \in I$ at time $t=1, \dots, T$ is given by the amount of items stored in SP $i$ from the previous time period, plus the items shipped from other \ac{sps} ($\forall \; i' \in I, i'\neq i$) to SP $i$ and the \ac{mdc} at the current period, and minus the items shipped from SP $i$ to other \ac{sps} at the current period. Constraint~\eqref{eq:det_form_cos-2} indicates that the total amount of items shipped from SP $i \in I$, at time $t=1, \dots, T$, to other \ac{sps} ($\forall \; i' \in I, i'\neq i$) cannot exceed the initial inventory level at the start of the current period. Constraint~\eqref{eq:det_form_cos-3} indicates that the amount of relief items stored in each SP $i \in I$ is nonnegative, and cannot exceed its capacity at any time. Constraint~\eqref{eq:det_form_cos-4} indicates that the total amount of items shipped from SP $i \in I$ at time period $t=1, \dots, T$ to all \ac{dps} ($\forall \;j \in J$) cannot exceed the inventory level at SP $i$ at the end of the same period. Constraint~\eqref{eq:det_form_cos-5} indicates that the total amount of items shipped to DP $j \in J$ should not exceed its demand level. Finally, constraints~\eqref{eq:det_form_cos-6} and~\eqref{eq:det_form_cos-7} impose the sign restrictions on decision variables $y$ and $f_{t}, \forall \; t=1, \dots, T$.

\subsection{Generic MSP models}
The proposed stochastic optimization models and methods are built based on the deterministic multiperiod network flow model~\eqref{eq:det_form}. Before we discuss them in the context of the \ac{hdrlp} problem, we first introduce a generic formulation for the \ac{msp} model.

\paragraph{Generic \ac{msp} formulation.} The starting point of the \ac{msp} models introduced for the \ac{hdrlp} problem is the following generic \textit{nested} formulation.
\begin{equation}
\label{eq:MSP}
\min_{a_1 \in \A_1(a_0,\xi_1)} z_1(a_1, \xi_{1})+ \E_{|\xi_{[1]}}\left[\rule{0cm}{0.4cm} \min_{a_2\in \A_2(a_1,\xi_2)} z_2(a_2,\xi_{2}) + \E_{|\xi_{[2]}} \left[\rule{0cm}{0.4cm}\cdots + \E_{|\xi_{[T-1]}}\left[\rule{0cm}{0.4cm} \min_{a_T\in \A_T(a_{T-1},\xi_T)}z_T(a_T,\xi_{T}) \right]\right]\right].
\end{equation}
As customary, we assume that vectors $a_0$ and $\xi_1$ are given as input data. Moreover, $\xi_{[t]} := (\xi_1, \xi_2, \dots, \xi_t)$ denotes the history of the stochastic process up to time $t$, and $\xi_{t}$ is a random vector with a known probability distribution, supported on a set $\Xi_{t} \subset \mathbb{R}^{n_{t}}$, for $t = 2,\ldots, T$. We also make the following assumption on the stochastic process $\{\xi_t\}_{t=1}^T$. 
\begin{assm}
\label{assm:markovian_process}
The stochastic process $\{\xi_t\}_{t=1}^T$ follows a discrete-time \ac{mc}, that is, $\mathbb{P}(\xi_{t} \mid \xi_1, \xi_2, \ldots, \xi_{t-1}) = \mathbb{P}(\xi_{t} \mid \xi_{t-1})$. Additionally, the state space $\Xi_t$ is finite and the one-step transition probability is given by a transition probability matrix $\mathbf{P}$, where the $(k,k')$-entry is defined by $\mathbf{P}_{k,k'} := \mathbb{P}(\xi_{t+1} = k' \mid \xi_{t} = k), \ \forall \; t=1,\dots,T$.
\end{assm}

Under Assumption~\ref{assm:markovian_process}, the conditional expectation $\E_{|\xi_{[t]}}[\cdot]$ will be replaced by $\E_{|\xi_{t}}[\cdot]$. With a slight abuse of notation, we denote the conditional probability of $\xi_{t+1}$ given $\xi_t$ as $\mathbf{P}_{\xi_{t+1}\mid \xi_t}$. 

In formulation~\eqref{eq:MSP}, the goal of the \ac{dm} is to optimize a policy, which is a mapping from the state of the system to an action. In particular, at every decision period $t=1, \dots, T$, the \ac{dm} makes a decision $a_t$ chosen from the feasible set $\A_t(a_{t-1},\xi_t)$. The decision vector $a_t$ can be split into two types of decisions: a vector of state variables $u_t$ and a vector of local (control) variables $v_t$, i.e., $a_t := (u_t, v_t)$. The vector of state variables $u_t$ is what links different stages together, and the vector of local variables $v_t$ participates exclusively in the optimization problem defined for each period $t$. Under Assumption~\ref{assm:linearity}, the feasible set $\A_t(a_{t-1},\xi_t)$ in each period $t$ is given by the set of linear constraints:
\begin{equation}
\label{eq:linear_constraints}
\A_t(a_{t-1},\xi_t) = \left\{a_t = (u_t,v_t) \mid A^{\xi_t}_{t}{u_{t}} + B^{\xi_t}_t{u_{t-1}} + C^{\xi_t}_t{v_t} = b^{\xi_t}_t\right\}.
\end{equation}
As such, the decision-making process in formulation~\eqref{eq:MSP} proceeds as follows. In each period $t=1, \dots, T$, the \ac{dm} observes the state of the system $s_t := S_t(u_{t-1}, \xi_t)$, which depends on: (i) the action from the previous period $a_{t-1}$ (and more specifically $u_{t-1}$), and (ii) the realization of random vector $\xi_t$ at the current period. When $t<T$, the \ac{dm} aims to find an optimal action $a_t$ which minimizes the {immediate} cost given by $z_t(a_t,\xi_t)$, plus the \textit{expected} {future} cost which is given by the nested expectation $\E_{|\xi_{t}}[\cdot]$. Finally, when $t=T$, no more future cost is to be paid, and the \ac{dm}'s goal is to minimize the {immediate} cost $z_T(a_T,\xi_T)$ only. A symbolic representation of the decision-making procedure in \ac{msp} models is shown in Figure~\ref{fig:MSP_models} in the Appendix.

We refer to formulation~\eqref{eq:MSP} as the \ac{famsp} model, as the optimal action $a_t$ can be adapted to the state of the system $s_t$ in \emph{every period} $t=1, \dots, T$. We care to note that, when referring to an \ac{famsp} formulation, we use the terms ``period'' and ``stage'' interchangeably. However, in other settings, as we describe later, a single stage may include multiple periods: decisions made in a single stage, which has the same level of adaptability, can include decisions for multiple periods in general. 

\paragraph{Dynamic programming equations.} A typical approach to proceed with the computation of an optimal policy for the \ac{famsp} model~\eqref{eq:MSP} is to use the so-called Bellman \ac{dpe}~\cite{bellman}: 
\begin{equation}
\label{eq:FOSDDP}
Q_t({a_{t-1}},\xi_{t}):= \displaystyle \underset{{a_t}}{\min}\left\{z_t({a_t},\xi_t) + \displaystyle \Q^{\xi_t}_{t+1}(a_t) \; \middle\vert \; a_t \in \A_t({a_{t-1}},\xi_{t})\right\}\quad \forall \; t=1,\dots, T,
\end{equation}
where $\Q^{\xi_t}_{t+1}(a_t)$ is referred to as the expected \textit{cost-to-go} function and is given by $\Q^{\xi_t}_{t+1}(a_t):= \E\left[Q_{t+1}(a_{t},\xi_{t+1}) \mid \xi_t\right], \forall \; t=1,\dots, T-1$ with $\Q^{\xi_T}_{T+1}(a_T) := 0$.
Since by Assumption~\ref{assm:markovian_process} the state space $\Xi_t$ is finite, we can write the expected cost-to-go function as $\Q^{\xi_t}_{t+1}(a_t) = \sum_{\xi_{t+1} \in \Xi_{t+1}} \mathbf{P}_{\xi_{t+1} \mid \xi_t} Q_{t+1}(a_{t},\xi_{t+1}), \forall \; t=1,\dots, T-1$ with $\Q^{\xi_T}_{T+1}(a_T) := 0$. We care to mention that \ac{msp} models where the associated stochastic process $\{\xi_t\}_{t=1}^T$ follows a discrete-time \ac{mc} are also referred to as the Markov-chain \ac{msp}, see, e.g.,~\cite{ding2019python,dowson2020policy}, for further discussions on this topic. A brief overview, and an algorithmic description, for how to solve the \ac{dpe}~\eqref{eq:FOSDDP}, using \textit{nested Benders decomposition}, is given in Section~\ref{sec:nested_benders_decompositions} in the Appendix.

\subsection{Stochastic Programming Models with Deterministic Landfall Time}
\label{subsec:determinsticT_MSP}
In this subsection, we assume that the hurricane landfall time $T$ is known a priori. First, we describe how the demand quantities at different \ac{dps} are determined in our model, and describe the stochastic process $\{\xi_t\}_{t=1}^T$ in the context of the \ac{hdrlp} problem. Second, we complete the description of the \ac{famsp} formulation by juxtaposing the components of the deterministic multiperiod \ac{hdrlp} problem~\eqref{eq:det_form} with the generic formulation provided in~\eqref{eq:MSP}. Finally, we discuss the alternative approximation policies.

\subsubsection{Modeling hurricane evolution and demand estimates}
\label{subsubsec:estimating_d_and_modelling_xi}

We assume that the demand for relief items at different \ac{dps} depends on two attributes of a hurricane at the time of its landfall: the \textit{intensity} $\alpha_T \in \mathrm{A}$ and the \textit{location} $\ell_T \in \mathrm{L}$. We also assume that this dependency structure is given by a deterministic mapping $D: \mathrm{A} \times \mathrm{L} \to \mathbb{R}^{|J|}$, that is, $d^{\xi_T} := (d^{\xi_T}_1, \dots, d^{\xi_T}_{|J|}) = D(\alpha_T, \ell_T)$. It is intuitive to assume that \ac{dps} within closer proximity to the hurricane's landfall location $\ell_T$ will prompt higher demand for relief items and vice versa. By contrast, hurricanes with higher intensity levels $\alpha_T$ will prompt higher demand for relief items and vice versa. As such, we assume that this mapping $D(\alpha_T, \ell_T)$ is \textit{increasing} in $\alpha_T$ and \textit{decreasing} in the distance between $\ell_T$ and the locations of \ac{dps}. The specifics of this function are provided, along with the remaining problem data, in Section~\ref{sec:problem_data} in the Appendix.

To model the stochastic evolution of the hurricane characteristics as an \ac{mc}, we assume that the hurricane's location and intensity are characterized by two \textit{independent} random variables $\ell_t \in \mathrm{L}$, $\alpha_t \in \mathrm{A}$ with finite state space, such that $\xi_t := (\alpha_t, \ell_t)$ and $\Xi_t := \mathrm{A} \times \mathrm{L}$, where $|\Xi_t|<\infty \forall \; t=1, \dots, T$. We further assume that $\ell_t$ and $\alpha_t$ follow an \ac{mc}. We denote the conditional probability of $\alpha_t$ given $\alpha_{t-1}$ as $\mathbf{P}_{\alpha_{t} \mid \alpha_{t-1}}$, and the conditional probability of $\ell_t$ given $\ell_{t-1}$ as $\mathbf{P}_{\ell_{t} \mid \ell_{t-1}}$. Consequently, $\mathbf{P}_{\xi_{t} \mid \xi_{t-1}} = \mathbf{P}_{\alpha_{t} \mid \alpha_{t-1}} \times \mathbf{P}_{\ell_{t} \mid \ell_{t-1}}$. We use matrices $\mathbf{P}^\alpha$ and $\mathbf{P}^\ell$ to denote the one-step transition probability matrices associated with the intensity $\alpha$ and location $\ell$, respectively.

\subsubsection{FA-MSP model with a deterministic time of landfall}
\label{subsubsec:determinsticT_MSP_stochastic_FA}

Following the definitions of the decision variables in~\eqref{eq:det_form}, at any point in time $t=1, \dots,T$, the state of the system in the \ac{hdrlp} problem is characterized by two components: (i) the informational state $\xi_t$ representing the hurricane attributes, and (ii) the physical state $x_{t-1}$ representing the relief item inventory levels at \ac{sps}. Hence, the state variables are given by $u_t := x_t, \forall \; t=1, \dots, T$ and the local (control) variables are given by $v_t := f_t, \forall \; t=1, \dots, T-1$, $v_T := (f_T,y)$ at the landfall stage $T$. That is, $a_t:= (x_t,f_t), \forall \; t=1, \dots, T-1,$ and $a_T := (x_T,f_T, y)$. Consequently, the policy is a mapping from the state of the system $(\xi_t, x_{t-1})$ to action $(x_t, f_{t})$ for $t=1, \dots, T-1$, and from $(\xi_T, x_{T-1})$ to $(x_T, f_{T}, y)$ for the landfall stage $t = T$.

We can then juxtapose the components of the \ac{cv} deterministic formulation~\eqref{eq:det_form} with the remaining components of formulation~\eqref{eq:MSP}; namely, $z_t(a_t,\xi_t)$ and $\A_t(a_{t-1},\xi_t)$. To that end, we remark that based on Assumption~\ref{assm:demand}, since the demand only occurs at the landfall stage $t = T$, one can remove the dependence on $\xi_t$ from $z_t(a_t,\xi_t)$ and $\A_t(a_{t-1},\xi_t)$ for every non-terminal stage $t = 1, \dots, T-1$, such that:
\begin{equation}
\label{eq:objective_det_Tminus}
z_t(a_t,\xi_t) := \displaystyle \sum_{i\in \{0\}\cup I}\sum_{i'\in I}c^b_{ii',t}f_{ii',t}+\sum_{i\in I}c^h_{i,t}x_{i,t} + h_t\sum_{i\in I}f_{0i,t},
\end{equation}
and
\begin{equation}
\label{eq:feasible_set_det_nonT}
    \A_t({a_{t-1}},\xi_{t}) :=
    \begin{cases} 
         \displaystyle x_{i,t-1} + \sum_{i'\in \{0\}\cup I, i'\neq i}f_{i'i,t} - \sum_{i'\in I, i'\neq i}f_{ii',t} = x_{i,t} & \forall \;i\in I, \\
         \displaystyle \sum_{i'\in I, i'\neq i}f_{ii',t} \leq x_{i,t-1} & \forall \;i\in I,\\
         \displaystyle 0 \leq x_{i,t} \leq u_i & \forall \;i\in I.
    \end{cases}
\end{equation}
At the terminal (landfall) stage $t=T$, we have:  
\begin{align}
z_T(a_T,\xi_T) &:= \displaystyle \sum_{i\in \{0\}\cup I}\sum_{i'\in I}c^b_{ii',T}f_{ii',T}+\sum_{i\in I}c^h_{i,T}x_{i,T} + h_T\sum_{i\in I}f_{0i,T}+\sum_{i\in I}\sum_{j\in J}c^a_{ij}y_{ij}\notag\\
&\quad \displaystyle + \sum_{j\in J} p(d^{\xi_T}_{j} - \sum_{i\in I }y_{ij}) + \sum_{i\in I} q(x_{i,T}- \sum_{j\in J}y_{ij}), \label{eq:objective_det_T}
\end{align}
and 
\begin{equation}
\label{eq:feasible_set_det_T}
\A_T({a_{T-1}},\xi_{T}) :=
    \begin{cases} 
        \displaystyle x_{i,T-1} + \sum_{i'\in \{0\}\cup I, i'\neq i}f_{i'i,T} - \sum_{i'\in I, i'\neq i}f_{ii',T} = x_{i,T} & \forall \;i\in I, \\
        \displaystyle \sum_{i'\in I, i'\neq i}f_{ii',T} \leq x_{i,T-1} & \forall \;i\in I,\\
        \displaystyle 0 \leq x_{i,T} \leq u_i & \forall \;i\in I,\\
        \displaystyle \sum_{j\in J}y_{ij} \leq x_{i,T} & \forall \;i\in I,\\
        \displaystyle \sum_{i \in I}y_{ij} \leq d^{\xi_T}_{j} & \forall \;j\in J,\\
        \displaystyle y_{ij} \geq 0 & \forall \;i\in I, j\in J.
    \end{cases}
\end{equation}
We refer to the \ac{famsp} model introduced in this section, hereafter, as the \ac{famspd}, to distinguish it from the one that we introduce later in Subsection~\ref{subsec:randomT_MSP}, where the hurricane's time of landfall $T$ is random.

\subsubsection{Alternative decision policies}
\label{subsubsec:alternative_approx_policies_T}

\ac{msp} models are well known to be computationally challenging to solve. These computational challenges might especially be of concern when the solution effort is restricted by a computational budget (e.g., when a strict time limit is imposed). Consider a situation where the forming of a tropical storm was not detected in its early stages and, consequently, the initial forecast advisories were overdue. In this case, it is vital to have a model that can provide a ``good'' disaster relief logistics policy within a reasonable amount of time. In the remainder of this subsection, we discuss two alternative decision policies for mitigating some of the computational burdens of solving fully-adaptive \ac{msp} models. The first is the \ac{s2ssp} model, where the logistics decisions are obtained by solving a \textit{single} \ac{2ssp} model and implemented statically. The second is the \ac{rh2ssp} approach, where an online policy is obtained by solving a \textit{sequence} of \ac{2ssp} models and implementing the corresponding solutions in an \ac{rh} fashion.

\paragraph{Two-stage approximation.} As mentioned earlier, two terms that are often used interchangeably in \ac{msp} models are \textit{stages} and \textit{periods}. Both terms refer to different time points in the planning horizon $\{1, \dots, T\}$ of an \ac{msp} problem. However, the two terms have the following subtle conceptual distinction which is of significant relevance to the \ac{2ssp} model. We use the term ``period'' to refer to a point in time when a new realization of the stochastic process is revealed; and we use the term ``stage'' to refer to a point in time where the \ac{dm} has the opportunity to adapt their decisions according to newly observed realizations of the stochastic process. 
As we discussed earlier, the two terms coincide in the \ac{famspd} model~\eqref{eq:FOSDDP} since the \ac{dm} can adapt their decisions in every period, gradually, as new information arrives. In other words, the decision-making process in the \ac{famspd} model takes the following form: $(u_1, v_1) \xrightarrow[]{\; \xi_2\;} (u_2, v_2) \xrightarrow[]{\; \xi_3\;} \quad \dots \dots \dots \quad \xrightarrow[]{\; \xi_T\;} (u_T, v_T)$. 

However, the level of adaptability in the model may be restricted to a limited number of stages. In such cases, a single-stage would consist of making decisions concerning multiple periods which are \textit{aggregated} together. In the most extreme case, this would correspond to a \ac{2ssp}, where the adaptability is restricted to two stages only. For instance, in the context of the \ac{hdrlp} problem, this would be a situation where the \ac{dm} has to decide, in advance, the amounts of relief items to be prepositioned at different \ac{sps} over the entire planning horizon. These decisions are considered first-stage decisions. Then, once the entire sample path of the stochastic process is revealed, the \ac{dm} makes a recourse decision on shipping the relief items from the \ac{sps} to the \ac{dps}, paying the penalty for unmet demand (if any), or salvaging the unused relief supplies (if any). These decisions are considered second-stage decisions. In this case, unlike the \ac{famspd} model, the prepositioning decisions are {static} because the \ac{dm} is unable to adapt their decisions when new information arrives. The corresponding decision-making process in the \ac{s2ssp} model takes the following form: $(u_1, \dots, u_T, v_1) \; \; \xrightarrow[]{\quad (\xi_2, \dots, \xi_T)\quad} \; \; (v_2, \dots, v_T)$. As such, given the initial data $x_0$ and $\xi_1$, we define the first-stage and second-stage problems in the \ac{s2ssp} model as follows.

The \textit{first-stage} problem: 
\begin{subequations}
\label{eq:twostage_detT_1st}
    \begin{align}
        \displaystyle \underset{x, f}{\min} \ & \displaystyle \sum_{t=1}^{T-1} \left(\sum_{i\in \{0\}\cup I}\sum_{i'\in I}c^b_{ii',t}f_{ii',t}+\sum_{i\in I}c^h_{i,t}x_{i,t}\right) + \sum_{t=1}^{T-1}h_t\sum_{i\in I}f_{0i,t} + \Q^{\xi_1}(x_{T-1})& \nonumber \\ 
        \text{s.t.} \quad & \displaystyle x_{i,t-1} + \sum_{i'\in \{0\}\cup I, i'\neq i}f_{i'i,t} - \sum_{i'\in I, i'\neq i}f_{ii',t} = x_{i,t}, & \forall \;i\in I,  \forall \;t = 1,\ldots, {T-1} \label{eq:twostage_detT_1st-1}\\
        & \displaystyle \sum_{i'\in I, i'\neq i}f_{ii',t} \leq x_{i,t-1}, & \forall \;i\in I, \ \forall \;t = 1,\ldots, {T-1} \label{eq:twostage_detT_1st-2}\\
        & \displaystyle 0 \leq x_{i,t} \leq u_i, & \forall \;i\in I, \forall \;t = 1,\ldots, {T-1}. \label{eq:twostage_detT_1st-3}
    \end{align}
\end{subequations}
Here, $\Q^{\xi_1}(x_{T-1}) := \sum_{\xi_T \in \Xi_T} \mathbf{P}^{T-1}_{\xi_{T} \mid \xi_1} Q(x_{T-1},\xi_T)$ gives the second-stage value function, where $\mathbf{P}^{T-1}$ gives the $(T-1)$-step transition probability matrix associated with $\xi$, and thus $\mathbf{P}^{T-1}_{\xi_{T} \mid \xi_1}$ corresponds to the conditional probability of $\xi_T$ given initial state $\xi_1$.

The \textit{second-stage} problem:
\begin{subequations}
\label{eq:twostage_detT_2nd}
    \begin{align}
        \displaystyle Q(x_{T-1}, \xi_T):= \underset{x_T, f_{T}, y}{\min} \quad & \displaystyle h_T\sum_{i\in I}f_{0i,T} + \sum_{i\in I} q\left(x_{i,T}- \sum_{j\in J}y_{ij}\right)+\sum_{i\in I}\sum_{j\in J}c^a_{ij}y_{ij} & \nonumber \\ 
        & +\sum_{i\in I}c^h_{i,T}x_{i,T} +\displaystyle \sum_{i\in I}\sum_{i'\in I}c^b_{ii',T}f_{ii',T}+\sum_{j\in J} p\left(d^{\xi_T}_{j} - \sum_{i\in I}y_{ij}\right) & \nonumber \\ 
        \text{s.t.} \quad & \displaystyle x_{i,T-1} + \sum_{i'\in \{0\}\cup I, i'\neq i}f_{i'i,T} - \sum_{i'\in I, i'\neq i}f_{ii',T} = x_{i,T}, & \forall \;i\in I, \label{eq:twostage_detT_2nd-1}\\
        &\quad \displaystyle \sum_{i'\in I, i'\neq i}f_{ii',T} \leq x_{i,T-1}, & \forall \;i\in I, \label{eq:twostage_detT_2nd-2}\\
        &\quad \displaystyle 0 \leq x_{i,T} \leq u_i, & \forall \;i\in I, \label{eq:twostage_detT_2nd-3}\\
        &\quad \displaystyle \sum_{j\in J}y_{ij} \leq x_{i,T}, & \forall \;i\in I \label{eq:twostage_detT_2nd-4}\\
        &\quad \displaystyle \sum_{i\in I}y_{ij} \leq d^{\xi_T}_{j}, & \forall \;j\in J \label{eq:twostage_detT_2nd-5}\\
        &\quad \displaystyle y_{ij} \geq 0, & \forall \;i\in I, j\in J. \label{eq:twostage_detT_2nd-6}
    \end{align}
\end{subequations}
We care to note that, other stage-aggregation approaches with less restrictive levels have also been proposed in the literature. For instance, in~\cite{zou2018partially}, the authors introduce a \textit{partially} adaptive model, in which the decisions are fully adaptive up to a certain period, and then follow a two-stage approach thereafter. 

While the \ac{s2ssp} model is easier to solve than the \ac{famspd} model, it is not difficult to see how the static nature of the \ac{s2ssp} could compromise the quality of the resulting policy. The policy obtained by the \ac{s2ssp} prescribes prepositioning specific levels of relief items, for each period $t=1, \dots, T$, before seeing how the stochastic process will unfold. As such, if the stochastic process $(\xi_1, \dots, \xi_T)$ end up evolving in a manner that is drastically different than the one initially predicted by the conditional probability $\mathbf{P}^{T-1}_{\xi_{T} \mid \xi_1}$, a large penalty might be incurred (or large quantities of relief items might be salvaged). To allow for some adaptability in the decision-making, we consider an \ac{rh} procedure, which we discuss next. Note that we shall refer to the \ac{s2ssp} model introduced in this section, hereafter, as the \ac{s2sspd}.

\paragraph{Rolling-horizon (\ac{rh}) procedure.} The \ac{rh} procedure injects adaptability directly into the policy by solving a new optimization problem (usually referred to as a look-ahead model) in every period $t \in \{1, \dots, T\}$. An optimal policy for the \ac{msp} formulation~\eqref{eq:FOSDDP} is induced by the expected cost-to-go functions $\Q^{\xi_t}_{t+1}(\cdot), \forall \; \xi_t \in \Xi_t,\, t=1, \dots,T$, which can be approximated by a collection of cutting planes, in an \textit{offline} fashion (see Algorithm~\ref{alg:nestedBenders_MC} in the Appendix). Instead, in the \ac{rh} procedure, the policy is obtained in an \textit{online} fashion, by solving a look-ahead model in each period of the planning horizon. In its simplest form, the \ac{rh} proceeds as follows. In every period $t=1, \dots, T$, the \ac{dm}: (i) defines an optimization problem with $\tau \in \{1,2\ldots, T-t+1\}$ stages to serve as the look-ahead model that approximates the fully adaptive $(T-t+1)$-stage problem at period $t$; (ii) solves this $\tau$-stage problem; (iii) implements \emph{only} the decisions corresponding to the same period $t$; and (iv) rolls forward one period. 

Several methods have been proposed in the literature on how to define the $\tau$-stage problem that serves as the look-ahead model in every roll (see e.g., \cite{alden1992rolling,chand2002forecast,sethi1991theory,siddig2021rolling,zou2018partially}). In this paper, we consider that $\tau=2$, i.e., we solve an \ac{s2sspd} model of the form~\eqref{eq:twostage_detT_1st} and~\eqref{eq:twostage_detT_2nd} in every roll. We refer to this approach that implements the \ac{s2sspd} in the \ac{rh} procedure as the \ac{rh2sspd}.

The advantage of using this framework is twofold. First, from a solution quality perspective, in each period $t$, by implementing the decisions pertaining to the current period $t$ only, the policy is naturally adaptive to the newly observed realizations of $\xi_t$, as opposed to the \ac{s2sspd} policy which commits to all decisions obtained from the model solved at the first period. Second, from a computational perspective, in each roll $t=1,\dots, T-1$, the \ac{dm} solves a \ac{2ssp} model which is easier to solve than the \ac{famspd} model. We investigate these claims in the context of the \ac{hdrlp} problem in our numerical study in Section~\ref{sec:case_study}. 

\subsection{Multi-stage Stochastic Programming Model with Random Landfall Time}
\label{subsec:randomT_MSP}
One crucial assumption of solving the \ac{msp} formulation~\eqref{eq:MSP} is that it has a fixed number of stages $T$, which is known a priori by the \ac{dm}. In the context of the \ac{hdrlp} problem, and in light of Assumption~\ref{assm:demand}, the value of $T$ corresponds to the landfall time of the hurricane. In practice, the time of landfall depends on several factors such as the hurricane's forward speed and trajectory, which are random over time. For example, a rapid surge in the hurricane's forward speed could cause the hurricane to arrive much earlier than anticipated -- rendering the \ac{dm} unprepared for its early arrival. There are other situations where the hurricane encounters circumstances hostile to its flourishing (e.g., cooler sea surface temperatures, shearing upper-level winds) which make it weaken or even dissipate altogether~\cite{bosart2000environmental}. In this subsection, we lift the assumption that the landfall time $T$ is fixed and known in advance. To that end, we extend the \ac{famspd}, \ac{s2sspd}, and \ac{rh2sspd} to the situation where $T$ is \textit{random}. We do this, first, by proposing an \ac{msp} model with a random number of stages in the spirit of the one presented in~\cite{guigues2021multistage}, which we briefly review next. 

~\cite{guigues2021multistage} presents an \ac{msp} model with a random number of stages $T$, where $\xi_t$'s are assumed to be stage-wise independent. Note that with the stage-wise independence assumption, when $T$ is fixed and known a priori, it is only necessary to define a single expected cost-to-go function $\Q_t(\cdot)$ for every stage in the planning horizon. To address the randomness of $T$, a parameter $T_{\max}<\infty$, referred to as the \textit{maximum} possible number of stages, is introduced to replace $T$ in the \ac{msp} formulation~\eqref{eq:MSP}. To facilitate \ac{dpe} analogous to~\eqref{eq:FOSDDP}, the following two features are introduced. First, a new state variable $\mathcal{G}_t = \mathbbm{1}_{\{t<T\}}$ is introduced in each stage $t = 1,\cdots, T_{\max}$, which indicates whether or not the terminal stage has occurred. Second, the expected cost-to-go function $\Q_t(a_{t-1})$ is replaced by $\Q_t(a_{t-1},\mathcal{G}_t)$. In other words, instead of using one expected cost-to-go function for each stage, the expected cost-to-go function is parameterized by whether or not the terminal stage has occurred yet, such that:
\begin{itemize}
    \item $\Q_t(a_{t-1}, 1)$ is the expected cost-to-go function when the terminal stage is yet to occur.
    \item $\Q_t(a_{t-1}, 0)$ is the expected cost-to-go function when the terminal stage has already occurred. This is also referred to as a \textit{null} function. 
\end{itemize}
Letting $p_t=\mathbb{P}(t\geq T)$, the idea then is to rewrite the \ac{dpe} as:
\begin{equation}
\label{eq:randomTDP}
\displaystyle Q_t(a_{t-1},\xi_t,\mathcal{G}_{t-1},\mathcal{G}_t) := \displaystyle \underset{{a_t}}{\min}\left\{\mathcal{G}_{t-1} z_t({a_t},\xi_t) + \Q_{t+1}(a_t,\mathcal{G}_t) \; \middle\vert \; a_t \in \A_t(a_{t-1},\xi_{t})\right\}\quad \forall \; t=1,\dots, T_{\max},
\end{equation}
where 
\begin{equation}
\label{eq:randomTDP_costtogo_t}
\displaystyle \Q_{t+1}(a_{t}, 1) = (1-p_{t+1})\mathbb{E}[Q_{t+1}(a_{t},\xi_{t+1},1,1)] + p_{t+1} \mathbb{E}[Q_{t+1}(a_{t},\xi_{t+1},1,0)],
\end{equation}
$\Q_{t+1}(a_{t}, 0) := 0$, for $t = 1,\cdots, T_{\max}-1$ and $\Q_{T_{\max}+1}(a_{T_{\max}}, \mathcal{G}_{T_{\max}}) := 0$. 

Note again that, since in~\cite{guigues2021multistage} the stochastic process is assumed to be stage-wise independent, the $\Q_{t}(\cdot)$ functions do not have the superscripts $\xi_{t-1}$. We extend formulation~\eqref{eq:randomTDP} to the case where the $\xi_t$'s are assumed to be stage-wise \textit{dependent}, following an \ac{mc} model.

\paragraph{\ac{dpe} for an \ac{msp} model with a random $T$ and stage-wise dependent $\xi_t$.} The key idea in our proposed model is to interpret the function $\Q_t(a_{t-1}, 1)$ in \eqref{eq:randomTDP_costtogo_t} as the expected cost-to-go function corresponding to a \textit{transient} state of the \ac{mc}: upon reaching this transient state, an immediate cost $z_t({a_t},\xi_t)$ is incurred, and the process evolves into a new state in the next period. By contrast, the null function $\Q_t(a_{t-1}, 0)$ can be interpreted as the expected cost-to-go function corresponding to an \textit{absorbing} state in the \ac{mc}: once the stochastic process reaches this absorbing state, it will remain there without incurring any future cost. Given this interpretation, we partition the state space into two sets $\mathcal{T}$ and $\mathcal{A}$, where $\mathcal{T}$ represents the set of transient states and $\mathcal{A}$ represents the set of absorbing states, such that $\Xi_t = \mathcal{T} \cup \mathcal{A}$. We then rewrite the \ac{dpe}~\eqref{eq:FOSDDP} as follows.
\begin{itemize}
    \item For a \textit{transient} state $\xi_t \in \mathcal{T}$:
    \begin{equation}
    \label{eq:FOSDDP_randomT_T}
        Q_t({a_{t-1}},\xi_{t}):= \displaystyle \underset{{a_t}}{\min}\left\{z_t({a_t},\xi_t) + \displaystyle \Q^{\xi_t}_{t+1}(a_t) \; \middle\vert \; a_t \in \A_t(a_{t-1},\xi_{t})\right\}\quad \forall \; t=1,\dots, T_{\max},
    \end{equation}
    where $\displaystyle \Q^{\xi_t}_{t+1}(a_t) := \sum_{\xi_{t+1}\in \Xi_{t+1}} \mathbf{P}_{\xi_{t+1} \mid \xi_t}  Q_{t+1}(a_t,\xi_{t+1}), \forall \;t = 1, \dots, T_{\max}-1,$ and $\Q^{\xi_{T_{\max}}}_{T_{\max}+1}(a_{T_{\max}}) := 0$.
    \item For an \textit{absorbing} state $\xi_t \in \mathcal{A}$:
    \begin{equation}
    \label{eq:FOSDDP_randomT_A}
       Q_t({a_{t-1}},\xi_{t}):= \displaystyle \underset{{a_t}}{\min}\left\{g_t({a_t},\xi_t) + \displaystyle \Q^{\xi_t}_{t+1}(a_t) \; \middle\vert \; a_t \in \A_t(a_{t-1},\xi_{t})\right\}\quad \forall \; t=1,\dots, T_{\max},
    \end{equation}
    where $g_t({a_t},\xi_t) := 0$ and $\Q^{\xi_t}_{t+1}(a_t) := 0, \forall \;t = 1, \dots, T_{\max}-1$, and $\Q^{\xi_{T_{\max}}}_{T_{\max}+1}(a_{T_{\max}}) := 0$.
\end{itemize}
Note that, while we set $g_t({a_t},\xi_t) = 0, \forall \;t\leq T_{\max}-1$ based on our assumptions for the hurricane relief logistics problem (Assumption 3), writing~\eqref{eq:FOSDDP_randomT_A} in this form allows for more general model, for example, when the absorbing states correspond to a cost that is paid indefinitely. 

We care to note that, with slight modification, the nested Benders decomposition, shown in Algorithm~\ref{alg:nestedBenders_MC} in the Appendix, can be applied to solve the \ac{dpe}~\eqref{eq:FOSDDP_randomT_T} and~\eqref{eq:FOSDDP_randomT_A}, and the convergence of the algorithm follows by the argument shown in~\cite{guigues2021multistage}. Nevertheless, we need to adapt our \ac{mc} model, which characterizes the hurricane evolution, to incorporate the newly introduced temporal dimension of the hurricane random time of landfall. Additionally, since we do not know in advance when the hurricane makes landfall, we need to adapt our definition of the decision variables, the objective functions, and the constraint sets in the \ac{msp} model~\eqref{eq:MSP} to reflect the same. We discuss these modifications next.

\begin{figure}[htbp]
\begin{center}
\begin{tikzpicture}[scale=1.25, transform shape]
\draw[step=0.5cm,gray,very thin] (-0.1,-3.1) grid (7.1,1.1);
\draw[->,line width=1.5pt] (0,0) -- (7.25,0);
\draw[->,line width=1.5pt] (0,0) -- (0,1.25);
\draw[->,line width=1.5pt] (0,0) -- (0,-3.25);
\draw[black]  (0,0) -- (7.00,0) -- (7.00,1.00) -- (0,1.00) -- (0,0);
\draw[black] (0,0) -- (0,-3) -- (7,-3) -- (7,0) -- (0,0);
\fill[line width=2pt,color=black,fill=black,fill opacity=0.1] (3.50,2.00) -- (4.00,2.50) -- (3.50,3.00) -- (3.00,2.50) -- cycle;
\fill[line width=2pt,color=black,fill=black,fill opacity=0.1] (1.20,1.80) -- (1.80,1.80) -- (1.80,1.20) -- (1.20,1.20) -- cycle;
\fill[line width=2pt,color=black,fill=black,fill opacity=0.1] (5.20,1.80) -- (5.80,1.80) -- (5.80,1.20) -- (5.20,1.20) -- cycle;

\draw [line width=1.25pt,color=black] (1.00,0.50) circle (0.4cm);
\draw [line width=1.25pt,color=black] (3.50,0.50) circle (0.4cm);
\draw [line width=1.25pt,color=black] (6.00,0.50) circle (0.4cm);
\draw [line width=1.25pt,color=black] (3.50,2.00)-- (4.00,2.50);
\draw [line width=1.25pt,color=black] (4.00,2.50)-- (3.50,3.00);
\draw [line width=1.25pt,color=black] (3.50,3.00)-- (3.00,2.50);
\draw [line width=1.25pt,color=black] (3.00,2.50)-- (3.50,2.00);
\draw [line width=1.25pt,color=black] (1.20,1.80)-- (1.80,1.80);
\draw [line width=1.25pt,color=black] (1.80,1.80)-- (1.80,1.20);
\draw [line width=1.25pt,color=black] (1.80,1.20)-- (1.20,1.20);
\draw [line width=1.25pt,color=black] (1.20,1.20)-- (1.20,1.80);
\draw [line width=1.25pt,color=black] (5.20,1.80)-- (5.80,1.80);
\draw [line width=1.25pt,color=black] (5.80,1.80)-- (5.80,1.20);
\draw [line width=1.25pt,color=black] (5.80,1.20)-- (5.20,1.20);
\draw [line width=1.25pt,color=black] (5.20,1.20)-- (5.20,1.80);
\draw [->,line width=1.25pt,color=black] (3.50,2.00) -- (1.80,1.80);
\draw [->,line width=1.25pt,color=black] (3.50,2.00) -- (5.20,1.80);
\draw [->,line width=1.25pt,color=black] (1.80,1.60) -- (5.20,1.6);
\draw [->,line width=1.25pt,color=black] (5.20,1.40) -- (1.80,1.4);
\draw [->,line width=1.25pt,color=black] (1.80,1.20) -- (3.1,0.5);
\draw [->,line width=1.25pt,color=black] (1.80,1.20) -- (5.60,0.6);
\draw [->,line width=1.25pt,color=black] (5.20,1.20) -- (5.625,0.7);
\draw [->,line width=1.25pt,color=black] (5.20,1.20) -- (3.9,0.5);
\draw [->,line width=1.25pt,color=black] (5.20,1.20) -- (1.4,0.6);
\draw [->,line width=1.25pt,color=black] (1.80,1.20) -- (1.375,0.7);
\begin{scriptsize}
\draw [fill=black] (0,0) circle (1pt);
\draw [fill=black] (1.00,0) circle (1pt);
\draw [fill=black] (2.00,0) circle (1pt);
\draw [fill=black] (3.00,0) circle (1pt);
\draw [fill=black] (4.00,0) circle (1pt);
\draw [fill=black] (5.00,0) circle (1pt);
\draw [fill=black] (6.00,0) circle (1pt);
\draw [fill=black] (7.00,0) circle (1pt);

\draw[color=black] (3.5,2.5) node {MDC};

\draw[color=black] (1.5,1.5) node {SP1};

\draw[color=black] (5.5,1.5) node {SP2};
\draw[color=black] (1,0.5) node {DP1};
\draw[color=black] (3.5,0.5) node {DP2};
\draw[color=black] (6,0.5) node {DP3};
\end{scriptsize}
\begin{tiny}

\draw[color=black] (1,-0.2) node {$100$};
\draw[color=black] (2.15,-0.2) node {$200$};
\draw[color=black] (3.125,-0.2) node {$300$};
\draw[color=black] (4.1,-0.2) node {$400$};
\draw[color=black] (5,-0.2) node {$500$};
\draw[color=black] (6,-0.2) node {$600$};
\draw[color=black] (7,-0.2) node {$700$};
\draw[color=black] (-0.3,0.5) node {$50$};
\draw[color=black] (-0.3,1) node {$100$};
\draw[color=black] (-0.4,-0.5) node {$-50$};
\draw[color=black] (-0.4,-1) node {$-100$};
\draw[color=black] (-0.4,-1.5) node {$-150$};
\draw[color=black] (-0.4,-2) node {$-200$};
\draw[color=black] (-0.4,-2.5) node {$-250$};
\draw[color=black] (-0.4,-3) node {$-300$};

\draw[color=black] (7.5,0) node {$\ell_{x,t}$};
\draw[color=black] (0,1.4) node {$\ell_{y,t}$};

\end{tiny}
\draw [fill=black] (0,-3) circle (2pt);
\draw [fill=black] (1,-2.5) circle (2pt);
\draw [-latex,line width=1.25pt,color=black] (0,-3) -- (1,-2.5);
\draw [fill=black] (1.5,-2) circle (2pt);
\draw [-latex,line width=1.25pt,color=black] (1,-2.5) -- (1.5,-2);
\draw [fill=black] (2,-1.5) circle (2pt);
\draw [-latex,line width=1.25pt,color=black] (1.5,-2) -- (2,-1.5);
\draw [fill=black] (3,-1) circle (2pt);
\draw [-latex,line width=1.25pt,color=black] (2,-1.5) -- (3,-1);
\draw [fill=black] (5,-0.5) circle (2pt);
\draw [-latex,line width=1.25pt,color=black] (3,-1) -- (5,-0.5);
\draw [fill=black] (6,0) circle (2pt);
\draw [-latex,line width=1.25pt,color=black] (5,-0.5) -- (6,0);
\draw [fill=black] (4,0) circle (2pt);
\draw [-latex,line width=1.6pt,color=black,dashed] (3,-1) -- (4,0);
\draw [fill=black] (3,0) circle (2pt);
\draw [-latex,line width=1.6pt,color=black,dashed] (2,-1.5) -- (3,0);
\end{tikzpicture}    
\end{center}
\caption{The hurricane trajectory towards the affected potential area with \textit{random} landfall time.}
\label{fig:network_randomT}
\end{figure}
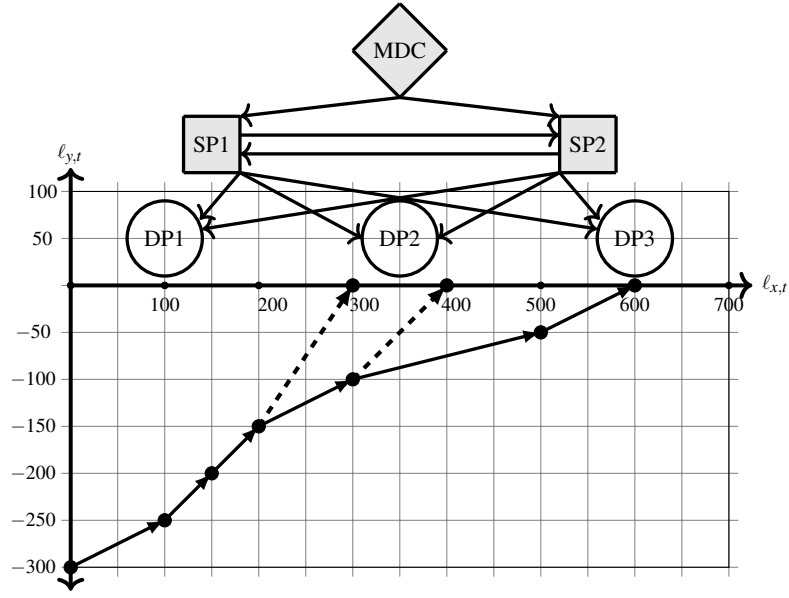

\paragraph{Incorporating the temporal component to the \ac{mc} model for hurricane evolution.}
We first split the random variable $\ell_t$ denoting the hurricane location at time $t$ into two independent random variables $\ell_{x,t} \in L_{x}$ and $\ell_{y,t} \in L_{y}$, representing the $x$-coordinate and $y$-coordinate of the hurricane location at time $t$, respectively. We further assume that the evolution of each $\ell_{x,t}$ and $\ell_{y,t}$ is given by an \ac{mc} model, where the one-step transition probability matrix is given by $\mathbf{P}_{\ell_{x,t} \mid \ell_{x,t-1}}$ and $\mathbf{P}_{\ell_{y,t} \mid \ell_{y,t-1}}$, respectively. For illustration, we extend the $y$-coordinate in the example shown in Figure~\ref{fig:network} to include $\ell_{y,t}<0$ values, discretize it, and let each point in the discretized space represent the location of the hurricane at time $t = 1,\cdots, T_{\max}$ (see Figure~\ref{fig:network_randomT}), such that:
\begin{itemize}
    \item $\ell_{y,t} < 0$ means that the hurricane has not made landfall yet at time $t$.
    \item $\ell_{y,t} = 0$ means that the hurricane makes landfall at time $t$.
    \item $\ell_{y,t} > 0$ means that the hurricane has already made landfall before time $t$.
\end{itemize}
We now define the state space of the modified \ac{mc} model to be $\Xi_t = \mathrm{A} \times \mathrm{L}_{x}\times \mathrm{L}_{y}$, and let $\xi_t = (\alpha_t, \ell_{x,t}, \ell_{y,t})$ denote the Markovian state at time $t$, where $\ell_{y,t} \leq 0 \implies \xi_t \in \mathcal{T}$, i.e., $\xi_t$ is a \textit{transient} state. By contrast, $\ell_{y,t} > 0 \implies \xi_t \in \mathcal{A}$, i.e., $\xi_t$ is an \textit{absorbing} state. Similarly to the \ac{mc} model with a deterministic time of landfall, we can now redefine the \textit{joint} probability distribution as $\mathbf{P}_{\xi_{t+1} \mid \xi_t} := \mathbf{P}_{\alpha_{t+1} \mid \alpha_t} \times \mathbf{P}_{\ell_{x,t+1} \mid \ell_{x,t}} \times \mathbf{P}_{\ell_{y,t+1} \mid \ell_{y,t}}$.

We care to note that, in addition to making the time of landfall random, this construction can also allow for modeling other situations. For instance, one can define an absorbing state corresponding to the situation where the hurricane dissipates. This, for instance, could be because the wind slows down and the storm subsides (i.e., the intensity $\alpha_t = 0$) or because the hurricane trajectory drifts outside of the \ac{apa}. Figure~\ref{fig:network_randomT} illustrates the \ac{mc} model with a random landfall time, by showing different hurricane trajectories towards the \ac{apa}. 

The last missing piece for formulating this model is to define an appropriate value for $T_{\max}$. This is, of course, case-by-case, and we discuss how to choose this value in Section~\ref{sec:problem_data} in the Appendix.

\subsubsection{FA-MSP model with a random time of landfall}
\label{subsubsec:randomT_MSP_stochastic_FA}
In this subsection, we present the \ac{famspr}. Recall that by Assumption 3, the stage-$t$ problems associated with the absorbing states are degenerate since, according to~\eqref{eq:FOSDDP_randomT_A}, the immediate cost functions are defined as $g_t({a_t},\xi_t) = 0$, and the expected cost-to-go functions $\Q^{\xi_t}_{t+1}(a_t) = 0, \; \forall \; \xi_t \in \mathcal{A}, \; t = 1, \dots, T_{\max}-1$. As such, we focus on the definition of the stage-$t$ problems associated with transient states only.

Since we assume the time of landfall $T$ to be random, the demand can happen at any point before $T_{\max}$. Hence, we need to index the demand quantities $d$ by the period $t$ for $t = 1,\cdots, T_{\max}$, such that $d_t := D(\alpha_t, \ell_{x,t})$, when $t=T$ and $d_t := 0, \forall \; t\neq T$. Additionally, since it is not possible to decouple the stage-$t$ problems, as we did in the \ac{famspd} into a set of non-terminal stages problems and a terminal-stage problem, all of the state and local variables need to be included in every stage-$t$ problem, $\forall \;\; t = 1,\cdots,T_{\max}$. To simplify the presentation, we introduce two auxiliary decision variables, $\underline{e}_t$ and $\overline{e}_t$ for $t = 1,\cdots, T_{\max}$, such that:
\begin{itemize}
    \item $\underline{e}_t := (\underline e_{j,t})_{j \in J}$ denotes the amount of unsatisfied demand (if any).
    \item $\overline{e}_t:= (\overline e_{i,t})_{i \in I}$ denotes the amount of salvaged items (if any).
\end{itemize}
Note that, when $d_{j,t}$ is known, we can adopt the deterministic model~\eqref{eq:det_form} as follows. 
\begin{subequations}
    \label{eq:det_form_randT}
    \begin{align}
    && \llap{$\underset{x, f, y, \underline{e}, \overline{e}}{\min} \quad \displaystyle \sum_{t=1}^{T_{\max}} \left(\sum_{i\in \{0\}\cup I}\sum_{i'\in I}c^b_{ii',t}f_{ii',t}+\sum_{i\in I}c^h_{i,t}x_{i,t}+ h_t\sum_{i\in I}f_{0i,t} +  \sum_{i\in I}\sum_{j\in J}c^a_{ij,t}y_{ij,t} + p\sum_{j \in J} \underline e_{j,t} + q\sum_{i \in I} \overline e_{i,t} \right)$} &   \nonumber \\    
    \text{s.t.} \quad & \displaystyle x_{i,t-1} + \sum_{i'\in \{0\}\cup I, i'\neq i}f_{i'i,t} - \sum_{i'\in I, i'\neq i}f_{ii',t} - \sum_{j\in J}y_{ij,t} - \overline e_{i,t} = x_{i,t}, & \forall \;i\in I, \; t = 1,\ldots, T_{\max}, \label{eq:det_form_randT-1} \\
    \quad & \displaystyle \sum_{i'\in I, i'\neq i}f_{ii',t} \leq x_{i,t-1}, & \forall \;i\in I, \; t = 1,\ldots, T_{\max}, \label{eq:det_form_randT-2} \\
    \quad & \displaystyle x_{i,t} \leq u_i, & \forall \;i\in I, \; t = 1,\ldots, T_{\max}, \label{eq:det_form_randT-3} \\
    \quad & \displaystyle \sum_{i \in I}y_{ij,t} + \underline e_{j,t}\geq d_{j,t}, & \forall \;j\in J,\; t = 1,\ldots, T_{\max}, \label{eq:det_form_randT-4} \\
    \quad & \displaystyle x_t, f_t, y_t, \underline{e}_t, \overline{e}_t \; {\geq 0}, & \forall \;t = 1,\ldots, T_{\max}. \label{eq:det_form_randT-5} 
   \end{align}
\end{subequations}
Letting $a_{t} :=(x_t, f_t, y_t, \underline{e}_t, \overline{e}_t)$, we can juxtapose the remaining components of the modified deterministic formulation~\eqref{eq:det_form_randT} with the components of formulation~\eqref{eq:FOSDDP_randomT_T} $\forall \;t = 1, \cdots, T_{\max}$, as follows. 
\begin{equation}
\label{eq:objective_randomT}
z_t({a_t},\xi_t) := \displaystyle \sum_{i\in \{0\}\cup I}\sum_{i'\in I}c^b_{ii',t}f_{ii',t}+\sum_{i\in I}c^h_{i,t}x_{i,t}+ h_t\sum_{i\in I}f_{0i,t} +  \sum_{i\in I}\sum_{j\in J}c^a_{ij,t}y_{ij,t} + p\sum_{j \in J} \underline e_{j,t} + q\sum_{i \in I} \overline e_{i,t},
\end{equation}
and
\begin{equation}
\label{eq:feasible_set_randomT}
\A_t({a_{t-1}},\xi_{t}) :=
    \begin{cases} 
    \displaystyle x_{i,t-1} + \sum_{i'\in \{0\}\cup I, i'\neq i}f_{i'i,t} - \sum_{i'\in I, i'\neq i}f_{ii',t} - \sum_{j\in J}y_{ij,t} - \overline e_{i,t} = x_{i,t}, & \forall \;i\in I, \\
    \displaystyle \sum_{i'\in I, i'\neq i}f_{ii',t} \leq x_{i,t-1}, & \forall \;i\in I, \\
    \displaystyle x_{i,t} \leq u_i, & \forall \;i\in I,\\
    \displaystyle \sum_{i \in I}y_{ij,t} + \underline e_{j,t}\geq d^{\xi_t}_{j,t}, & \forall \;j\in J, \\
    \displaystyle x_t, f_t, y_t, \underline{e}_t, \overline{e}_t \; {\geq 0}. &
    \end{cases}
\end{equation}

\subsubsection{Alternative decision policies for cases with random landfall time}
\label{subsubsec:alternative_approx_policies_Trand}
We also consider alternative decision policies analogous to the  \ac{s2sspd} and \ac{rh2sspd} discussed in Subsection~\ref{subsubsec:alternative_approx_policies_T}.

\paragraph{Static2SSP model with random landfall time $T$.}
The \ac{s2sspr} presented in this section is similar to the \ac{s2sspd} presented in Subsection~\ref{subsec:determinsticT_MSP}: they both introduce a formulation where different periods across the planning horizon are aggregated into two stages only and they both render a static policy which entails making all the prepositioning decisions, in advance, before the hurricane making landfall. Specifically, in the \ac{s2sspd}, the certainty about the hurricane's time of landfall facilitates for decoupling the \ac{famspd} problem~\eqref{eq:FOSDDP} into a first-stage problem~\eqref{eq:twostage_detT_1st} that includes all periods before the landfall time, where decisions concerning the preparation of the hurricane arrival are made; and a second-stage problem~\eqref{eq:twostage_detT_2nd} that corresponds to the landfall time, where the \ac{dm}'s reaction decisions are made once the uncertain demand is revealed. However, when $T$ is random, the static nature of a two-stage approximation has some delicate implications on defining the first-stage and second-stage problems. 

The main concern is that, when $T$ is random, if the prepositioning decisions that are previously made cannot be adapted, it is possible that the \ac{dm} might arrive at a static decision which prescribes prepositioning relief items after the hurricane's landfall. That is, it could happen that $x_t > 0$ and $f_t> 0$ for some $T < t \leq T_{\max}$. To see how this could be case, let $t'$ be an arbitrary period in the planning horizon such that $1 < t' \leq T_{\max}$, and suppose that we solve an \ac{s2ssp} model with $T_{\max}$ stages, such that the resulting policy prescribes to have the following prepositioning decisions: $x_t > 0$ and $f_t > 0$ for some $t \leq t'$, and $x_t = 0$ and $f_t = 0, \; \forall \; t' < t$. In other words, the \ac{dm} starts (and continues) the prepositioning at any point in time before $t \leq t'$, and stops after time $t > t'$. Since $T$ is assumed to be random, there are two possible outcomes for the hurricane's time of landfall: (i) the hurricane makes landfall \textit{after/when} all the prepositioning is finished, i.e., ${t' \leq T}$. This situation is similar to the one discussed for the deterministic landfall time case; (ii) the hurricane makes landfall \textit{before} the prepositioning is finished, i.e., ${T < t'}$. This second situation has two further implications. First, as is the case for the model with a deterministic $T$, upon arriving at the landfall the \ac{dm} pays penalty for demand shortage (if any) or salvages the remaining items (if any). The second implication, which is of more relevance, is to start or continue prepositioning after the hurricane has made landfall. In light of Assumption~\ref{assm:demand}, the act of starting or continuing to preposition relief items despite knowing that the hurricane has already made landfall is somewhat pathological. To circumvent this potential shortcoming, we introduce a set of ``reimbursement'' parameters $r_t, \forall \; t = 1, \cdots, T_{\max}$ to the second-stage problem to compensate the \ac{dm} for all the logistics costs that were counted in the first stage for every period $t > T$, such that
\begin{equation}
\label{eq:reimbursement_var}
    r_t :=
    \begin{cases} 
        \displaystyle -\left(\sum_{i\in \{0\}\cup I}\sum_{i'\in I}c^b_{ii',t}{f}_{ii',t}+\sum_{i\in I}c^h_{i,t}{x}_{i,t}+ h_t\sum_{i\in I}{f}_{0i,t}\right) & \forall \;\; t > T, \\
        0                                                                                                                                                   & \text{otherwise}.
    \end{cases}
\end{equation}

Putting everything together, given initial data $x_0$ and $\xi_1$, we define the first-stage and second-stage problems in the \ac{s2sspr} model as follows.

The \textit{first-stage} problem:
\begin{subequations}
    \label{eq:2SSP_rand_1st-stage}
    \begin{align}
    \displaystyle \underset{x, f}{\min}\quad & \displaystyle \sum_{t=1}^{T_{\max}} \left(\sum_{i\in \{0\}\cup I}\sum_{i'\in I}c^b_{ii',t}f_{ii',t}+\sum_{i\in I}c^h_{i,t}x_{i,t}+ h_t\sum_{i\in I}f_{0i,t}\right) + \Q^{\xi_1}(x, f) & \nonumber \\ 
    \text{s.t.} 
     & \displaystyle \sum_{j\in J, j\neq i}f_{ij,t} \leq x_{i,t-1}, & \forall \;i\in I,\; t = 1,\cdots, T_{\max}\label{eq:2SSP_rand_1st-stage-1} \\
     & \displaystyle x_{i,t} \leq u_i, & \forall \;i\in I,\; t = 1,\cdots, T_{\max} \label{eq:2SSP_rand_1st-stage-2} \\
     & \displaystyle x_t, f_t\geq 0, & \forall \;t = 1,\cdots, T_{\max} \label{eq:2SSP_rand_1st-stage-3},
   \end{align}
\end{subequations}
where $x = (x_1, \dots, x_{T-1})$, $f = (f_1, \dots, f_{T-1})$, and $\Q^{\xi_1}(x, f) := \sum_{\xi_{T_{\max}} \in \Xi_{T_{\max}}} \mathbf{P}^{T_{\max}-1}_{\xi_{T_{\max}} \mid \xi_1} Q(x, f, \xi_{[1:T_{\max}]})$. Note that, since we do not know in advance when the hurricane makes landfall, the second-stage value function $Q(\cdot)$ is now a function of $\xi_{[1:T_{\max}]}$ which gives the entire history of the stochastic process from $t=1$ to $t=T_{\max}$.

To define the second-stage problem, we let $\mathbbm{1}_{\{t>T^{\xi_{[1:T_{\max}]}}\}}$ be an indicator function on whether or not period $t$ is after the hurricane's landfall time $T^{\xi_{[1:T_{\max}]}}$ associated with trajectory $\xi_{[1:T_{\max}]}$.

The \textit{second-stage} problem:
\begin{subequations}
    \label{eq:2SSP_rand_2nd-stage}
    \begin{align}
    \displaystyle Q(x, f, \xi_{[1:T_{\max}]}) := \displaystyle \underset{y, \underline{e}, \overline{e}}{\min} \quad & \displaystyle \sum_{t=1}^{T_{\max}} \left(\sum_{i\in I}\sum_{j\in J}c^a_{ij,t}y_{ij,t} + p\sum_{j \in J} \underline e_{j,t} + q\sum_{i \in I} \overline e_{i,t}\right) & \nonumber \\
    - &\displaystyle \sum_{t=1}^{T_{\max}}\left(\sum_{i\in \{0\}\cup I}\sum_{i'\in I}c^b_{ii',t}{f}_{ii',t}+\sum_{i\in I}c^h_{i,t}{x}_{i,t}+ h_t\sum_{i\in I}{f}_{0i,t}\right)\mathbbm{1}_{\{t>T^{\xi_{[1:T_{\max}]}}\}} & \nonumber \\
    \text{s.t.} \quad  & \displaystyle \sum_{j\in J}y_{ij,t} + \overline e_{i,t} =  x_{i,t-1} + \sum_{j\in \{0\}\cup I, j\neq i} f_{ji,t} - \sum_{j\in I, j\neq i} f_{ij,t} -  x_{i,t}, & \forall \;i\in I,\; t=1, 2, \cdots, T_{\max} \label{eq:2SSP_rand_2nd-stage-1} \\
      & \displaystyle \sum_{i \in I}y_{ij,t} + \underline e_{j,t}\geq d^{\xi_t}_{j,t}, & \forall \;j\in J,\; t=1, 2, \cdots, T_{\max} \label{eq:2SSP_rand_2nd-stage-2} \\
      & \displaystyle y_t, \underline{e}_t , \overline{e}_t \geq 0, &\; \forall \;  t=1, 2, \cdots, T_{\max}. \label{eq:2SSP_rand_2nd-stage-3} 
   \end{align}
\end{subequations}

\begin{figure}[htbp]
\begin{center}
\begin{tikzpicture}[scale=1.15, transform shape]

\draw[step=1cm,black,very thin] (0.75,-0.1) grid (9.25,0.1);

\draw[->, line width=1.5pt, color=black] (0.75,0) -- (9.5,0);

\begin{scriptsize}
\draw[color=black] (1,-0.3) node {$1$};
\draw[color=black] (2,-0.3) node {$t-1$};
\draw[color=black] (3,-0.3) node {$t$};
\draw[color=black] (4,-0.3) node {$\dots$};
\draw[color=black] (5,-0.3) node {$T$};
\draw[color=black] (6,-0.3) node {$\dots$};
\draw[color=black] (7,-0.3) node {$t'-1$};
\draw[color=black] (8,-0.3) node {$t'$};
\draw[color=black] (9,-0.3) node {$T_{\max}$};
\end{scriptsize}

\exclude{
\draw [color=black, line width=2pt] (0.5,1) circle (3.5pt); 
\draw[line width=2pt, color=black] (0.5,0.9) -- (0.5,0.4);  
\draw[line width=2pt, color=black] (0.5,0.7) -- (0.25,0.5); 
\draw[line width=2pt, color=black] (0.5,0.7) -- (0.75,0.5); 
\draw[line width=2pt, color=black] (0.5,0.4) -- (0.3,0.2); 
\draw[line width=2pt, color=black] (0.5,0.4) -- (0.7,0.2); 
\draw[color=black] (0.5,1.3) node {\scriptsize DM}; 
}

\draw [line width=2pt, decorate,decoration={brace,amplitude=10pt}] (1,0.25) -- (8,0.25);
\draw[color=black] (4.5,1.7) node {\scriptsize \textbf{1st-stage:} prepositioning the relief items};
\draw[color=black] (4.5,1) node {\scriptsize $\displaystyle \sum_{t=1}^{t'} \left(\sum_{i\in \{0\}\cup I}\sum_{i'\in I}c^b_{ii',t}f_{ii',t}+\sum_{i\in I}c^h_{i,t}x_{i,t}+ h_t\sum_{i\in I}f_{0i,t}\right)$};

\draw [-latex,line width=0.75pt,color=black] (5,0) -- (3,-1);
\draw [fill=black] (5,0) circle (2.5pt);
\draw[color=black] (0.5,-1) node {\scriptsize \textbf{2nd-stage-(i):} hurricane landfall};
\draw[color=black] (1,-1.75) node {\scriptsize $\displaystyle \left(\sum_{i\in I}\sum_{j\in J}c^a_{ij,T}y_{ij,T} + p\sum_{j \in J} \underline e_{j,T} + q\sum_{i \in I} \overline e_{i,T}\right)$};

\draw [line width=2pt, decorate,decoration={brace,amplitude=10pt,mirror}] (5,-0.45) -- (8,-0.45);
\draw[color=black] (6.5,-1) node {\scriptsize \textbf{2nd-stage-(ii):} reimbursement};
\draw[color=black] (8,-1.75) node {\scriptsize $\quad\displaystyle -\sum_{t=T+1}^{t'}\left(\sum_{i\in \{0\}\cup I}\sum_{i'\in I}c^b_{ii',t}{f}_{ii',t}+\sum_{i\in I}c^h_{i,t}{x}_{i,t}+ h_t\sum_{i\in I}{f}_{0i,t}\right)$};
\end{tikzpicture}    
\end{center}
\caption{Symbolic representation of the \ac{s2sspr} with reimbursement.}
\label{fig:reimbursement}
\end{figure}
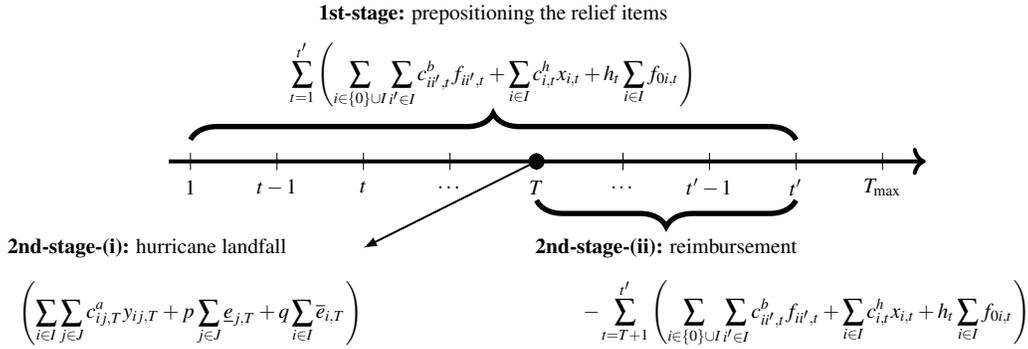

Figure~\ref{fig:reimbursement} shows a symbolic representation of the proposed \ac{s2sspr} formulation. 

\paragraph{\ac{rh2ssp} approach with random landfall time $T$.} In general, the \ac{rh2sspr} presented here is similar to the \ac{rh2sspd} discussed in Subsection~\ref{subsec:determinsticT_MSP} for the deterministic landfall time case: in each period $t$, the \ac{dm} solves an \ac{s2ssp} model that serves as the look-ahead model, implements the decisions corresponding to period $t$ only, and then rolls forward. Nonetheless, there is one, albeit small, remark we care to highlight for the sake of completeness, which we present in Section~\ref{sec:rh-2ssp-random} in the Appendix. 

\FloatBarrier 
\section{Numerical Results}
\label{sec:case_study}
In this section, we summarize the performances of different approaches presented in Section~\ref{sec:MSP_models}. In our experiments, we implemented eight different approaches. The first four correspond to the case when the hurricane's landfall time $T$ is deterministic: (i) \ac{cvd} described in~\eqref{eq:det_form}; (ii) the \ac{famspd} described in~\eqref{eq:FOSDDP},~\eqref{eq:objective_det_Tminus},~\eqref{eq:feasible_set_det_nonT},~\eqref{eq:objective_det_T}, and~\eqref{eq:feasible_set_det_T}; (iii) the \ac{s2sspd} described in~\eqref{eq:twostage_detT_1st} and~\eqref{eq:twostage_detT_2nd}; and (iv) the \ac{rh2sspd} described in Section~\ref{subsubsec:alternative_approx_policies_T}. The other four correspond to the case when the hurricane's landfall time $T$ is random: (v) the \ac{cvr} described in~~\eqref{eq:det_form_randT}; (vi) the \ac{famspr} described in~\eqref{eq:FOSDDP_randomT_T},~\eqref{eq:objective_randomT} and~\eqref{eq:feasible_set_randomT}; (vii) the \ac{s2sspr} described in~\eqref{eq:2SSP_rand_1st-stage} and~\eqref{eq:2SSP_rand_2nd-stage}; and (ix) the \ac{rh2sspr} discussed in Subsection~\ref{subsubsec:alternative_approx_policies_Trand} (and Section~\ref{sec:rh-2ssp-random} in the Appendix).

We benchmark the performances of the different policies using the following metrics.
\begin{itemize}
    \item $\hat{z}$: the sample mean for the performance of a policy calculated using the out-of-sample evaluation (see~\eqref{eq:sample_mean} in Appendix~\ref{sec:implement_details}).
    \item $\pm1.96 \hat{\sigma}/\sqrt{N}$: the width of the $95\%$-\ac{ci} on the mean performance. Here, $N$ is the number of scenarios used in the out-of-sample and $\hat{\sigma}$ is the sample standard deviation (see~\eqref{eq:sample_std} in Appendix~\ref{sec:implement_details}).
    \item \textit{gap}: the relative gap in the mean performance of each policy ($\hat{z}$) compared to the one obtained by the \ac{cv} solution. 
    \item \textit{time}: time (in seconds) required for training each policy. We use notation ``$-$'' to indicate that the time-limit of three hours is reached. Note that this statistics is only applicable to the \ac{famsp} and \ac{s2ssp} models. Moreover, the policy evaluation time for \ac{famsp} and \ac{s2ssp} was always less than 10 seconds and, hence, we do not report it.
\end{itemize}
We note that the reason for not reporting the time for training policies or performance evaluation for \ac{rh} approaches is due to the fact that the decision policy is obtained in an online fashion during the out-of-sample evaluation. It is not a fair comparison between the computational time associated with the \ac{rh} approaches and the time spent on training offline policies by the \ac{famsp} and \ac{s2ssp} models.

In all of our experiments, we use programming language \emph{Julia} version 1.6.1, utilizing \emph{JuMP} version 0.22.1 package~\cite{dunning2017jump}, with commercial solver \emph{Gurobi}, version 9.5.0 in our implementation. All of the numerical experiments are conducted on Clemson University's primary high-performance computing cluster, the \emph{Palmetto cluster}, where we used an \emph{R830 Dell Intel Xeon} ``big memory'' compute node with 2.60GHz, 1.0TB memory, and 24 cores. The problem data used in our implementation can be found in Appendix~\ref{sec:problem_data} and the remaining implementation details can be found in Appendix~\ref{sec:implement_details}. All of the codes can be found at Githb~\cite{Github}.

\subsection{Main results for the case of deterministic landfall time}
\label{subsubsec:main_results_deterministic}
In Table~\ref{tab:results_det}, we report the numerical results for \ac{famspd}, \ac{s2sspd} and \ac{rh2sspd}. The first three columns describe the test instances based on different values of $\nu$, $|I|$ and $|J|$, and columns 4-7, 8-10 and 11-14 show the performance metrics results for \ac{famspd}, \ac{rh2sspd} and \ac{s2sspd}, respectively. Additionally, Figure~\ref{fig:detT_allavg} depicts the $\hat{z}$ values of different policies, \textit{averaged} across all the instances for different $\nu$ values (see also Figure~\ref{fig:detT_avg} in the Appendix).
\exclude{

\begin{table}
\centering
\footnotesize
\begin{adjustbox}{width=\textwidth}
\begin{tabular}{@{}lllcccccccccccccc@{}}
\toprule
 & & \multicolumn{4}{c}{FA-MSP-D} & \phantom{abc}& \multicolumn{4}{c}{RH2SSP-D} & \phantom{abc} & \multicolumn{4}{c}{static2SSP-D}\\
\cmidrule{4-7} \cmidrule{9-12} \cmidrule{14-17}
  $\nu$ & $|I|\quad |J|$ && $\hat{z}$ & $\pm1.96 \hat{\sigma}/\sqrt{N}$ & time & gap && $\hat{z}$ & $\pm1.96 \hat{\sigma}/\sqrt{N}$ & time & gap && $\hat{z}$ & $\pm1.96 \hat{\sigma}/\sqrt{N}$ & time & gap\\ \midrule
  
0.001 & 3 $\quad$    10 && 2333.37 &180.31 &269.68 &200.80\%    && 2333.62 &180.38 &1624.59 &200.83\%  && 7015.46 &1143.26 &4.21 &804.38\%   \\ 
      &    $\quad\quad$ 20 && 4440.41 &314.50    &282.67 &186.57\%  && 4440.39 &314.51 &2328.68 &186.57\%  && 13174.82 &2009.98 &6.77 &750.26\%   \\ 
      &    $\quad\quad$ 30 && 6942.63 &509.68 &301.93 &198.95\%    && 6943.24 &509.79 &3004.19 &198.98\%  && 19584.82 &2692.96 &10.30   &743.32\%   \\ 
      & 6 $\quad$    10 && 2490.31 &278.27 &291.10 &154.87\%  && 2490.18 &278.26 &2862.97 &154.86\%  && 7053.26 &1143.83 &7.83 &621.87\%    \\ 
      &    $\quad\quad$ 20 && 4471.15 &335.79 &357.81 &176.77\%    && 4470.92 &335.78 &4936.60    &176.76\%   && 13242.86 &2009.08 &17.11 &719.76\%  \\ 
      &    $\quad\quad$ 30 && 7148.76 &633.26 &350.09 &176.54\%    && 7149.89 &633.30   &6130.73 &176.59\%   && 19744.29 &2724.49 &31.45 &663.79\%  \\ 
      & 9 $\quad$    10 && 2393.83 &216.69 &395.43 &178.31\%    && 2394.08 &216.70     &3894.24 &178.34\%   && 7054.97 &1144.17 &13.76 &720.22\%   \\ 
      &    $\quad\quad$ 20 && 4454.36 &319.36 &426.58 &188.01\%    && 4454.17 &319.35 &5944.89 &188.00\%  && 13245.42 &2009.77 &31.65 &756.44\%  \\ 
      &    $\quad\quad$ 30 && 7012.40   &532.57 &389.18 &192.77\%  && 7013.29 &532.65 &7441.00 &192.81\%   && 19743.99 &2724.87 &62.97 &724.31\%  \\ \\ 
0.6 & 3 $\quad$        10 && 4696.86 &423.90 &990.69 &308.58\%  && 7349.11 &181.66 &655.69 &539.30\%   && 9008.89 &1282.18 &9.43 &683.68\%    \\ 
      &    $\quad\quad$ 20 && 8839.00 &670.49 &937.68 &286.04\%  && 13933.37 &340.05 &886.40   &508.53\%   && 16932.93 &2236.03 &15.59 &639.53\%  \\ 
      &    $\quad\quad$ 30 && 13773.02 &1156.63 &1149.39 &305.30\% && 20813.01 &503.96 &1047.21 &512.47\% && 25784.98 &3342.88 &24.37 &658.78\%  \\ 
      & 6 $\quad$    10 && 4783.55 &475.90 &1047.17 &259.42\%    && 7362.41 &277.68 &748.56 &453.19\%   && 9030.75 &1303.43 &18.93 &578.55\%   \\ 
      &    $\quad\quad$ 20 && 8812.41 &676.11 &1261.63 &280.69\%   && 13922.8 &346.65 &1119.26 &501.45\%  && 16920.49 &2233.02 &40.74 &630.95\%  \\ 
      &    $\quad\quad$ 30 && 13840.24 &1209.13 &1038.13 &284.14\% && 20851.20  &597.72 &1346.35 &478.74\%  && 25790.33 &3351.70     &74.14 &615.82\%   \\ 
      & 9 $\quad$    10 && 4742.91 &447.86 &1421.52 &285.22\%   && 7339.79 &220.86 &923.35 &496.14\%   && 9043.80 &1304.17 &34.09 &634.54\%    \\ 
      &    $\quad\quad$ 20 && 8855.21 &682.44 &1418.54 &290.29\%   && 13981.64 &351.88 &1523.87 &516.24\% && 16941.19 &2234.26 &85.30   &646.68\%   \\ 
      &    $\quad\quad$ 30 && 13813.02 &1174.42 &1416.75 &301.25\% && 20856.74 &530.58 &2296.02 &505.86\% && 25813.97 &3353.10 &164.31 &649.86\%  \\ \\ 
5 & 3 $\quad$        10 && 9130.21 &1412.26 &1494.01 &322.02\%  && 9553.51 &940.09 &635.30   &341.59\%    && 10032.96 &1397.54 &10.51 &363.75\%  \\ 
      &    $\quad\quad$ 20 && 17223.85 &2440.65 &2292.23 &305.57\% && 17948.61 &1401.07 &774.83 &322.64\% && 18921.42 &2392.37 &15.67 &345.55\%  \\ 
      &    $\quad\quad$ 30 && 25589.33 &3675.81 &2419.88 &323.72\% && 26649.21 &2297.14 &1033.55 &341.27\%&& 28465.84 &3561.99 &26.40   &371.35\%   \\ 
      & 6 $\quad$    10 && 8978.55 &1422.53 &1814.95 &306.86\%  && 9392.60    &960.44 &805.39 &325.62\%    && 9849.95 &1395.57 &20.76 &346.34\%   \\ 
      &    $\quad\quad$ 20 && 16893.07 &2434.97 &3067.77 &325.17\% && 17603.28 &1402.02 &1128.86 &343.04\%&& 18564.48 &2388.98 &44.09 &367.24\%  \\ 
      &    $\quad\quad$ 30 && 25166.01 &3681.05 &2267.39 &332.38\% && 26155.27 &2300.44 &1210.80 &349.37\% && 27992.50    &3557.89 &80.53 &380.94\%   \\ 
      & 9 $\quad$    10 && 9106.33 &1477.99 &3678.81 &308.68\%  && 9518.79 &1032.71 &925.95 &327.19\%  && 9905.71 &1401.51 &30.04 &344.56\%   \\ 
      &    $\quad\quad$ 20 && 17113.32 &2526.51 &3859.10 &318.66\%  && 17778.34 &1534.64 &1338.55 &334.93\%&& 18638.58 &2397.41 &85.52 &355.98\%  \\ 
      &    $\quad\quad$ 30 && 25217.53 &3698.00 &4061.64 &330.59\%  && 26260.81 &2351.03 &1480.29 &348.40\%&& 28066.72 &3570.47 &158.72 &379.24\% \\  
\bottomrule
\end{tabular}
\end{adjustbox}
\caption{Performance of the \ac{famspd}, \ac{s2sspd} and \ac{rh2sspd} policies.}
\label{tab:results_det}
\end{table}

}

\begin{table}
\centering
\begin{adjustbox}{width=\textwidth}
\begin{tabular}{@{}lllccccccccccccc@{}}
\toprule
 & & \multicolumn{4}{c}{FA-MSP-D} & \phantom{abc}& \multicolumn{3}{c}{RH2SSP-D} & \phantom{abc} & \multicolumn{4}{c}{static2SSP-D}\\
\cmidrule{4-7} \cmidrule{9-11} \cmidrule{13-16}
  $\nu$ & $|I|\quad |J|$ && $\hat{z}$ & $\pm1.96 \hat{\sigma}/\sqrt{N}$ & time & gap 
                      && $\hat{z}$ & $\pm1.96 \hat{\sigma}/\sqrt{N}$ & gap 
                      && $\hat{z}$ & $\pm1.96 \hat{\sigma}/\sqrt{N}$ & time & gap\\ \midrule
0.001 & 3 $\quad$    10  && 2333.37  & 180.31  & 262.22  & 200.80\% && 2333.62  & 180.38  & 200.83\% && 7015.46  & 1143.26 & 0.53  & 804.38\% \\
      & $\quad\quad$ 20  && 4440.41  & 314.50  & 272.70  & 186.57\% && 4440.39  & 314.51  & 186.57\% && 13174.82 & 2009.98 & 0.74  & 750.26\% \\
      & $\quad\quad$ 30  && 6942.63  & 509.68  & 287.93  & 198.95\% && 6943.24  & 509.79  & 198.98\% && 19584.82 & 2692.96 & 0.96  & 743.32\% \\
      & 6 $\quad$    10  && 2490.31  & 278.27  & 273.53  & 154.87\% && 2490.18  & 278.26  & 154.86\% && 7053.26  & 1143.83 & 0.80  & 621.87\% \\
      & $\quad\quad$ 20  && 4471.15  & 335.79  & 330.32  & 176.77\% && 4470.92  & 335.78  & 176.76\% && 13242.86 & 2009.08 & 1.41  & 719.76\% \\
      & $\quad\quad$ 30  && 7148.76  & 633.26  & 307.56  & 176.54\% && 7149.89  & 633.30  & 176.59\% && 19744.29 & 2724.49 & 2.13  & 663.79\% \\
      & 9 $\quad$    10  && 2393.83  & 216.69  & 345.52  & 178.31\% && 2394.08  & 216.70  & 178.34\% && 7054.97  & 1144.17 & 0.89  & 720.22\% \\
      & $\quad\quad$ 20  && 4454.36  & 319.36  & 356.05  & 188.01\% && 4454.17  & 319.35  & 188.00\% && 13245.42 & 2009.77 & 1.35  & 756.44\% \\
      & $\quad\quad$ 30  && 7012.40  & 532.57  & 286.20  & 192.77\% && 7013.29  & 532.65  & 192.81\% && 19743.99 & 2724.87 & 1.84  & 724.31\% \\ \\
      
0.6   & 3 $\quad$    10  && 4696.86  & 423.90  & 975.49  & 308.58\% && 7349.11  & 181.66  & 539.30\% && 9008.89  & 1282.18 & 2.14  & 683.68\% \\
      & $\quad\quad$ 20  && 8839.00  & 670.49  & 915.29  & 286.04\% && 13933.37 & 340.05  & 508.53\% && 16932.93 & 2236.03 & 2.73  & 639.53\% \\
      & $\quad\quad$ 30  && 13773.02 & 1156.63 & 1118.92 & 305.30\% && 20813.01 & 503.96  & 512.47\% && 25784.98 & 3342.88 & 3.39  & 658.78\% \\
      & 6 $\quad$    10  && 4783.55  & 475.90  & 1009.75 & 259.42\% && 7362.41  & 277.68  & 453.19\% && 9030.75  & 1303.43 & 3.19  & 578.55\% \\
      & $\quad\quad$ 20  && 8812.41  & 676.11  & 1199.41 & 280.69\% && 13922.8  & 346.65  & 501.45\% && 16920.49 & 2233.02 & 4.93  & 630.95\% \\
      & $\quad\quad$ 30  && 13840.24 & 1209.13 & 946.69  & 284.14\% && 20851.20 & 597.72  & 478.74\% && 25790.33 & 3351.70 & 6.72  & 615.82\% \\
      & 9 $\quad$    10  && 4742.91  & 447.86  & 1299.33 & 285.22\% && 7339.79  & 220.86  & 496.14\% && 9043.80  & 1304.17 & 3.55  & 634.54\% \\
      & $\quad\quad$ 20  && 8855.21  & 682.44  & 1250.37 & 290.29\% && 13981.64 & 351.88  & 516.24\% && 16941.19 & 2234.26 & 7.08  & 646.68\% \\
      & $\quad\quad$ 30  && 13813.02 & 1174.42 & 1181.76 & 301.25\% && 20856.74 & 530.58  & 505.86\% && 25813.97 & 3353.10 & 8.08  & 649.86\% \\ \\
      
5     & 3 $\quad$    10  && 9130.21  & 1412.26 & 1476.25 & 322.02\% && 9553.51  & 940.09  & 341.59\% && 10032.96 & 1397.54 & 2.66  & 363.75\% \\
      & $\quad\quad$ 20  && 17223.85 & 2440.65 & 2266.71 & 305.57\% && 17948.61 & 1401.07 & 322.64\% && 18921.42 & 2392.37 & 3.41  & 345.55\% \\
      & $\quad\quad$ 30  && 25589.33 & 3675.81 & 2382.45 & 323.72\% && 26649.21 & 2297.14 & 341.27\% && 28465.84 & 3561.99 & 4.43  & 371.35\% \\
      & 6 $\quad$    10  && 8978.55  & 1422.53 & 1772.19 & 306.86\% && 9392.60  & 960.44  & 325.62\% && 9849.95  & 1395.57 & 5.23  & 346.34\% \\
      & $\quad\quad$ 20  && 16893.07 & 2434.97 & 3003.52 & 325.17\% && 17603.28 & 1402.02 & 343.04\% && 18564.48 & 2388.98 & 9.66  & 367.24\% \\
      & $\quad\quad$ 30  && 25166.01 & 3681.05 & 2169.74 & 332.38\% && 26155.27 & 2300.44 & 349.37\% && 27992.50 & 3557.89 & 9.66  & 380.94\% \\
      & 9 $\quad$    10  && 9106.33  & 1477.99 & 3562.41 & 308.68\% && 9518.79  & 1032.71 & 327.19\% && 9905.71  & 1401.51 & 4.49  & 344.56\% \\
      & $\quad\quad$ 20  && 17113.32 & 2526.51 & 3682.39 & 318.66\% && 17778.34 & 1534.64 & 334.93\% && 18638.58 & 2397.41 & 7.71  & 355.98\% \\
      & $\quad\quad$ 30  && 25217.53 & 3698.00 & 3814.67 & 330.59\% && 26260.81 & 2351.03 & 348.40\% && 28066.72 & 3570.47 & 11.20 & 379.24\% \\
\bottomrule
\end{tabular}
\end{adjustbox}
\caption{Performances of decision policies given by \ac{famspd}, \ac{s2sspd} and \ac{rh2sspd}.}
\label{tab:results_det}
\end{table}
From these results, we can see that on \textit{average}, the relative gap in $\hat{z}$ values to the \ac{cvd} solution is $264.01\%, 340.73\%$ and $573.99\%$ for \ac{famspd}, \ac{rh2sspd} and \ac{s2sspd}, respectively: specifically, $183.73\%, 183.75\%$ and $722.71\%$ when $\nu = 0.001$; $288.99\%, 501.32\%$ and $637.60\%$ when $\nu = 0.6$; and $319.29\%$, $337.12\%$ and $361.66\%$ when $\nu = 5$.

In summary, the \ac{famspd} has the best overall performance, followed by the \ac{rh2sspd}, and finally the \ac{s2sspd}. This is expected due to the level of adaptability associated with the corresponding decision policy for each approach. We also see the trade-off between the solution quality and the computational effort, which is reflected in the computational time spent on policy training. One important observation is that, compared to when $\nu = 0.6$, the difference in the average performance by the three policies shrinks when $\nu=5$, and the difference in the average performance between \ac{famspd} and \ac{rh2sspd} almost vanishes when $\nu = 0.001$. We discuss this observation in more details during our discussion in the sensitivity analysis in Subsection~\ref{subsubsec:sensitivity_analysis}.
\exclude{
From these results, we make the following observations:
\begin{itemize}
\item \emph{Upper bound estimates.} The $\hat{z}$ value of the \ac{famspd} model is less than that of the \ac{rh2sspd}  approach, which is less than that of the \ac{s2sspd} model. This observation is consistent across all instances -- except for a few, when $\nu = 0.001$, where the $\hat{z}$ value of the \ac{rh2sspd}  approach is smaller by some $\epsilon$ than that of the \ac{famspd} -- and the overall average of $\hat{z}$ across different instances is $10305.99, 12183.71$ and $16502.42$ for \ac{famspd}, \ac{rh2sspd} and \ac{s2sspd}, respectively. {\color{blue} XXX This seems repetitive as you have the last bullet point, shall we just remove this one and move the last bullet point here? In my view that is the most important observation. XXX}

\item \emph{Confidence intervals.} On average, the $\pm1.96 \hat{\sigma}/\sqrt{N}$ value of the \ac{rh2sspd}  approach has the smallest \ac{ci}, followed by the \ac{famspd} model, and then the \ac{s2sspd} model -- with values $773.75$, $1222.49$ and $2233.59$, respectively. In fact, we can see that the \ac{ci} of the \ac{rh2sspd}  is consistently less than (or equal to) that of the \ac{famspd} and \ac{s2ssp} across all of the instances. Whereas, the \ac{ci} of the \ac{famspd} is smaller than that of \ac{s2ssp} across all of the instances -- except those where $\nu=5$, in which case, the opposite is true. {\color{blue} XXX Yes, you have done an excellent job getting these statistics and making the comparisons, but what do we learn from these numbers? I don't see much information or explanation here. On the other hand, I see this more like a validation that the size of out-of-sample scenarios is acceptable -- the width of \ac{ci} is small compared to the sample mean. XXX}

\item \emph{Training and evaluation times.} In all of the instances, the \textit{time} required to train and evaluate the \ac{s2sspd} model is smaller than that of the \ac{famspd} and \ac{rh2sspd}. Additionally, while the time required to train and evaluate the \ac{famspd} model is significantly smaller than that of the \ac{rh2ssp} approach in all of the instances where $\nu=0.001$, the opposite is true in all of the remaining instances. {\color{blue} XXX It might be difficult to analyze if the training time and the evaluation time is not separated. XXX}

\item \emph{Performance relative to the \ac{cvd} solution.} On average, the relative gap in $\hat{z}$ values to the \ac{cvd} solution is $264.01\%, 340.73\%$ and $573.99\%$ for the \ac{famspd}, \ac{rh2sspd} and \ac{s2sspd}, respectively. In particular, the average relative gap in $\hat{z}$ values to the \ac{cvd} solution is
\begin{itemize}
    \item $183.73\%, 183.75\%$ and $722.71\%$ for \ac{famspd}, \ac{rh2sspd} and \ac{s2sspd}, respectively -- when $\nu = 0.001$;
    \item $288.99\%, 501.32\%$ and $637.60\%$ for \ac{famspd}, \ac{rh2sspd} and \ac{s2sspd}, respectively -- when $\nu = 0.6$.
    \item $319.29\%$, $337.12\%$ and $361.66\%$ for \ac{famspd}, \ac{rh2sspd} and \ac{s2sspd}, respectively -- when $\nu = 5$.
\end{itemize}
\end{itemize}

\paragraph{Summary.} In sum, the \ac{famspd} has the best overall performance, followed by the \ac{rh2sspd}, and finally the \ac{s2sspd}. This is somewhat expected due to the level of adaptability associated with the corresponding decision policy for each approach. We also see the tradeoff between the solution quality and the computational effort. One important observation is that, compared to when $\nu = 0.6$, the difference in the average performance by the three policies shrinks when $\nu=5$, and the difference in the average performance between \ac{famspd} and \ac{rh2sspd}  almost vanishes when $\nu = 0.001$. We discuss this observation in more details via a sensitivity analysis in Subsection~\ref{subsubsec:sensitivity_analysis}.
}

\subsection{Main results for the case of random landfall time}
\label{subsubsec:main_results_random}
\begin{figure}[htbp]
    \centering
    \subfigure[$\hat{z}$ values averaged deterministic.]{
    \includegraphics[scale=0.2]{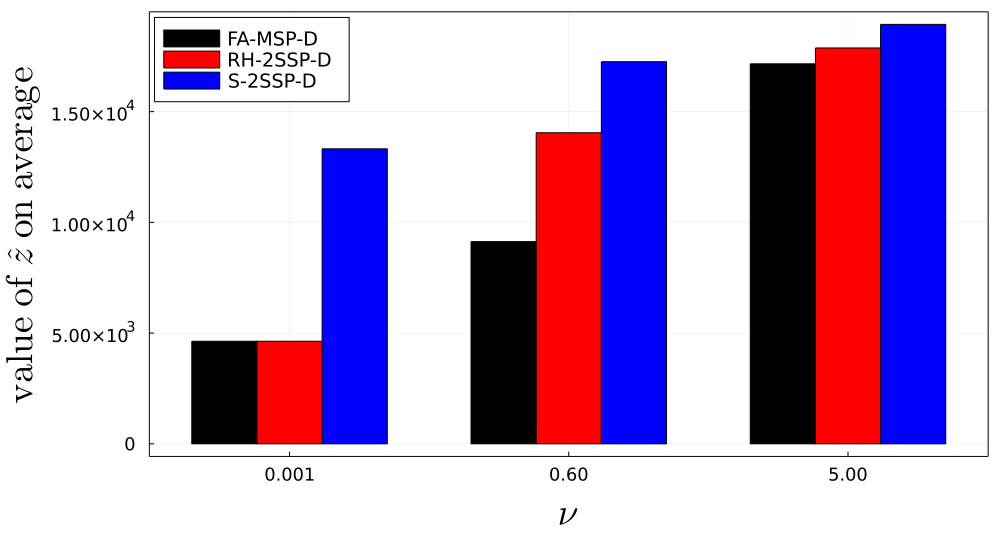}
    \label{fig:detT_allavg}} 
    \hspace{0.1cm}
    \subfigure[$\hat{z}$ values averaged random.]{
    \includegraphics[scale=0.2]{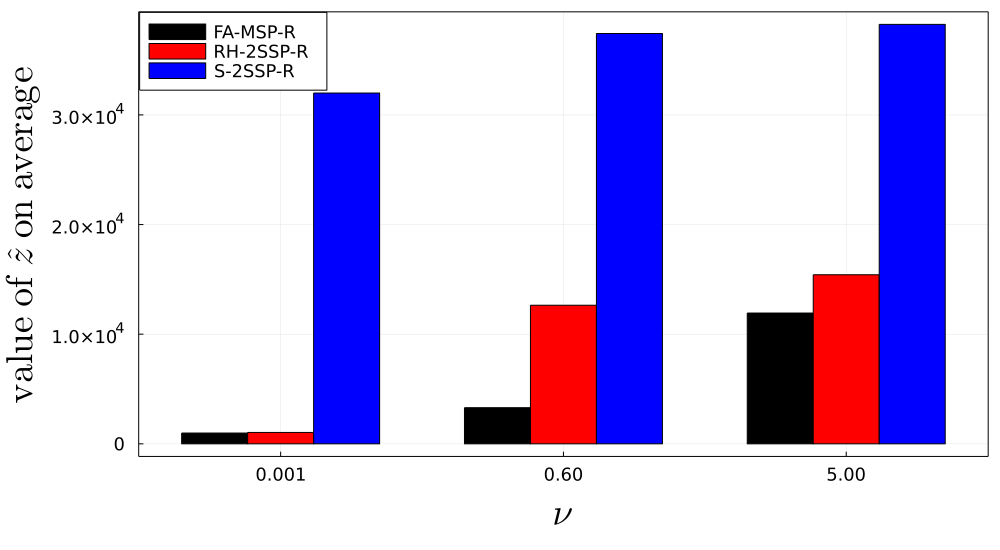}
    \label{fig:randT_allavg}} 
    \caption{$\hat{z}$ values averaged across all test instances for $\nu \in \{0.001, 0.6, 5\}$ when $T$ is deterministic (left) and random (right).}
    \label{fig:main_all_avg}
\end{figure}

In Table~\ref{tab:results_rand}, we report the numerical results for \ac{famspr}, \ac{s2sspr} and \ac{rh2sspr}. The first three columns describe the test instances based on different values of $\nu$, $|I|$ and $|J|$, and columns 4-7, 8-10 and 11-14 show the performance metrics results of \ac{famspr}, \ac{rh2sspr} and \ac{s2sspr}, respectively. Additionally, Figure~\ref{fig:randT_allavg} depicts the $\hat{z}$ values of different policies, \textit{averaged} across all the instances for different $\nu$ values (see also Figure~\ref{fig:randT_avg} in the Appendix).
\exclude{
\begin{table}
\centering
\footnotesize
\begin{adjustbox}{width=\textwidth}
\begin{tabular}{@{}lllcccccccccccccc@{}}
\toprule
 & & \multicolumn{4}{c}{FA-MSP-R} & \phantom{abc}& \multicolumn{4}{c}{RH2SSP-R} & \phantom{abc} & \multicolumn{4}{c}{static2SSP-R}\\
\cmidrule{4-7} \cmidrule{9-12} \cmidrule{14-17}
  $\nu$ & $|I|\quad |J|$ && $\hat{z}$ & $\pm1.96 \hat{\sigma}/\sqrt{N}$ & time & gap && $\hat{z}$ & $\pm1.96 \hat{\sigma}/\sqrt{N}$ & time & gap && $\hat{z}$ & $\pm1.96 \hat{\sigma}/\sqrt{N}$ & time & gap\\ \midrule
0.001 & 3 $\quad$    10    && 496.55   & 65.66 &950.83 &- && 496.55 &65.66 &1828.8 &- && 15106.07 &3057.67 &3.26 &2942.20\%  \\ 
      &    $\quad\quad$ 20 && 1001.85  & 126.01 &1070.06 &- && 1001.85 &126.01 &2981.48 &- && 30341.73 &5912.77 &4.25 &2928.56\%  \\ 
      &    $\quad\quad$ 30 && 1494.69  & 187.17 &1145.96 &- && 1494.69 &187.17 &3830.38 &- && 43727.69 &8252.10 &5.56 &2825.53\%  \\ 
      & 6  $\quad$      10    && 480.03   & 63.03 &2056.33 &- && 1001.85 &126.01 &4657.99 &108.71\% && 30341.73 &5912.77 &6.29 &6220.76\%  \\ 
      &    $\quad\quad$ 20 && 976.82   & 122.51 &2518.16 &- && 976.82 &122.51 &11257.16 &- && 29273.03 &5691.50 &8.26 &2896.77\%  \\ 
      &    $\quad\quad$ 30 && 1466.41  & 183.57 &2261.39 &- && 1466.41 &183.57 &14419.96 &- && 47311.64 &9027.15 &9.96 &3126.35\%  \\ 
      & 9  $\quad$      10    && 479.99   & 63.02 &2301.61 &- && 479.99 &63.02 &9708.90 &- && 15352.76 &3075.57 &9.50 &3098.58\%  \\ 
      &    $\quad\quad$ 20 && 976.76   & 122.51 &2200.61 &- && 976.76 &122.51 &15670.59 &- && 29938.36 &5856.80 &11.22 &2965.07\%  \\ 
      &    $\quad\quad$ 30 && 1466.35  & 183.56 &2870.88 &- && 1466.35 &183.56 &19751.16 &- && 46582.57 &9028.02 &14.31 &3076.77\%  \\ 
0.6  & 3  $\quad$      10    && 1690.89  & 236.74  & 10812.06 & 78.00\%  && 7265.66  & 293.95  & 25200.30 & 666.26\% && 19012.41 & 3211.42 & 6.84   & 1905.10\% \\
      & $\quad\quad$ 20    && 3339.33  & 439.97  & 10808.89 & 76.00\%  && 13417.20 & 606.42  & 37630.41 & 605.45\% && 37481.26 & 6321.66 & 9.80   & 1870.68\% \\
      & $\quad\quad$ 30    && 4971.89  & 651.25  & 10815.62 & 77.00\%  && 19469.59 & 936.16  & 49679.79 & 592.34\% && 55697.42 & 9545.94 & 12.72  & 1880.60\% \\
      & 6  $\quad$      10    && 1643.65  & 229.14  & 10815.07 & 76.00\%  && 6632.39  & 189.85  & 31177.71 & 609.25\% && 19043.77 & 3211.18 & 15.52  & 1936.49\% \\
      & $\quad\quad$ 20    && 3273.83  & 431.23  & 10807.32 & 75.00\%  && 12206.76 & 401.22  & 47359.21 & 554.15\% && 37524.74 & 6321.17 & 25.96  & 1910.91\% \\
      & $\quad\quad$ 30    && 4901.15  & 641.19  & 10811.98 & 77.00\%  && 17657.47 & 625.82  & 30218.39 & 539.14\% && 55748.90 & 9545.47 & 30.05  & 1917.93\% \\
      & 9  $\quad$      10    && 1646.45  & 229.75  & 10813.10 & 73.00\%  && 6744.97  & 207.28  & 36451.10 & 609.60\% && 19047.42 & 3211.36 & 20.15  & 1903.88\% \\
      & $\quad\quad$ 20    && 3275.83  & 431.24  & 10811.25 & 74.00\%  && 12430.11 & 436.63  & 46379.39 & 558.80\% && 37533.51 & 6321.39 & 37.38  & 1889.29\% \\
      & $\quad\quad$ 30    && 4896.78  & 641.64  & 10809.04 & 76.00\%  && 17962.17 & 675.50  & 42203.94 & 545.04\% && 55760.81 & 9545.60 & 46.27  & 1902.41\% \\
5  & 3  $\quad$      10    && 6275.95  & 1063.71 & 10806.65 & 223.00\% && 8629.04  & 775.24  & 6696.79  & 344.05\% && 19519.25 & 3215.24 & 14.05  & 904.46\%  \\
      & $\quad\quad$ 20    && 12218.53 & 1930.03 & 10811.83 & 223.00\% && 16181.22 & 1467.88 & 9744.21  & 327.34\% && 38460.07 & 6325.58 & 17.67  & 915.71\%  \\
      & $\quad\quad$ 30    && 18191.26 & 2786.20 & 10811.91 & 239.00\% && 23555.40 & 2185.96 & 12391.69 & 338.44\% && 56996.25 & 9547.83 & 25.36  & 960.87\%  \\
      & 6  $\quad$      10    && 6002.48  & 1027.01 & 10811.06 & 214.00\% && 8068.62  & 706.93  & 38267.95 & 321.72\% && 19461.32 & 3216.37 & 27.68  & 917.19\%  \\
      & $\quad\quad$ 20    && 11738.97 & 1871.23 & 10801.86 & 224.00\% && 15110.21 & 1314.32 & 30071.90 & 317.28\% && 38342.80 & 6327.21 & 51.70  & 958.85\%  \\
      & $\quad\quad$ 30    && 17561.68 & 2720.56 & 10804.93 & 243.00\% && 21866.83 & 1919.71 & 38733.54 & 327.02\% && 56843.74 & 9550.53 & 72.01  & 1010.05\% \\
      & 9 $\quad$    10    && 6028.92  & 1033.08 & 10807.62 & 202.00\% && 8125.81  & 723.59  & 25855.58 & 306.38\% && 19476.23 & 3216.57 & 33.34  & 874.03\%  \\
      & $\quad\quad$ 20    && 11776.55 & 1881.90 & 10808.76 & 216.00\% && 15237.02 & 1352.90 & 45690.74 & 309.30\% && 38355.39 & 6327.42 & 72.91  & 930.31\%  \\
      & $\quad\quad$ 30    && 17588.52 & 2741.80 & 10807.99 & 239.00\% && 22006.18 & 1975.57 & 57731.73 & 323.82\% && 56852.46 & 9551.04 & 101.32 & 994.94\%  \\ 
\bottomrule
\end{tabular}
\end{adjustbox}
\caption{Performance of the \ac{famspr}, \ac{s2sspr} and \ac{rh2sspr} policies.}
\label{tab:results_rand}
\end{table}
}

\begin{table}
\centering
\begin{adjustbox}{width=\textwidth}
\begin{tabular}{@{}lllccccccccccccc@{}}
\toprule
 & & \multicolumn{4}{c}{FA-MSP-R} & \phantom{abc}& \multicolumn{3}{c}{RH2SSP-R} & \phantom{abc} & \multicolumn{4}{c}{static2SSP-R}\\
\cmidrule{4-7} \cmidrule{9-11} \cmidrule{13-16}
  $\nu$ & $|I|\quad |J|$ && $\hat{z}$ & $\pm1.96 \hat{\sigma}/\sqrt{N}$ & time & gap 
                      && $\hat{z}$ & $\pm1.96 \hat{\sigma}/\sqrt{N}$ & gap 
                      && $\hat{z}$ & $\pm1.96 \hat{\sigma}/\sqrt{N}$ & time & gap\\ \midrule
0.001 & 3 $\quad$    10    && 496.55   & 65.66   & 948.54   & 00.00\%        && 496.55   & 65.66   & 00.00\%        && 15106.07 & 3057.67 & 1.63  & 2942.20\% \\
      & $\quad\quad$ 20    && 1001.85  & 126.01  & 1067.22  & 00.00\%        && 1001.85  & 126.01  & 00.00\%        && 30341.73 & 5912.77 & 2.18  & 2928.56\% \\
      & $\quad\quad$ 30    && 1494.69  & 187.17  & 1143.10  & 00.00\%        && 1494.69  & 187.17  & 00.00\%        && 43727.69 & 8252.10 & 3.21  & 2825.53\% \\
      & 6  $\quad$   10    && 480.03   & 63.03   & 2053.36  & 00.00\%        && 1001.85  & 126.01  & 108.71\% && 30341.73 & 5912.77 & 3.40  & 6220.76\% \\
      & $\quad\quad$ 20    && 976.82   & 122.51  & 2511.37  & 00.00\%        && 976.82   & 122.51  & 00.00\%        && 29273.03 & 5691.50 & 4.34  & 2896.77\% \\
      & $\quad\quad$ 30    && 1466.41  & 183.57  & 2257.64  & 00.00\%        && 1466.41  & 183.57  & 00.00\%        && 47311.64 & 9027.15 & 5.26  & 3126.35\% \\
      & 9  $\quad$   10    && 479.99   & 63.02   & 2298.09  & 00.00\%        && 479.99   & 63.02   & 00.00\%        && 15352.76 & 3075.57 & 5.91  & 3098.58\% \\
      & $\quad\quad$ 20    && 976.76   & 122.51  & 2196.55  & 00.00\%        && 976.76   & 122.51  & 00.00\%        && 29938.36 & 5856.80 & 6.49  & 2965.07\% \\
      & $\quad\quad$ 30    && 1466.35  & 183.56  & 2861.96  & 00.00\%       && 1466.35  & 183.56  & 00.00\%        && 46582.57 & 9028.02 & 8.51  & 3076.77\% \\ \\

0.6   & 3  $\quad$   10    && 1690.89  & 236.74  & - & 78.00\%  && 7265.66  & 293.95  & 666.26\% && 19012.41 & 3211.42 & 4.19  & 1905.10\% \\
      & $\quad\quad$ 20    && 3339.33  & 439.97  & - & 76.00\%  && 13417.20 & 606.42  & 605.45\% && 37481.26 & 6321.66 & 6.53  & 1870.68\% \\
      & $\quad\quad$ 30    && 4971.89  & 651.25  & - & 77.00\%  && 19469.59 & 936.16  & 592.34\% && 55697.42 & 9545.94 & 8.80  & 1880.60\% \\
      & 6  $\quad$   10    && 1643.65  & 229.14  & - & 76.00\%  && 6632.39  & 189.85  & 609.25\% && 19043.77 & 3211.18 & 12.43 & 1936.49\% \\
      & $\quad\quad$ 20    && 3273.83  & 431.23  & - & 75.00\%  && 12206.76 & 401.22  & 554.15\% && 37524.74 & 6321.17 & 22.03 & 1910.91\% \\
      & $\quad\quad$ 30    && 4901.15  & 641.19  & - & 77.00\%  && 17657.47 & 625.82  & 539.14\% && 55748.90 & 9545.47 & 25.48 & 1917.93\% \\
      & 9  $\quad$   10    && 1646.45  & 229.75  & - & 73.00\%  && 6744.97  & 207.28  & 609.60\% && 19047.42 & 3211.36 & 16.62 & 1903.88\% \\
      & $\quad\quad$ 20    && 3275.83  & 431.24  & - & 74.00\%  && 12430.11 & 436.63  & 558.80\% && 37533.51 & 6321.39 & 32.48 & 1889.29\% \\
      & $\quad\quad$ 30    && 4896.78  & 641.64  & - & 76.00\%  && 17962.17 & 675.50  & 545.04\% && 55760.81 & 9545.60 & 40.47 & 1902.41\% \\ \\
      
5     & 3  $\quad$   10    && 6275.95  & 1063.71 & - & 223.00\% && 8629.04  & 775.24  & 344.05\% && 19519.25 & 3215.24 & 10.10 & 904.46\%  \\ 
      & $\quad\quad$ 20    && 12218.53 & 1930.03 & - & 223.00\% && 16181.22 & 1467.88 & 327.34\% && 38460.07 & 6325.58 & 13.33 & 915.71\%  \\
      & $\quad\quad$ 30    && 18191.26 & 2786.20 & - & 239.00\% && 23555.40 & 2185.96 & 338.44\% && 56996.25 & 9547.83 & 20.04 & 960.87\%  \\
      & 6  $\quad$   10    && 6002.48  & 1027.01 & - & 214.00\% && 8068.62  & 706.93  & 321.72\% && 19461.32 & 3216.37 & 23.13 & 917.19\%  \\
      & $\quad\quad$ 20    && 11738.97 & 1871.23 & - & 224.00\% && 15110.21 & 1314.32 & 317.28\% && 38342.80 & 6327.21 & 46.29 & 958.85\%  \\
      & $\quad\quad$ 30    && 17561.68 & 2720.56 & - & 243.00\% && 21866.83 & 1919.71 & 327.02\% && 56843.74 & 9550.53 & 65.45 & 1010.05\% \\
      & 9 $\quad$    10    && 6028.92  & 1033.08 & - & 202.00\% && 8125.81  & 723.59  & 306.38\% && 19476.23 & 3216.57 & 28.21 & 874.03\%  \\
      & $\quad\quad$ 20    && 11776.55 & 1881.90 & - & 216.00\% && 15237.02 & 1352.90 & 309.30\% && 38355.39 & 6327.42 & 66.29 & 930.31\%  \\
      & $\quad\quad$ 30    && 17588.52 & 2741.80 & - & 239.00\% && 22006.18 & 1975.57 & 323.82\% && 56852.46 & 9551.04 & 92.36 & 994.94\%  \\
\bottomrule
\end{tabular}
\end{adjustbox}
\caption{Performances of decision policies given by \ac{famspr}, \ac{s2sspr} and \ac{rh2sspr}.}
\label{tab:results_rand}
\end{table}

From these results, we can see that on average, the relative gap in $\hat{z}$ values to the \ac{cvr} solution is $100.14\%, 307.56\%$ and $2061.64\%$ for the \ac{famspr}, \ac{rh2sspr} and \ac{s2sspr}, respectively: specifically, $0.00\%$,	$12.08\% $	and $3342.29\%$ when $\nu = 0.001$; $75.78\%$, $586.67\%$	and $1901.92\%$ when $\nu = 0.6$; and $224.63\%$, $323.93\%$	and $940.71\%$ when $\nu = 5$.

In summary, most of the observations made in Subsection \ref{subsubsec:main_results_deterministic}, for the case of deterministic landfall time, carry over to the case of random landfall time: \ac{famspr} outperforms \ac{rh2sspr}, which in turn outperforms \ac{s2sspr}. Additionally, the difference in the $\hat{z}$ values between \ac{famspr} and \ac{rh2sspr} almost vanishes when $\nu = 0.001$ and gets smaller when $\nu=5$. Nevertheless, the additional layer of uncertainty (randomness in $T$) gives rise to additional insights. In particular, the randomness in $T$ emphasizes the influence of adaptability on the quality of the decision policy. As we can see from Table~\ref{tab:results_rand} and Figure~\ref{fig:randT_allavg}, the difference in the performance between \ac{s2sspr} and \ac{famspr} (or \ac{rh2sspr}) is amplified. More importantly, the performance of \ac{famspr} is as good as the clairvoyant (\ac{cvr} solution) across all of the instances when $\nu=0.001$; this is also the case for the \ac{rh2sspr} approach, except for one instance. We suspect that this is because when $\nu=0.001$, the increase in the logistics cost over time is so small that it encourages a ``wait-and-see'' policy, where the \ac{dm} delays all the prepositioning of the relief items until the hurricane's attributes at landfall become clear. Nevertheless, the adaptability in the decision is incumbent on such a ``wait-and-see'' policy, that is, the \ac{dm} must have the leverage of being able to \textit{wait} and \textit{adapt} their decision according to what they observe over time. We investigate these claims further in the next subsection. 
\exclude{
From these results, we make the following observations. {\color{blue} XXX Please see my previous comments for the deterministic $T$ case and apply them here, too. XXX}
\begin{itemize}
\item \emph{Upper bound estimates.} As is the case in the policies with deterministic time of landfall, in terms of the $\hat{z}$ values, the \ac{famspr} outperforms the \ac{rh2sspr} approach, which in turn outperforms the \ac{s2sspr}. This observation is consistent across all instances -- without any exceptions -- and the overall average of $\hat{z}$ across different instances is $5402.30, 9701.03$ and $35893.83$ for \ac{famspd}, \ac{rh2sspd} and \ac{s2sspd}, respectively.

\item \emph{Confidence intervals.} On average, the $\pm1.96 \hat{\sigma}/\sqrt{N}$ value of the \ac{rh2sspr} approach has the smallest \ac{ci}, followed by the \ac{famspr} model, and then the \ac{s2sspr} model -- with values $665.74$, $818.69$ and $6308.42$, respectively. This observation, which is based on the average values, is consistent across all instances when comparing the \ac{s2sspr} with the \ac{famspr} and \ac{rh2sspr} (i.e., \ac{ci} of \ac{s2sspr} > \ac{ci} of \ac{rh2sspr}, and \ac{ci} of \ac{s2sspr} > \ac{ci} of \ac{famspr} ). However, it is not consistent across all instances when comparing the \ac{famspr} with \ac{rh2sspr}. Specifically, when $\nu = 0.001$, the \ac{ci} values of the \ac{famspr} are equal to that that of the \ac{rh2sspr} across all instances, except for two instances where \ac{famspr} is smaller. Whereas, when $\nu = 5.00$, the \ac{ci} values of the \ac{famspr} are consistently larger than that of the \ac{rh2sspr}. For the case when $\nu = 0.60$, no one policy prevails, with the \ac{famspr} being smaller in most of the instances.

\item \emph{Training and evaluation times.} On one hand, the \textit{time} required to train and evaluate the \ac{s2sspr} model is significantly smaller than that of the \ac{famspr} and \ac{rh2sspr} in all of the instances. On the other hand, except for the smallest two instances when $\nu=5.00$, the \textit{time} required to train and evaluate \ac{famspr} is consistently smaller than that of the \ac{rh2sspr}, with the difference peaking at the instances where $\nu=0.60$. However, as we can see from~\ref{tab:results_rand}, the \textit{time} values for the \ac{famspr} exceed $3\times 60^2$ in all of the instances where $\nu \in \{0.60, 5.00\}$. In those instances, the training was terminated when reaching the $3\times 60^2$ seconds time-limit, and the difference between the reported \textit{time} and $3\times 60^2$ is how long it took to evaluate the policy for the $N$ sample paths. Nevertheless, the \textit{time} values pertaining to the \ac{rh2sspr} reflects the cumulative time for solving the \ac{s2ssp} model for every roll and every sample paths, i.e., $N\times T$. Note, however, that \textit{time} values pertaining to the \ac{s2sspr} gives an \ac{ub} on the time required for each of the $N\times T$ models solved in every roll. This is because, the model solved in the \ac{s2sspr} coincides with the first roll of the \ac{rh2sspr}, and the number of aggregated periods recedes as we roll forward. 

\item \emph{Performance relative to the \ac{cvr} solution.} On average, the relative gap in $\hat{z}$ values to the \ac{cvr} solution is $100.14\%, 307.56\%$ and $2061.64\%$ for the \ac{famspd}, \ac{rh2ssp} and \ac{s2ssp}, respectively. In particular, the average relative gap in $\hat{z}$ values to the \ac{cvd} solution is
\begin{itemize}
    \item $0.00\%$,	$12.08\% $	and $3342.29\%$ for \ac{famspr}, \ac{rh2sspr}, and \ac{s2sspr}, respectively, when $\nu = 0.001$;
    \item $75.78\%$,	$586.67\%$	and $1901.92\%$ for \ac{famspr}, \ac{rh2sspr}, and \ac{s2sspr}, respectively, when $\nu = 0.6$.
    \item $224.63\%$,	$323.93\%$	and $940.71\%$ for \ac{famspr}, \ac{rh2sspr}, and \ac{s2sspr}, respectively, when $\nu = 5$.
\end{itemize}
\end{itemize}

\paragraph{Summary.} 
In summary, most of the conclusions in Subsection \ref{subsubsec:main_results_deterministic}, when $T$ is deterministic, carry over to the case where $T$ is random: \ac{famspr} outperforms \ac{rh2sspr}, which in turn outperforms \ac{s2sspr}. Additionally, the difference in the $\hat{z}$ values between the \ac{famspr} and \ac{rh2sspr} vanishes when $\nu = 0.001$ and gets smaller when $\nu=5.00$. Nevertheless, the additional layer of uncertainty (randomness in $T$) gives rise to additional insights. In particular, the randomness in $T$ emphasizes the influence of adaptability on the quality of the decision policy. As we can see from Table~\ref{tab:results_rand} and Figure~\ref{fig:randT_allavg}, the difference in the performance between \ac{s2sspr} and \ac{famspr} (or \ac{rh2sspr}) is amplified. More importantly, the performance of \ac{famspr} is as good as the clairvoyant (\ac{cvr} solution) across all of the instances when $\nu=0.001$; this is also the case for the \ac{rh2sspr} approach, except for one instance. The reason for this is that, when $\nu=0.001$, the increase in the logistics cost over time is so small that it encourages a ``wait-and-see'' policy, where the \ac{dm} delays all the prepositioning of the relief items until the hurricane's attributes at landfall become clear. Nevertheless, adaptability in the decision is incumbent on such a reactive ``wait-and-see'' policy. That is, the \ac{dm} must have the leverage of being able to \textit{wait and adapt} their decision according to what they observe over time. We investigate these observations further via a sensitivity analysis in Subsection~\ref{subsubsec:sensitivity_analysis}. 
}

\subsection{Sensitivity analysis}
\label{subsubsec:sensitivity_analysis}
In this subsection, we present and discuss the results for the sensitivity analysis performed on different policies for instances with different cost-scaling factor $\nu$ and different initial intensity level of the hurricane. We discuss this for both the deterministic and random landfall time cases, according to the following.
\begin{itemize}
    \item The relative \textit{gaps} in the out-of-sample performance between solutions obtained by different policies and the corresponding \ac{cv} solutions, as a function of $\nu$, shown in Figures~\ref{fig:detT_gaps_combined} and~\ref{fig:randT_gaps_combined}.
    \item The average amount of relief items procured from the \ac{mdc} in each period $t$, denoted by $\bar{f}_t$, as a function of $\nu$, where
    \begin{equation}
    \label{eq:avg_procurement}
        \bar{f}_t = \sum_{n \in N} \frac{1}{N}\sum_{i \in I} f^n_{0i,t}
    \end{equation}
   and $f^n_{0i,t}$ is the amount of relief items shipped from the \ac{mdc} to the SP $i \in I$ at time $t=1, \dots, T$, when evaluating the policy on the out-of-sample scenario $n=1, \dots, N$. We evaluate each policy with initial hurricane's intensity level $\alpha_1 \in \{1, 3, 5\}$ for \ac{famspd} in Figure~\ref{fig:detT_sensitivity}, and for \ac{famspr} in Figure~\ref{fig:randT_sensitivity} (see also Tables~\ref{tab:sens_results_det} and~\ref{tab:sens_results_random} in the Appendix).
\end{itemize}

\subsubsection{Sensitivity analysis with respect to $\nu$.} 
\begin{figure}[htbp]
    \centering
    \subfigure[Relative gap with \ac{cv}-D when $T$ is deterministic.]{
    \includegraphics[scale=0.225]{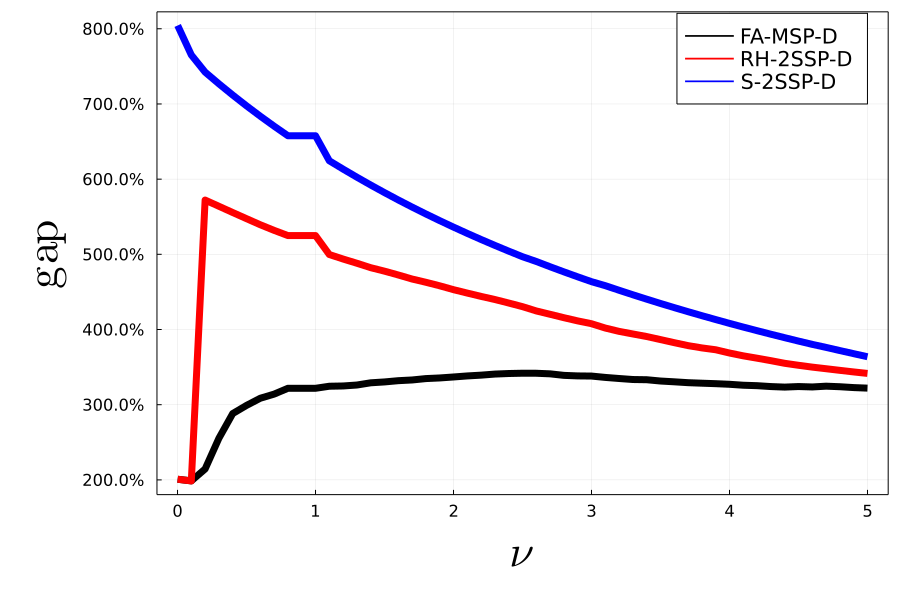}
    \label{fig:detT_gaps_combined}} 
    \hspace{0.1cm}
    \subfigure[Relative gap with \ac{cv}-R when $T$ is random.]{
    \includegraphics[scale=0.225]{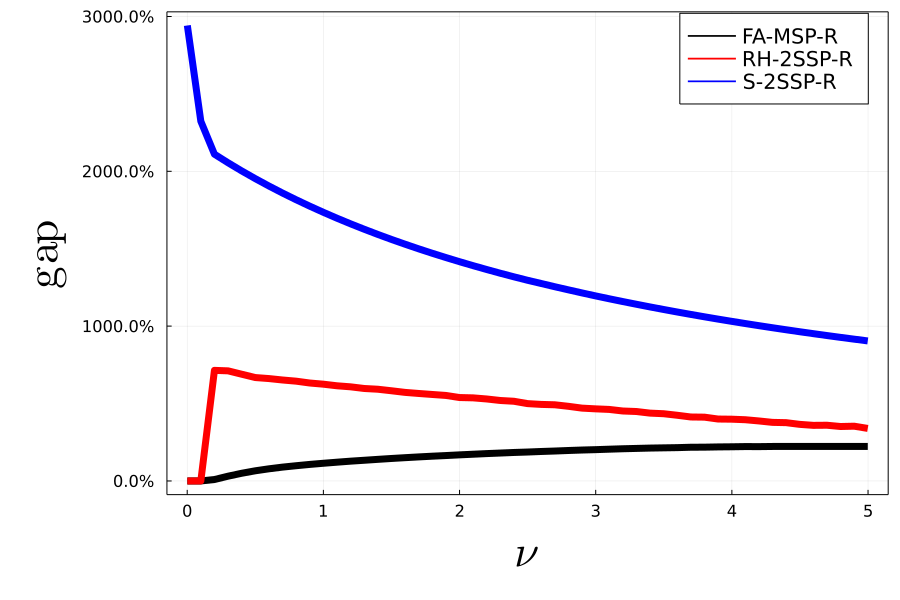}
    \label{fig:randT_gaps_combined}} 
    \caption{Relative gap in $\hat{z}$ values compared to the CV solutions as a function of $\nu$ when $T$ is deterministic (left) and random (right).}
    \label{fig:sens-gaps-combined}
\end{figure}

\begin{figure}[htbp]
    \centering
    \subfigure[$\alpha_1 = 1$.]{
    \includegraphics[scale=0.175]{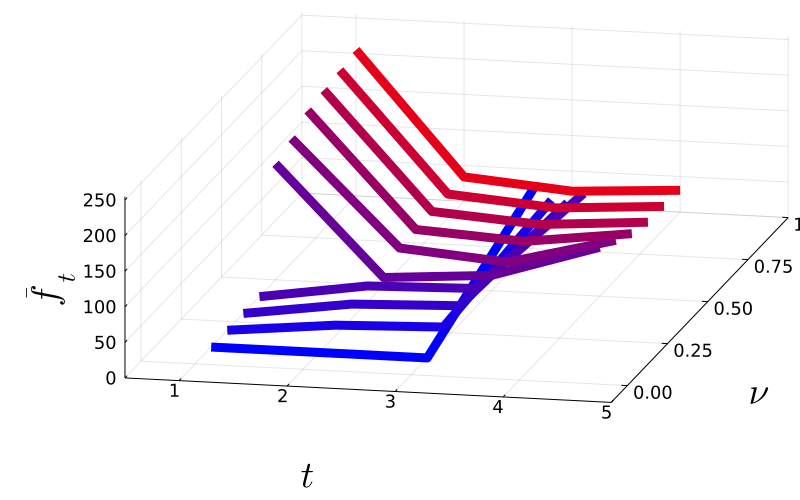}
    \label{fig:detT_sensitivity_p1}} 
    \hspace{0.1cm}
    \subfigure[$\alpha_1 = 3$.]{
    \includegraphics[scale=0.175]{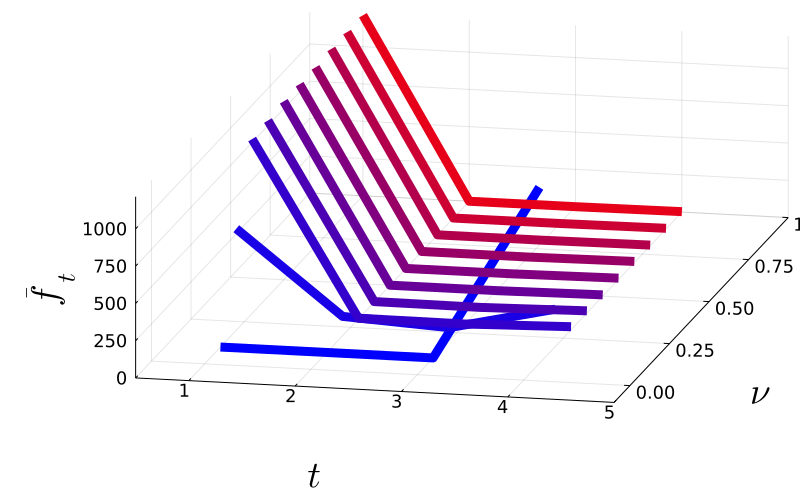}
    \label{fig:detT_sensitivity_p3}} 
    \hspace{0.1cm}
    \subfigure[$\alpha_1 = 5$.]{
    \includegraphics[scale=0.175]{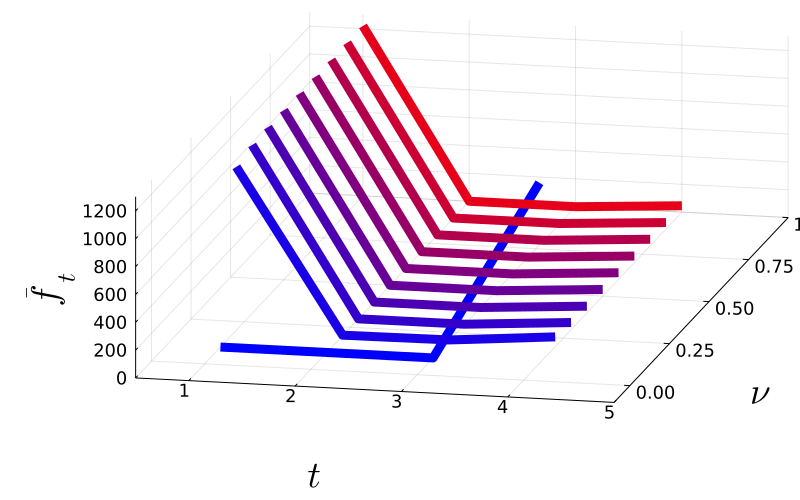}
    \label{fig:detT_sensitivity_p5}}
    \caption{$\bar{f}_t$~\eqref{eq:avg_procurement} obtained using the \ac{famspd} model as a function of $\nu$ for $t=1, \dots, T$, and $\alpha_1 \in \{1, 3, 5\}$.}
    \label{fig:detT_sensitivity}
\end{figure}

\begin{figure}[htbp]
    \centering
    \subfigure[$\alpha_1 = 1$.]{
    \includegraphics[scale=0.175]{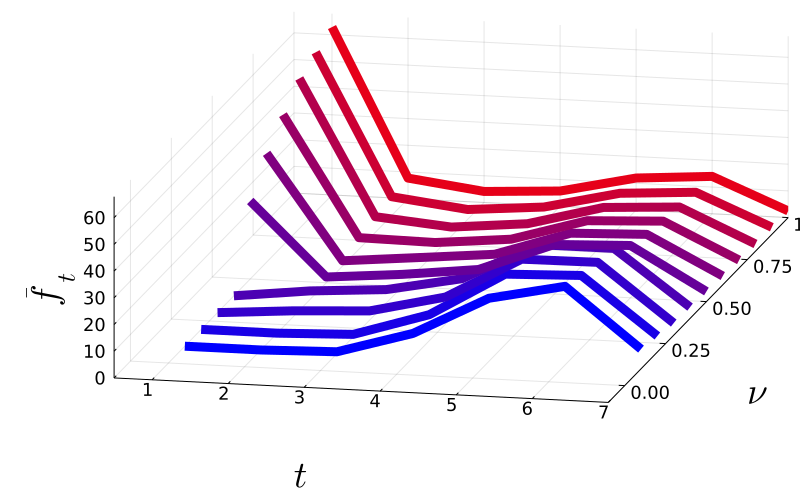}
    \label{fig:randT_sensitivity_p1}}
    \hspace{0.1cm}
    \subfigure[$\alpha_1 = 3$.]{
    \includegraphics[scale=0.175]{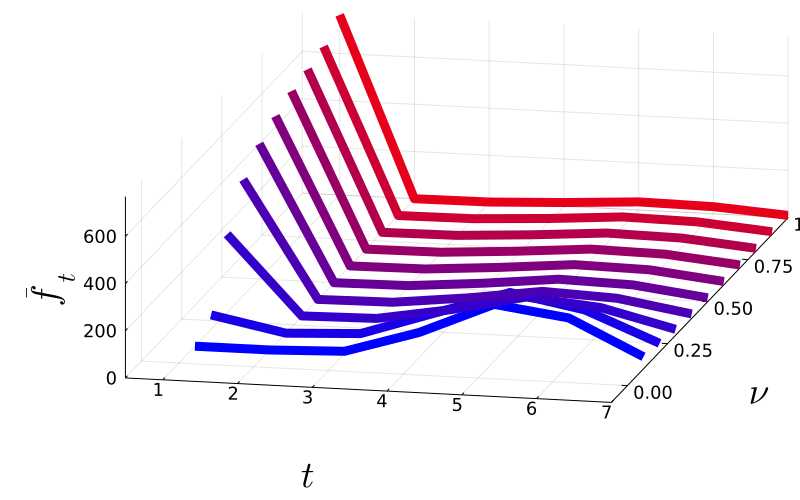}
    \label{fig:randT_sensitivity_p3}}
    \hspace{0.1cm}
    \subfigure[$\alpha_1 = 5$.]{
    \includegraphics[scale=0.175]{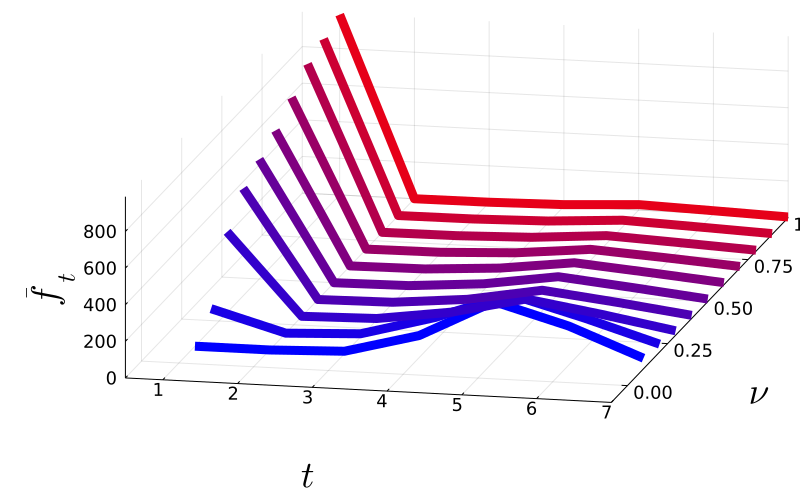}
    \label{fig:randT_sensitivity_p5}}
    \caption{$\bar{f}_t$~\eqref{eq:avg_procurement} obtained using the \ac{famspr} model as a function of $\nu$ for $t=1, \dots, T$, and $\alpha_1 \in \{1, 3, 5\}$.}
    \label{fig:randT_sensitivity}
\end{figure}

As we can see from Figures~\ref{fig:detT_gaps_combined} and~\ref{fig:randT_gaps_combined}, the decision policies given by \ac{famspd} and \ac{famspr} models are superior to their \ac{rh2sspd} and \ac{rh2sspr} counterparts over the entire range of $\nu$ values. Similarly, \ac{rh2sspd} and \ac{rh2sspr} maintain their superiority over their \ac{s2sspd} and \ac{s2sspr} counterparts over the entire range of $\nu$ values. It is important to note, however, that their differences are not always paramount, especially in the deterministic landfall case. On one hand, as discussed in Subsection~\ref{subsubsec:main_results_random}, the insignificance of a small cost scaling factor (e.g., $\nu=0.001$) encourages a ``wait-and-see'' policy where all of the prepositioning decisions are postponed until the hurricane's attributes at landfall become more clear. From Figures~\ref{fig:detT_sensitivity} and~\ref{fig:randT_sensitivity}, we can see that this corresponds to postponing the procurement decisions until time $T-1$. On the other hand, we can see from Figures~\ref{fig:detT_gaps_combined} and~\ref{fig:randT_gaps_combined} that when $\nu \to 5$, not only the \ac{famsp} and \ac{rh2ssp} policies yield similar performances, but the \ac{s2ssp} solution yields competitive performance as well. The reason for this is precisely opposite to what happens at the other end of the spectrum (when the value of $\nu$ is small): as we can see from Figures~\ref{fig:detT_sensitivity} and~\ref{fig:randT_sensitivity}, when $\nu$ gets large, most (if not all) of the relief items are procured at the first period. In other words, the \ac{dm} commits to a policy where all the prepositioning decisions are made as early as possible to hedge against the steep increase in the operational cost over time due to a large $\nu$ value. This leads to similar performances obtained by different policies. Although small differences in their performances can still be observed, these differences can be attributed to the level of adaptability inherited in each policy which can be leveraged to modify the initial procurement decision slightly, depending on how the stochastic process unfolds. 

An alternative perspective on the above assertions is that, for very large or very small values of $\nu$, the problem becomes similar to a two-stage optimization problem, where the first-stage corresponds to the period where all the procurement decisions are made, and the second-stage corresponds to the recourse decisions. From this perspective, this second-stage period is always defined as the landfall period $t=T$, while the definition of the first-stage period depends on the value of $\nu$. Specifically: when $\nu \to 0$, the first-stage period $t \to T-1$ (the``wait-and-see'' policy); and when $\nu \to \infty$, the first-stage period $t \to 1$ (the ``early-commitment'' policy).

At the first glance, these assertions could seem rooted in the $\nu$ values only. However, these observations are concurrently owing to Assumption~\ref{assm:demand}, which supposes that the demand occurs only at the landfall time, i.e., ${d}_{j,t} = 0, \forall \; j \in J,\; t=1, \dots, T-1$. This can be substantiated by comparing how well \ac{famspd} performs against \ac{rh2sspd} and \ac{s2sspd} with how well \ac{famspr} performs against \ac{rh2sspr} and \ac{s2sspr}. As we can see from Figure~\ref{fig:sens-gaps-combined}, when $\nu$ gets close to $0$, \ac{famspd} and \ac{rh2sspd}  have very similar performances, and so does \ac{famspr} and \ac{rh2sspr}. Additionally, as $\nu$ gets larger, the difference between \ac{famspd}, \ac{rh2sspd} and \ac{s2sspd} gradually gets smaller, and so does the difference between \ac{famspr}, \ac{rh2sspr} and \ac{s2sspr}. However, we can see that as $\nu$ gets larger, this difference is shrinking at a slower pace in the random landfall time case (Figure~\ref{fig:randT_gaps_combined}) than the deterministic (Figure~\ref{fig:detT_gaps_combined}) landfall time case. In fact, we can reinterpret Assumption~\ref{assm:demand} as $d_{j,t} = 0, \forall \; j \in J, \; t=1, \dots, T-1$ with probability $1$ in the deterministic landfall time case, while in the random landfall time case, there is a positive chance that $d_{j,t} > 0$ for multiple periods $t < T_{\max}$. In other words, the situation that $T$ is random resembles the situation where demand for relief items may occur in multiple stages. Hence, even when the cost-scaling factor $\nu$ is large, we cannot treat the problem as a two-stage optimization problem due to the uncertainty of the timing of demand realization. 

\subsubsection{Sensitivity analysis with respect to $\alpha_1$.}
From Figures~\ref{fig:detT_sensitivity} and~\ref{fig:randT_sensitivity}, we can see that this experiment on the value of initial intensity level $\alpha_1$ amplifies the impact that the cost-scaling factor $\nu$ has on the performances of different policies. For instance, when $\nu$ is small (e.g., $0.1 < \nu < 0.5$), we can see that when $\alpha_1 = 1$ (very low intensity), most of the prepositioning decisions are postponed to the end. However, as the value of $\alpha_1$ increases, less relief items are procured towards the end of the planning horizon and more are procured early on. In other words, a higher initial intensity level of the hurricane expedites the process of transitioning from the ``wait-and-see'' policy to the ``early-commitment'' policy as the value of $\nu$ increases. Furthermore, we can also see that this transition happens more swiftly in the case of a deterministic landfall time (see Figure~\ref{fig:detT_sensitivity}) in comparison with the case of a random landfall time (see Figure~\ref{fig:randT_sensitivity}). At the highest level of intensity $\alpha_1=5$, when $\nu = 0.001$, $\bar f_t = 0$ for all the non-terminal stages except $t=T-1$. However, in the random landfall case, this transition happens more slowly. We suspect that this can be attributed to the additional layer of uncertainty about the timing of the demand realization, which undermines the impact that $\nu$ and $\alpha_1$ has on the sequential nature of the decision making process.

\section{Conclusion}\label{sec:conclusion}

In this paper, we considered the problem of prepositioning relief items prior to an impending hurricane landfall, which we referred to as the hurricane disaster relief logistics planning (\ac{hdrlp}) problem. \exclude{In this problem, the goal of the \ac{dm} is to preposition relief items at different \ac{sps} such that, when/if the hurricane makes landfall, the prepositioned items are delivered from the \ac{sps} to the affected individuals (\ac{dps}) at a minimum cost. Here, the cost consists of logistics costs for purchasing, delivering, and storing the items, as well as a penalty for failing to satisfy the demand in case any shortage is present.} We have studied this problem in two different settings, depending on if the time of landfall is deterministic or random. In each setting, we assume that the demand for relief items can be derived from the hurricane's characteristics at the time of its landfall, namely its intensity and landfall location. We also assume that the evolution of the hurricane characteristics, including the timing of its landfall, is modeled as a discrete time \ac{mc} model. 

To solve the \ac{hdrlp} problem, we have proposed two \ac{famsp} models, which we refer to as \ac{famspd} and \ac{famspr} when the time of landfall is deterministic and random, respectively. Both \ac{famsp} models induce optimal offline decision policies that allow the \ac{dm} to make adaptive logistics operational decisions in every period of the planning horizon, depending on the most recently observed hurricane characteristics. In particular, our proposed \ac{famspr} model provides a novel extension to multistage stochastic programming with a random number of stages where the underlying stochastic process is assumed to be stage-wise dependent. Moreover, due to the computational challenges in exactly solving the \ac{msp} models, we have introduced alternative approximate decision policies. Specifically, we have introduced a \ac{2ssp} model where the corresponding decisions are implemented in a static manner, and an \ac{rh} approach that employs the \ac{2ssp} model as the look-ahead model in each period of the planning horizon to achieve an online decision policy. These approaches are referred to, respectively, as \ac{s2sspd} and \ac{rh2sspd} when the time of landfall is deterministic; and \ac{s2sspr} and \ac{rh2sspr} when the time of landfall is random.

\exclude{Our numerical results have shown that although optimal solutions to the \ac{2ssp} models can be obtained efficiently, the corresponding static decision policies were consistently outperformed by the offline decision policies given by \ac{famspd} and \ac{famspr} models and the online decision policies obtained by the \ac{rh2sspd} and \ac{rh2sspr} approaches.} 
We have conducted an extensive set of numerical experiments to validate the proposed approaches. One key insight of our numerical results and sensitivity analysis is how the performances by different policies depend on how fast the unit logistics cost increases over the planning horizon. On one hand, when the cost increase is insignificant over time, the performances of decision policies associated with \ac{famsp} and \ac{rh2ssp} are similar for both situations when the time of landfall is deterministic and random. This is due to the fact that, insignificant increase in the cost encourages a latent ``wait-and-see'' policy where almost all of the logistics decisions are postponed until the hurricanes gets close to its landfall, so that its characteristics becomes more clear. However, this ``wait-and-see'' policy must be met with the \ac{dm}'s ability to adapt their decisions at the later stages when the hurricane characteristics becomes clearer, which justifies why the same cannot be said about the \ac{s2sspd} and \ac{s2sspr} models. On the other hand, a similar observation can also be found when the unit logistics cost scales up significantly over time, where the performances of all different policies become more similar. This is due to the fact that the significant cost increase renders a \ac{dm} who commits to making almost all the logistics decisions early on in the planning horizon (a ``here-and-now'' policy), before the logistics costs become untenable. Overall, the performances of different policies clearly depend on the levels of adaptability inherited in different approaches: with \ac{famspd} and \ac{famspr} having the best performances, followed by \ac{rh2sspd} and \ac{rh2sspr}, and finally \ac{s2sspd} and \ac{s2sspr}.

We have identified several directions to pursue for future research. One interesting avenue is to investigate cases where the demand for relief items may occur over multiple periods, rather than just occurring at the landfall stage only. Having to satisfy the demand over multiple periods will likely amplify the merits of the \ac{famsp} models. In addition, it would be interesting to consider the more realistic situations where the \ac{dm} has to make decisions regarding the selection of \ac{sps} and the timings of their activation. Finally, from a risk preference perspective, it would be of interest to study risk-averse \ac{msp} models for the \ac{hdrlp} problem.

\paragraph{Acknowledgements.} Clemson University is acknowledged for generous allotment of compute time on Palmetto cluster. The authors acknowledge partial support by the National Science Foundation [Grant CMMI 2045744]. Any opinions, findings, and conclusions or recommendations expressed in this material are those of the authors and do not necessarily reflect the views of the National Science Foundation.
\bibliographystyle{plain}  
\bibliography{references}
\clearpage
\begin{appendix}
\section{Nested Benders decomposition.}
\label{sec:nested_benders_decompositions}
The \ac{dpe}~\eqref{eq:FOSDDP} decomposes the nested formulation~\eqref{eq:MSP} over different stages $t$, which serves as a basis for the so-called nested Benders decomposition~\cite{birge2011introduction}. In a nutshell, the basic idea of nested Benders decomposition is to iteratively approximate the expected cost-to-go functions $\Q^{\xi_t}_{t+1}(\cdot)$ by a collection of cutting-plane approximations. The algorithmic machinery for this procedure consists of the following two main steps.
\begin{enumerate}
    \item \textit{Forward} step: which simulates a trajectory of Markovian states by going forward in time from $t = 1$ to $t = T$, solves the respective stage-$t$ problem based on the current approximation of the expected cost-to-go function $\check \Q^{\xi_t}_{t+1}(\cdot)$, and carries the corresponding optimal solution $\check a_{t}$ (also known as trial points) over to the subsequent stage (if any).
    \item \textit{Backward} step: which goes backward in time starting from $t=T$, and places cuts on $\check\Q^{\xi_{t-1}}_{t}(\cdot)$ which are obtained by the subgradient of $\partial\check\Q^{\xi_{t-1}}_{t}(\cdot)$ at the points collected during the forward step. In particular, given a sequence of trial points $(\check a_1, \dots, \check a_T)$, one places a cut on $\check \Q^{\xi_{t-1}}_t(\cdot)$ given by $q_t(a_{t-1}) := \beta_t^{\top} a_{t-1} + \pi_t$, where
        \begin{equation}
        \label{eq:cut_subgradient}
            \beta_t \in \E[\partial{\check{Q}}_{t}(\check a_{t-1},\xi_{t})],\; \forall \;t=T, \dots, 2
        \end{equation}
    and
        \begin{equation}
        \label{eq:cut_coeff}
            \pi_t :=
            \begin{cases} 
                \E[b^{\xi_t\top}_t \check \lambda_t] & \text{for}\; t=T \\
                \E\left[b^{\xi_T\top}_t \check \lambda_t+\sum_{c\in C^{\xi_t}_{t+1}} \pi_{t+1, c} \check \rho_{t, c}\right] & \text{for}\; t=T-1, T-2, \dots, 2
            \end{cases}
        \end{equation}
    where
    \begin{itemize}
        \item $C^{\xi_t}_{t+1}$ is a collection of all the cutting planes placed on $\check\Q^{\xi_t}_{t+1}(\cdot),\; \forall \;t=1, \dots, T-1$,
        \item $\check \rho_{t, c}$ are the optimal dual vectors associated with $\beta_{t+1, c}^{\top} x_{t} + \pi_{t+1, c} \leq \check \Q^{\xi_t}_{t+1}\; \forall \;\; c \in C_{t+1}$ and $t=2, \dots, T-1$, 
        \item $\check \lambda_t$ are the optimal dual vectors associated with constraints in $\A_t({a_{t-1}},\xi_{t}),\; \forall \;t=2, \dots, T$.
    \end{itemize}
\end{enumerate}
Due to the polyhedral structure of the epigraph of $\check \Q^{\xi_t}_{t+1}(\cdot)$, this iterative procedure builds a piece-wise linear lower approximation to $\check \Q^{\xi_t}_{t+1}(\cdot)$ that converges in a finite number of iterations with probability one~\cite{shapiro2011analysis}.

We care to note that, the cuts that one seeks to collect by the machinery of the nested Benders decomposition are generally categorized into two types of cuts: \textit{optimality} cuts and \textit{feasibility} cuts. In particular, in the forward step, when solving ${\check{Q}}_{t}(\check a_{t-1},\xi_{t})$ given $\check a_{t-1}$ and $\xi_t$: if the feasible set $\A_t(\check a_{t-1},\xi_t)=\emptyset$, a {feasibility} cut must be added; otherwise, the cut added to $\check \Q^{\xi_{t-1}}_{t}(\cdot)$ is referred to as an {optimality} cut. In the \ac{hdrlp} problem of our interest, we have what is known as the \textit{relatively complete recourse} property which qualifies the following assumption.
\begin{assm}
\label{assm:relatively_complete_recourse}
We assume that for any given value $\check a_{t-1}$ and a realization of the random process $\xi_t \in \Xi_t$, $\A_t(\check a_{t-1},\xi_t) \neq \emptyset$ and there always exists a feasible action $\check a_t= (\check u_t,\check v_t)$ such that $A^{\xi_t}_{t}{\check u_{t}} + B^{\xi_t}_t{\check u_{t-1}} + C^{\xi_t}_t{\check v_t} = b^{\xi_t}_t$.
\end{assm}
Assumption~\ref{assm:relatively_complete_recourse} eliminates the need to consider adding feasibility cuts. In Algorithm~\ref{alg:nestedBenders_MC}, we provide a pseudocode description of the nested Benders decomposition algorithm. 

\begin{algorithm}[htb]
\footnotesize
\caption{Nested Benders Decomposition.}
\label{alg:nestedBenders_MC}
\begin{algorithmic}
\State \textit{Input:} $a_0$, $\xi_1$, $M \in \mathbb{R}$ which is a \ac{lb} on $\Q^{\xi_t}_{t+1}, \forall \;\xi_t \in \Xi_t, \;t=1, \dots, T-1$, and termination criterion (e.g., time-limit).
\vspace{0.2cm}
\State \textbf{STEP 0:} (\textit{Initialization}). 
\begin{itemize}
    \item Set $C^{\xi_{t}}_{t+1} = \emptyset,\; \forall \; \xi_t \in \Xi_t, \; t=1, \dots, T-1$, $n=0$ be the iteration count.
    \item Let $Q^{n}_t(\cdot,\cdot), \forall \; t=1, \dots, T$, be the problem solved at the $n$th iteration and is defined according~\eqref{eq:FOSDDP}.
    \item Let $\theta^{\xi_t,n}_{t+1}$ represent the epigraph of $\Q^{\xi_t}_{t+1}(\cdot)$, and add constraint $\theta^{\xi_t,n}_{t+1} \geq M$ to $Q^{n}_t(\cdot,\cdot), \forall \xi_t \in \Xi_t, t=1, \dots, T-1$.
    \item Initialize $t=1$, $\xi^n_t=\xi_1$, $\check a^n_{t-1} = a_0$, $(\check a^n_1, \dots, \check a^n_T)$ for the trial points, and a sample path $(\xi^n_1, \dots, \xi^n_T)$.
\end{itemize}
\vspace{0.1cm}
\State \textbf{STEP 1:} (\textit{Forward Pass}).
\begin{itemize}
    \item Solve the current problem $Q^{n}_t(\check a^n_{t-1},\xi^n_t)$ to obtain $a^*_t$.
    \item Set $\check a^n_t = a^*_t$ and let $t \gets t+1$.
    \item If $t<T$, sample a new state $\xi^n_{t}$ given $\xi^n_{t-1}$ and repeat \textbf{STEP 1}. Otherwise, go to \textbf{STEP 2}.
\end{itemize}
\vspace{0.1cm}
\State \textbf{STEP 2:} (\textit{Termination check}). 
\begin{itemize}
    \item If termination criterion is met, \textit{STOP}; otherwise, go to \textbf{STEP 3}.
\end{itemize}
\State \textbf{STEP 3:} (\textit{Backward Pass}).
\begin{itemize}
    \item For every $\xi_t \in \Xi_t$ solve the current problem $Q^{n}_t(\check a^n_{t-1},\xi)$.
    \item Generate a cut (if any) to $\theta^{\xi^n_{t-1},n}_{t}$ following the description in~\eqref{eq:cut_subgradient} and~\eqref{eq:cut_coeff}.
    \item Add the cut (if any) to $C^{\xi_{t-1}}_{t-1}$ and let $t \gets t-1$.
    \item If $t = 1$, let $n \gets n+1$, set $\check a^n_{t-1} = a_0$, $\xi^n_t=\xi_1$ and go to \textbf{STEP 1}. Otherwise, repeat \textbf{STEP 3}.
\end{itemize}
\vspace{0.2cm}
\State \textit{Output:} The collection of cutting planes $\left( \{C^{\xi_{1}}_{2}\}, \{C^{\xi_{2}}_{3}\}_{\xi_2 \in \Xi_2}, \dots, \{C^{\xi_{T-1}}_{T}\}_{\xi_{T-1} \in \Xi_{T-1}}\right)$.
\end{algorithmic}
\end{algorithm}

The nested Benders decomposition algorithm not only provides an optimal first-stage solution to be implemented here-and-now, but also induces a decision policy for the \ac{famsp} model. Contrary to the \textit{online} \ac{rh} policy, discussed in Subsection~\ref{subsubsec:alternative_approx_policies_T}, policies induced by the expected cost-to-go functions $\Q^{\xi_t}_{t+1}(\cdot)$ are referred to as \textit{offline} policies. As the name suggests, {offline} policies are constructed in advance, before the unraveling of the stochastic process, whereas the {online} policies are constructed during the sequential realization of the stochastic process.
\section{RH-2SSP approach with random landfall time $T$} 
\label{sec:rh-2ssp-random}
In this section we describe the RH-2SSP approach with random landfall time $T$. In the deterministic landfall time case, the \ac{s2sspd} problem we solve in every roll is almost identical between different periods. However, because the planning horizon recedes as we roll forward, the aggregated number of periods shrinks by one period from one roll to the next. Note that the first-stage problem~\eqref{eq:twostage_detT_1st} does not contain any of the decision variables pertaining to the aftermath of the hurricane landfall, i.e., the delivery, shortage and demand for relief items. However, recall that in the random landfall time case, the first-stage problem in a given roll could be the stage when the hurricane makes landfall and the demand for relief items occurs. Therefore, in addition to reducing the length of the planning horizon in each roll, we move decision variables $y, \underline{e}, \overline{e}$ of the first period in the current roll from the second-stage problem to the first-stage problem. That is, we move $y_{t_{\text{roll}}}, \underline{e}_{t_{\text{roll}}}, \overline{e}_{t_{\text{roll}}}$ to the first-stage problem, where $t_{\text{roll}}$ is the time index of the first period in the current roll. We then reformulate~\eqref{eq:2SSP_rand_1st-stage} and~\eqref{eq:2SSP_rand_2nd-stage} as~\eqref{eq:2SSP_rand_1st-stage_RH} and~\eqref{eq:2SSP_rand_2nd-stage_RH}, as follows.

\textit{First-stage:}
\begin{subequations}
    \label{eq:2SSP_rand_1st-stage_RH}
    \begin{align}
    \displaystyle \underset{x, f, y_{t_{\text{roll}}}, \underline{e}_{t_{\text{roll}}}, \overline{e}_{t_{\text{roll}}}}{\min} \quad  & \displaystyle \sum_{t=t_{\text{roll}}}^{T_{\max}} \left(\sum_{i\in \{0\}\cup I}\sum_{i'\in I}c^b_{ii',t}f_{ii',t}+\sum_{i\in I}c^h_{i,t}x_{i,t}+ h_t\sum_{i\in I}f_{0i,t}\right) & \nonumber \\
    & \displaystyle - \left(\sum_{i\in I}\sum_{j\in J}c^a_{ij,t_{\text{roll}}}y_{ij,t_{\text{roll}}} + p\sum_{j \in J} \underline e_{j,t_{\text{roll}}} + q\sum_{i \in I} \overline e_{i,t_{\text{roll}}}\right) + \Q^{\xi_1}(x, f) & \nonumber \\
    \text{s.t.} \quad & \displaystyle \sum_{i'\in I, i'\neq i}f_{ii',t} \leq x_{i,t-1}, \hspace{6.2cm} \forall \;i\in I,\; t=t_{\text{roll}}, \dots, T_{\max} & \label{eq:2SSP_rand_1st-stage-1_RH} \\
     & \displaystyle x_{i,t} \leq u_i, \hspace{8.15cm} \forall \;i\in I,\; t=t_{\text{roll}}, \dots, T_{\max} & \label{eq:2SSP_rand_1st-stage-2_RH} \\
& \displaystyle \sum_{j\in J}y_{ij,t_{\text{roll}}} + \overline e_{i,t_{\text{roll}}} =  x_{i,t_{\text{roll}}-1} + \sum_{j\in \{0\}\cup I, j\neq i} f_{ji,t_{\text{roll}}} - \sum_{j\in I, j\neq i} f_{ij,t_{\text{roll}}} -  x_{i,t_{\text{roll}}}, \hspace{1.75cm} \forall \;i\in I, & \label{eq:2SSP_rand_1st-stage-3_RH} \\
      & \displaystyle \sum_{i \in I}y_{ij,t_{\text{roll}}} + \underline e_{j,t_{\text{roll}}}\geq d^{\xi_t}_{j,t_{\text{roll}}}, \hspace{7.75cm} \forall \;j\in J, & \label{eq:2SSP_rand_1st-stage-4_RH} \\
     & \displaystyle x_t, f_t\geq 0, \hspace{8cm} \forall \;t=t_{\text{roll}}, \dots, T_{\max} & \label{eq:2SSP_rand_1st-stage-5_RH}\\
     & \displaystyle y_{t_{\text{roll}}}, \underline{e}_{t_{\text{roll}}}, \overline{e}_{t_{\text{roll}}}\geq 0. &  \label{eq:2SSP_rand_1st-stage-6_RH}
     \end{align}
\end{subequations}
where $\Q^{\xi_1}(x, f) := \sum_{\xi_{T_{\max}} \in \Xi_{T_{\max}}} \mathbf{P}^{T_{\max}-1}_{\xi_{T_{\max}} \mid \xi_1} Q(x, f, \xi_{[t_{\text{roll}}:T_{\max}]})$, and $\xi_{[t_{\text{roll}}:T_{\max}]}$ gives the trajectory of the stochastic process between $t=t_{\text{roll}}$ to $t=T_{\max}$.

\textit{Second-stage:}
\begin{subequations}
\label{eq:2SSP_rand_2nd-stage_RH}
    \begin{align}
        \displaystyle Q(x, f, \xi_{[t_{\text{roll}}:T_{\max}]}) := \displaystyle \underset{y, \underline{e}, \overline{e}}{\min} \quad & \displaystyle \sum_{t=t_{\text{roll}}+1}^{T_{\max}} \left(\sum_{i\in I}\sum_{j\in J}c^a_{ij,t}y_{ij,t} + p\sum_{j \in J} \underline e_{j,t} + q\sum_{i \in I} \overline e_{i,t}\right) & \nonumber \\
        - & \displaystyle \sum_{t=t_{\text{roll}}+1}^{T_{\max}}\left(\sum_{i\in \{0\}\cup I}\sum_{i'\in I}c^b_{ii',t}{f}_{ii',t}+\sum_{i\in I}c^h_{i,t}{x}_{i,t}+ h_t\sum_{i\in I}{f}_{0i,t}\right)\mathbbm{1}_{\{t>T^{\xi_{[t_{\text{roll}}:T_{\max}]}}\}} & \nonumber \\
        \text{s.t.} \quad  & \displaystyle \sum_{j\in J}y_{ij,t} + \overline e_{i,t} =  x_{i,t-1} + \sum_{j\in \{0\}\cup I, j\neq i} f_{ji,t} - \sum_{j\in I, j\neq i} f_{ij,t} -  x_{i,t}, \forall \;i\in I,\; t=t_{\text{roll}}+1, \dots, T_{\max} & \label{eq:2SSP_rand_2nd-stage-1_RH} \\
        & \displaystyle \sum_{i \in I}y_{ij,t} + \underline e_{j,t}\geq d^{\xi_t}_{j,t}, \hspace{4.25cm} \forall \;j\in J,\; t=t_{\text{roll}}+1, \dots, T_{\max} & \label{eq:2SSP_rand_2nd-stage-2_RH} \\
        & \displaystyle y_t, \underline{e}_t , \overline{e}_t \geq 0, \hspace{6.15cm} \forall \; t=t_{\text{roll}}+1, \dots, T_{\max}. &  \label{eq:2SSP_rand_2nd-stage-3_RH} 
    \end{align}
\end{subequations}

\section{Problem Data} 
\label{sec:problem_data}
This section is organized as follows. First, we describe how we generate the configuration of the disaster relief logistics network. Second, we describe how the logistics cost parameters in the multi-period network flow model~\eqref{eq:det_form} are specified. Third, we discuss how the stochastic process on the evolution of the hurricane's attributes is modeled. Fourth, we demonstrate how we define the appropriate $T_{\max}$ parameter for the case when the hurricane's landfall time is random. Finally, we discuss how the hurricane's attributes are mapped to the demand for relief items at different \ac{dps}. We note that many of the parameter settings follow those given in~\cite{galindo2013prepositioning,pacheco2016forecast}.

\paragraph{The coordinates of the \ac{mdc}, \ac{sps} and \ac{dps}.} We consider an \ac{apa} similar to the one depicted in Figure~\ref{fig:network} when $T$ deterministic, and the one depicted in Figure~\ref{fig:network_randomT} when $T$ is random. In both cases, we randomly generate the DP locations by letting the $x$-coordinates of the \ac{sps} follow a uniform distribution $U(0,700)$, and independently, letting the $y$-coordinates follow a uniform distribution $U(0,100)$. Similarly, the $x$-coordinates of the \ac{dps} locations are randomly generated according to a uniform distribution $U(0,700)$, and independently, their $y$-coordinates are randomly generated according to a uniform distribution $U(100,200)$. Using these distributions, we generate the coordinates for different number of \ac{dps} $|J| \in \{10, 20, 30\}$ and \ac{sps} $|I| \in \{3, 6, 9\}$ to create network instances with different sizes. In all instances, we set the coordinates of the \ac{mdc} to be $(350,450)$.

\paragraph{Defining the logistics operational costs.} We assume that any transportation cost is a multiple of a unit cost $\omega = 0.0038$, scaled by an Euclidean distance-based metric. The unit cost of transporting a relief item from node $i \in \{0\}\cup I$ to node $i' \in I$ in period $t$ is given by: 
\begin{equation}
    \label{eq:ship_cost}
    c^b_{ii',t} = \omega(1+\nu(t-1))||(i_x, i_y) - (i'_x, i'_y)||_2,\; \quad \forall i \in \{0\}\cup I, i' \in I, t=1,2, \cdots, T_{\max},
\end{equation}
where $\nu$ is  the scaling factor that captures the increasing operational cost over time as we discussed in Assumption 5, and ($i_x,i_y$) and ($i'_x, i'_y$) are the coordinates of SP $i \in \{0\} \cup I$ and SP $i' \in I$, respectively. Similarly, the unit cost of transporting a relief item from SP $i \in \{0\}\cup I$ to DP $j \in J$ in period $t$ is given by: 
\begin{equation}
    \label{eq:deliver_cost}
    c_{ij,t}^a =  \omega(1+\nu(t-1))||(i_x, i_y) - (j_x, j_y)||_2,\; \quad \forall i \in I, j \in J, t=1, 2, \cdots, T_{\max},
\end{equation}
where ($i_x,i_y$) and ($j_x,j_y$) are the coordinates of SP $i \in I$ and DP $j \in J$, respectively. Moreover, we define the remaining parameters $h, c_t^h, p$ and $q$ as a multiple of a base cost $\beta = 5$: $c_{i,t}^h = 0.2\beta, \forall \; i \in I, t=1,2, \cdots, T_{\max}$, $h_t= \beta (1+\nu(t-1))\; \forall t=1, \dots, T_{\max}$, $p = 80\beta$, and $q = -0.05\beta$. 
Finally, the capacities of the \ac{sps} are generated randomly according to a uniform distribution $U(0.05\bar{d}\frac{|J|}{|I|}, 0.5\bar{d}\frac{|J|}{|I|})$, where $\bar d$ is the total maximum possible demand for relief items at any DP $j \in J$, and in our test instances we set $\bar d = 400$. The rational behind this setup is that the capacity of each SP $i \in I$ is (i) proportional to the maximum possible demand that can occur at a given DP; (ii) increasing in the total number of \ac{dps} in the \ac{apa} (i.e., $|J|$); and (iii) decreasing in the total number of \ac{sps} in the \ac{apa} that could supply different \ac{dps} with the relief items (i.e., $|I|$).

\paragraph{Modeling the evolution of the hurricane's attributes.} We first specify the state space $\Xi_t$ of the \ac{mc} model for the case of \textit{deterministic} landfall time. Recall that when $T$ is deterministic and known a priori, we only consider the hurricane's horizontal location and define $\Xi_t = A\times L_x$, where $A$ and $L_x$ are the state spaces associated with the hurricane's intensity and location along the $x$-coordinate, respectively. In our test instances, we set $A = \{0, 1, \dots, 5\}$, where $0$ is an absorbing state corresponding to the situation that the hurricane dissipates, and states $1$ to $5$ correspond to the five categories pertaining to the \ac{sshws}. The one-step transition probability matrix for the \ac{mc} associated with the intensity level $A$ is given by:
\begin{equation}
\label{eq:intensity_MC}
\mathbf{P}^{\alpha} = 
\begin{bmatrix}
1.00 & 0.00 & 0.00 & 0.00 & 0.00 & 0.00 \\ 
0.11 & 0.83 & 0.06 & 0.00 & 0.00 & 0.00 \\ 
0.00 & 0.15 & 0.60 & 0.25 & 0.00 & 0.00 \\ 
0.00 & 0.00 & 0.04 & 0.68 & 0.28 & 0.00 \\ 
0.00 & 0.00 & 0.00 & 0.18 & 0.79 & 0.03 \\ 
0.00 & 0.00 & 0.00 & 0.00 & 0.50 & 0.50
\end{bmatrix}.
\end{equation}
Additionally, we define the state space for the hurricane's landfall location along the $x$-coordinate as: 
\[
L_{x} = \{[0,100), [100,200), \dots, [600,700)\},
\]
and the associated one-step transition probability matrix is given by: 
\begin{equation}
\label{eq:hlocation_MC}
\mathbf{P}^{\ell_{x}} = 
\begin{bmatrix}
0.004 & 0.300 & 0.395 & 0.198 & 0.049 & 0.038 & 0.016 \\ 
0.150 & 0.202 & 0.249 & 0.222 & 0.117 & 0.033 & 0.027 \\ 
0.198 & 0.249 & 0.029 & 0.168 & 0.206 & 0.099 & 0.051 \\ 
0.099 & 0.222 & 0.169 & 0.012 & 0.150 & 0.198 & 0.150 \\ 
0.025 & 0.117 & 0.206 & 0.15 & 0.004 & 0.150 & 0.348 \\ 
0.019 & 0.033 & 0.098 & 0.198 & 0.150 & 0.004 & 0.498 \\ 
0.008 & 0.019 & 0.025 & 0.098 & 0.198 & 0.150 & 0.502
\end{bmatrix}.
\end{equation}
Consequently, when $T$ is deterministic, the stochastic process associated with the hurricane's evolution is characterized by an \ac{mc}, where the transition probability matrix $\mathbf{P}$ is given by $\mathbf{P}_{\xi_{t} \mid \xi_{t-1}} = \mathbf{P}_{\alpha_{t} \mid \alpha_{t-1}} \times \mathbf{P}_{\ell_{x,t} \mid \ell_{x,t-1}}$ with $\mathbf{P}^{\alpha}$ and $\mathbf{P}^{\ell_{x}}$ specified above.

We extended this to the case where $T$ is \textit{random}, as we discussed in Section~\ref{subsec:randomT_MSP}, by introducing a temporal state that corresponds to the $y$-coordinate $\ell_{y,t} \in L_{y}$, where $L_y$ is defined as:
\[
L_{y} = \left\{[-350,-300), [-300,-250), \dots, [-50,0), [0,\infty)\right\},
\]
and the associated one-step transition probability matrix is given by:
\begin{equation}
\label{eq:vlocation_MC}
\mathbf{P}^{\ell_{y}} = 
\begin{bmatrix}
0.0 & 0.6 & 0.3 & 0.1 & 0.0 & 0.0 & 0.0 & 0.0 \\ 
0.0 & 0.0 & 0.6 & 0.3 & 0.1 & 0.0 & 0.0 & 0.0 \\ 
0.0 & 0.0 & 0.0 & 0.6 & 0.3 & 0.1 & 0.0 & 0.0 \\ 
0.0 & 0.0 & 0.0 & 0.0 & 0.6 & 0.3 & 0.1 & 0.0 \\ 
0.0 & 0.0 & 0.0 & 0.0 & 0.0 & 0.6 & 0.3 & 0.1 \\ 
0.0 & 0.0 & 0.0 & 0.0 & 0.0 & 0.0 & 1.0 & 0.0 \\ 
0.0 & 0.0 & 0.0 & 0.0 & 0.0 & 0.0 & 0.0 & 1.0 \\ 
0.0 & 0.0 & 0.0 & 0.0 & 0.0 & 0.0 & 0.0 & 1.0 \\
\end{bmatrix}.
\end{equation}
Consequently, when $T$ is random, the stochastic process associated with the hurricane's evolution is characterized by an \ac{mc} whose the transition probability matrix $\mathbf{P}$ is given by $\mathbf{P}_{\xi_{t+1} \mid \xi_t} := \mathbf{P}_{\alpha_{t+1} \mid \alpha_t} \times \mathbf{P}_{\ell_{x,t+1} \mid \ell_{x,t}} \times \mathbf{P}_{\ell_{y,t+1} \mid \ell_{y,t}}$, with matrices $\mathbf{P}^{\alpha}$, $\mathbf{P}^{\ell_{x}}$ and $\mathbf{P}^{\ell_{y}}$ specified above. Note that state $\ell_{y,t}=[-50,0)$ corresponds to the landfall event and state $\ell_{y,t} = [0,\infty)$ is an absorbing state indicating that the hurricane has already made landfall. The rational behind this setup is to have an \ac{mc} where the probability that the hurricane goes backward or remains in the same place from one period to the next is zero, i.e., $\mathbbm{P}(\ell^{-}_{y,t+1} \leq \ell^{-}_{y,t} | \ell_{y,t}) = 0$, where $\ell^{-}_{y,t}$ indicates the \ac{lb} of the interval associated with $\ell_{y,t}$. Additionally, we care to note that a dissipating hurricane corresponds to the state $\alpha_t = 0$. Accordingly, we have that:
\begin{itemize}
    \item the set of transient states $\mathcal{T}:= \{\xi_t = (\alpha_t, \ell_{x,t}, \ell_{y,t}) \; | \; \alpha_t \neq 0, \ell^{-}_{y,t} < 0\}$, and
    \item the set of absorbing states $\mathcal{A}:= \{\xi_t = (\alpha_t, \ell_{x,t}, \ell_{y,t}) \; | \; \alpha_t = 0$ or $ \ell^{-}_{y,t} \geq 0\}$.
\end{itemize}

\paragraph{Defining $T_{\max}$ for the case of random landfall time.} To define the value of $T_{\max}$, we need to identify the largest number of steps required to transition from any initial state $\ell_{y,1}$ to state $[0,\infty)$. Given the one-step transition probability matrix $\mathbf{P}^{\ell_{y}}$ specified in in~\eqref{eq:vlocation_MC}, since $\mathbbm{P}(\ell^{-}_{y,t+1} \leq \ell^{-}_{y,t} | \ell_{y,t}) = 0$, the largest number of steps required to transition from state $\ell_{y,1}$ to state $[0,\infty)$ corresponds to the situation where $\ell^{-}_{y,t+1} = \ell^{+}_{y,t}, \forall t=1, \dots, T_{\max}$, where $\ell^{+}_{y,t}$ indicate the \ac{ub} of the interval associated with $\ell_{y,t}$. In other words, in each period, the \ac{mc} visits every subsequent intervals without skipping any of them. Assuming that $\ell_{y,1}=[-350,-300)$, we can set $T_{\max} = |\mathrm{L}_{y}|= 8$. Moreover, we define the value of $T$ used in our test instances for the case where the time of landfall is deterministic, from the one-step transition probability matrix $\mathbf{P}^{\ell_{y}}$ as well. We do this by assuming that $T = \floor{\hat T} = 5$, where $\hat T$ is the \textit{expected} number of steps before reaching state $\ell_{y,t} = [0,\infty)$, starting from initial state $\ell_{y,1}=[-350,-300)$. \exclude{Note that $\hat T$ is obtained by rearranging the components of $\mathbf{P}_{\ell_{y,t} \mid \ell_{y,t-1}}$ into a new matrix
\begin{equation*}
\mathbf{P}'_{\ell_{y,t} \mid \ell_{y,t-1}} = 
\begin{bmatrix}
P_{\mathcal{T}} & P_{\mathcal{A}}\\
\mathbf{0}      & I_{\mathcal{A}}\\ 
\end{bmatrix},
\end{equation*}
where $P_{\mathcal{T}} \in \mathbb{R}^{n\times n}$ is a submatrix of~\eqref{eq:vlocation_MC} that describes the probability of transitioning from some transient state $\ell_{y,t} \in L_{y}$ to another $\ell'_{y,t} \in L_{y}$; $P_{\mathcal{A}}\in \mathbb{R}^{n\times m}$ is a submatrix of~\eqref{eq:vlocation_MC} that describes the probability of transitioning from some transient state $\ell_{y,t} \in L_{y}$ to some absorbing state $\ell'_{y,t} \in L_{y}$;  $\mathbf{0} \in \mathbb{R}^{m \times n}$ zero matrix; and $I_{\mathcal{A}} \in \mathbb{R}^{m\times m}$ identity matrix. We can then let $\hat T$ be the first entry of the vector $(I_{\mathcal{A}} - P_{\mathcal{T}})^{-1} \mathbf{1}$, which gives the expected number of steps before being absorbed, when starting in transient state $k$, which is the $k$th entry of the vector, and $\mathbf{1}$ is a length-$n$ column vector, whose entries are all ones~\cite{ross2014introduction}.}

\paragraph{Modeling the demand at \ac{dps}.} To define the demand value $d_{j,t}$ associated with the hurricane's state $\xi_t \in \Xi_t$ at time $t$ for DP $j \in J$ with location $(j_x, j_y)$, we do the following. First, we discretize each of the possible landfall locations $\ell_{x,t}$ \textit{uniformly} into a set of points $\{\ell^m_{x,t}\}_{m=1}^{M}$, where each point $\ell^m_{x,t}$ represents a realization of the hurricane's landfall location (again, we consider the $x$-coordinates only and the $y$-coordinates are $0$), $\forall \; m=1,2,\cdots,M$. Hence, for a given state $\ell_{x,t}$, the $m$-th realization for $x$-coordinate of the landfall location is given by:
\begin{equation}
    \label{eq:exact_coor}
    \ell^m_{x,t} := \ell^{-}_{x,t}+\frac{M}{2}+M(m-1) \equiv \ell^{+}_{x,t}-\frac{M}{2}+M(m-1),
\end{equation}
where $\ell^{-}_{x,t}$ and $\ell^{+}_{x,t}$ indicate the lower and upper ends of the interval associated with $\ell_{x,t}$, respectively. We further assume that each of the $M$ realizations is equally likely. In our experiments, we assume that $M=10$, and as an example, state $\ell_{x,t} = [0,100)$ is discretized into a set of realizations $\{5, 15, \dots, 95\}$. We then define the deterministic mapping $D(\alpha_T, \ell_T)$ as
\[
D(\alpha_T, \ell_T) = \bar d \times \left(1-\frac{\delta}{\bar \delta}\right) \times \frac{\alpha_T^2}{(|\mathrm{A}|-1)^2} 
\]
where $\delta$ gives the distance between a DP $j \in J$ and $\ell_T$, and $\bar \delta$ is the maximum possible distance to a DP where demand for relief items occurs. Hence, given the hurricane's attributes $\xi_t = (\alpha_t,\ell^{m}_{x,t}, \ell_{y,t})$, we assume that the demand associated with DP $j \in J$ with location $(j_x, j_y)$ is given by:
\begin{equation}
    \label{eq:demand_function}
    d^{\xi_t}_{j,t} = 
    \begin{cases} 
      \bar d \times \left(1-\frac{\delta(\ell^{m}_{x,t},j)}{\bar \delta}\right) \times \frac{\alpha_t^2}{(|\mathrm{A}|-1)^2}  & \text{if}\;  \delta(\ell^{m}_{x,t},j) \leq \bar \delta \; \text{and} \; \ell^{+}_{y,t} = 0 \\
      0 & \text{otherwise.}
   \end{cases}
\end{equation}
where $\delta(\ell^{m}_{x,t},j) := ||(\ell^{m}_{x,t},0) - (j_x, j_y)||_2$ gives the Euclidean distance between the hurricane's landfall location and the location of DP $j \in J$. In our test instances, we set $\bar \delta = 300$. As we can see from~\eqref{eq:demand_function}, in accordance with the discussion in Section~\ref{sec:MSP_models}, the demand function $d^{\xi_t}_{j,t}$ is \textit{increasing} in the hurricane's intensity $\alpha_t$ and \textit{decreasing} in the distance $\delta(\ell^{m}_{x,t},j)$ between the hurricane's landfall location and the location of a DP.

\paragraph{Initial values.} We assume that every SP is empty at the start of the planning horizon, that is, $x_{0} = 0$. Additionally, for our initial experiments, we assume that $\alpha_{1} = 1$ and $\ell_{1,x} = [100,200)$ for the case of deterministic landfall time, and assume in addition that $\ell_{1,y} = [-350,-300)$ for the case of random landfall time. \exclude{We investigate other initial states in our sensitivity analysis in Subsection~\ref{subsec:numerical_results}.}

We end this Section by noting that the above problem data is only considered for the purpose of our numerical experiments, and is by no means derived from historical data. As a summary, all parameters used in our numerical experiments are provided in Table~\ref{tab:data} in Appendix~\ref{sec:append_tables}.

\section{Implementation Details}
\label{sec:implement_details}
We organize the presentation of our implementation details as follows. First, we discuss the performance evaluation of these approaches. Next, we discuss how the stochastic optimization models are solved by: (i) describing the choice of initial values for some of the parameters; (ii) describing how the training sample paths are used for different stochastic optimization models; and (iii) discussing the termination criteria used for solving these models. 

\paragraph{Initial values.} For the \ac{famspd} , \ac{famspr} , and \ac{s2sspd} (including the \ac{s2sspd} solved in every roll of the \ac{rh2sspd} ), we assume that the initial \ac{lb} for the approximate expected cost-to-go function to be $0$. For the \ac{s2sspr} (including the one solved in every roll of the \ac{rh2sspr}), because of the reimbursement variable~\eqref{eq:reimbursement_var}, it is possible that we arrive at a negative recourse value, and thus we set the initial \ac{lb} to be $-10^{10}$. 

\paragraph{Sample paths for out-of-sample evaluation.} To evaluate the performances of policies obtained from different approaches, we create a set of out-of-sample scenarios (sample paths) with a sample size of $N = 1000$. We compute a sample mean
\begin{equation}
    \label{eq:sample_mean}
    \hat{z}=\frac{1}{N}\sum_{n=1}^N z(\hat{\xi}^n),
\end{equation}
and a sample standard deviation
\begin{equation}
    \label{eq:sample_std}
    \hat{\sigma}=\sqrt{\frac{1}{N-1}\sum_{n=1}^N \left(z(\hat{\xi}^n)-\hat{z}\right)^2},
\end{equation}
where $\hat{\xi}^n := (\hat{\xi}^n_2, \hat{\xi}^n_3, \cdots, \hat{\xi}^n_T)$ is the realized trajectory of the stochastic process corresponding to sample path $n$, and $z(\hat{\xi}^n)$ is the cumulative objective value accrued by the employed decision policy at sample path $n$. We compute a $95\%$ \ac{ci} for the performance of the corresponding policy on the set of out-of-sample scenarios, which is given by $\hat{z}_{\pm} := [\hat{z} - 1.96 \hat{\sigma}/\sqrt{N}, \hat{z} + 1.96 \hat{\sigma}/\sqrt{N}]$. 

\paragraph{Training scenarios.} We use nested Benders algorithm (see Algorithm~\ref{alg:nestedBenders_MC}) to obtain \textit{offline} decision policies associated with the \ac{famspd} and \ac{famspr} models. In principle, during the training process of the nested Benders algorithm, every iteration of Algorithm~\ref{alg:nestedBenders_MC} can be carried over a set of multiple sample paths. That is, in the forward step, one can collect trial points over multiple sample paths and, analogously, generate cuts for each of those sample paths when going backward. However, in all of our implementation, we observed that using a \textit{single} sample path per forward/backward step works the best. 

To solve the \ac{s2sspd} and \ac{s2sspr} models, we employ a standard implementation of the L-shaped method~\cite{birge2011introduction}. The basic idea of the L-shaped method is to successively build an outer-approximation of the recourse function $\Q^{\xi_1}(x_{T-1})$ in~\eqref{eq:twostage_detT_1st} (or $\Q^{\xi_1}(x, f)$ in~\eqref{eq:2SSP_rand_1st-stage}). To do this, we can generate cuts in the spirit of the one generated at the terminal stage, in the backward step of Algorithm~\ref{alg:nestedBenders_MC}. One crucial distinction here is that, the recourse function for at the terminal stage in Algorithm~\ref{alg:nestedBenders_MC} reflects the cost of the terminal period only, and to generate the cut, one needs to solve a scenario subproblem for each state $\xi_t \in \Xi_t$. However, the recourse function in~\eqref{eq:twostage_detT_1st} (or in~\eqref{eq:2SSP_rand_1st-stage}) reflects the aggregated cost of the periods between $t=2$ and $t=T$ (or $t=T_{\max}$ when $T$ is random). Hence, to generate a cut, one needs to solve the second-stage problem ~\eqref{eq:twostage_detT_2nd} (or in~\eqref{eq:2SSP_rand_2nd-stage}) for every possible sample path $(\xi_2, \dots, \xi_{T})$ (or $(\xi_2, \dots, \xi_{T_{\max}})$ when $T$ is random). According to our problem data, this would correspond to a maximum of $\left(|\mathrm{A}| \times |L_{x}|\right)^{T-1}$ (or $\left(|\mathrm{A}| \times |L_{x}| \times |L_{y}| \right)^{T_{\max}-1}$ when $T$ is random) sample paths. This exponentially large number of scenarios is especially prohibitive in the \ac{rh2sspd} and \ac{rh2sspr} approaches, where one solves the \ac{s2sspd} (or \ac{s2sspr}) model at each period $t=1, \dots, T-1$. We circumvent this by generating $100$ equally likely scenarios to serve as the training scenarios for the two-stage stochastic programming models. The same number of scenarios is also used in every roll of the \ac{rh2sspd}  and \ac{rh2sspr} approaches. We performed an in-sample stability test in our experiments that shows that the size of $100$ training samples is appropriate. 

\paragraph{Termination criteria.} We use the following {termination criteria} for the nested Benders decomposition algorithm: (i) \textit{the maximum number of iterations} is set to be $10^5$; (ii) the \textit{time limit} is set to be three hours; and (iii) a \textit{stability test} which keeps track of the progress of the \ac{lb} for the optimal objective value (denoted by $\underline{z}$): if $\underline{z}$ does not progress by more than $\epsilon = 10^{-5}$ (in the relative term) for more than $\bar{m} = 500$ consecutive iterations, i.e., $(\underline{z}^n- \underline{z}^{n-\bar{m}})/\underline{z}^n < \epsilon$, we terminate the algorithm. We use the same termination criteria when solving the two-stage SP model, both the for \ac{rh2sspd}  and \ac{rh2sspr} in every roll of the rolling-horizon procedure and for \ac{s2sspd} and \ac{s2sspr}. Although criteria (i) and (ii) are very unlikely to be reached for \ac{s2sspd} and \ac{s2sspr}, we maintain them for the sake of consistency. 
\section{Acronyms}
\printacronyms
\clearpage
\section{Figures}
\label{sec:append_figures}
\begin{figure}[htbp]
    \centering
    \subfigure[$\nu = 0.001$.]{
    \includegraphics[scale=0.25]{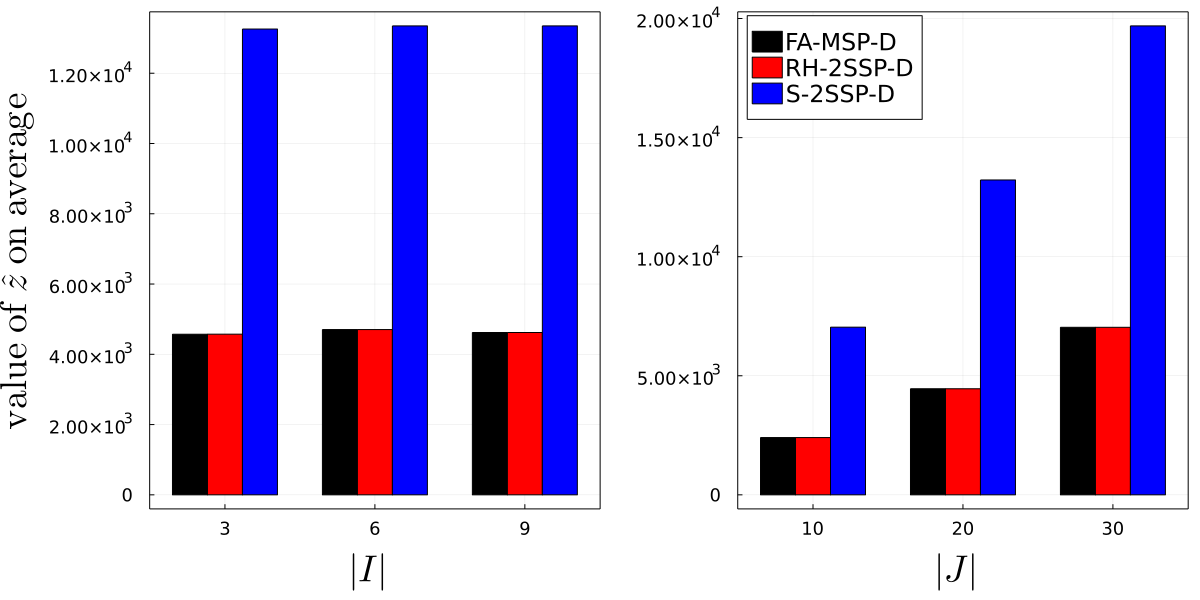}
    \label{fig:detT_avg_0.001}}
    \subfigure[$\nu = 0.60$.]{
    \includegraphics[scale=0.25]{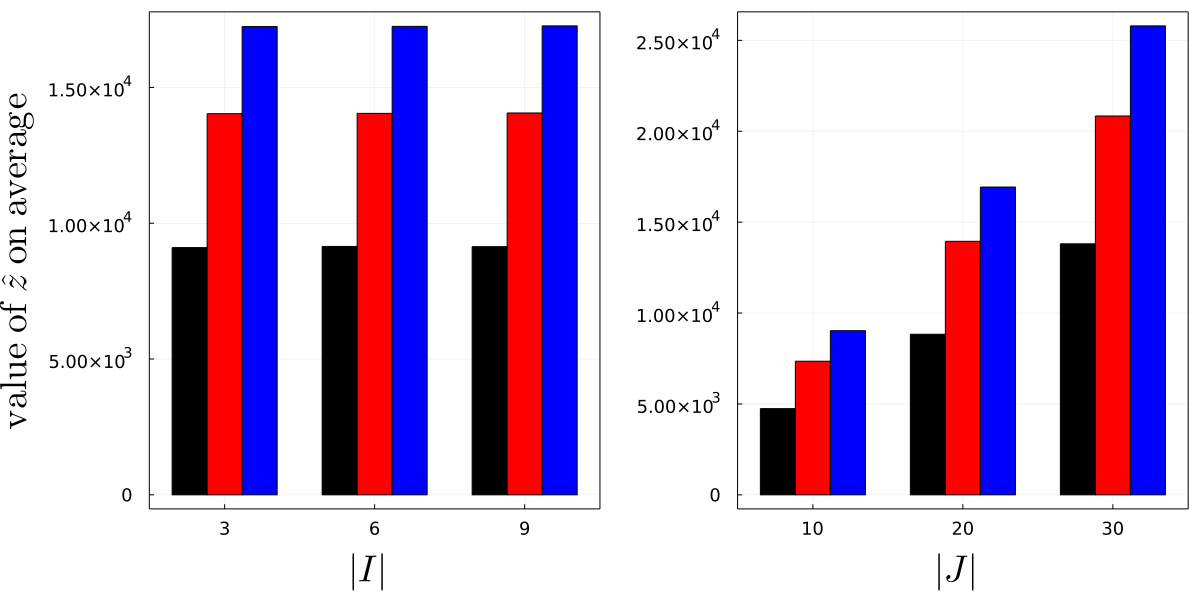}
    \label{fig:detT_avg_0.60}}
    \subfigure[$\nu = 5.00$.]{
    \includegraphics[scale=0.25]{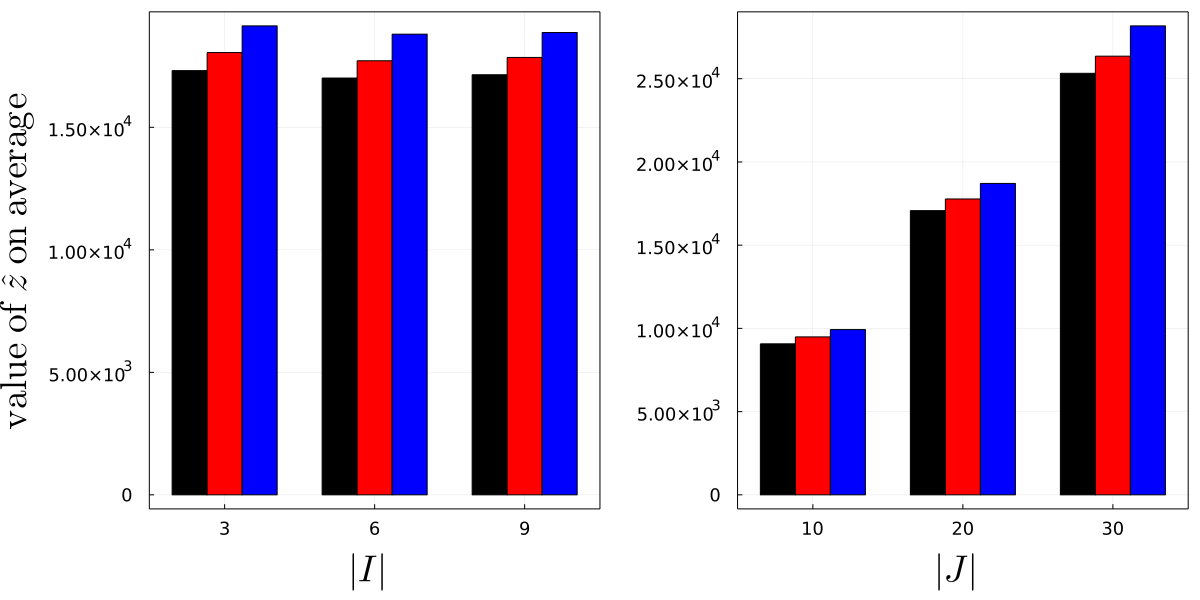}
    \label{fig:detT_avg_5.00}}
    \caption{ $\hat{z}$ values averaged across the different SPs and DPs for $\nu \in \{0.001, 0.6, 5.00\}$ when $T$ is \textit{deterministic}.}
    \label{fig:detT_avg}
\end{figure}
\begin{figure}[htbp]
    \centering
    \subfigure[$\nu = 0.001$.]{
    \includegraphics[scale=0.25]{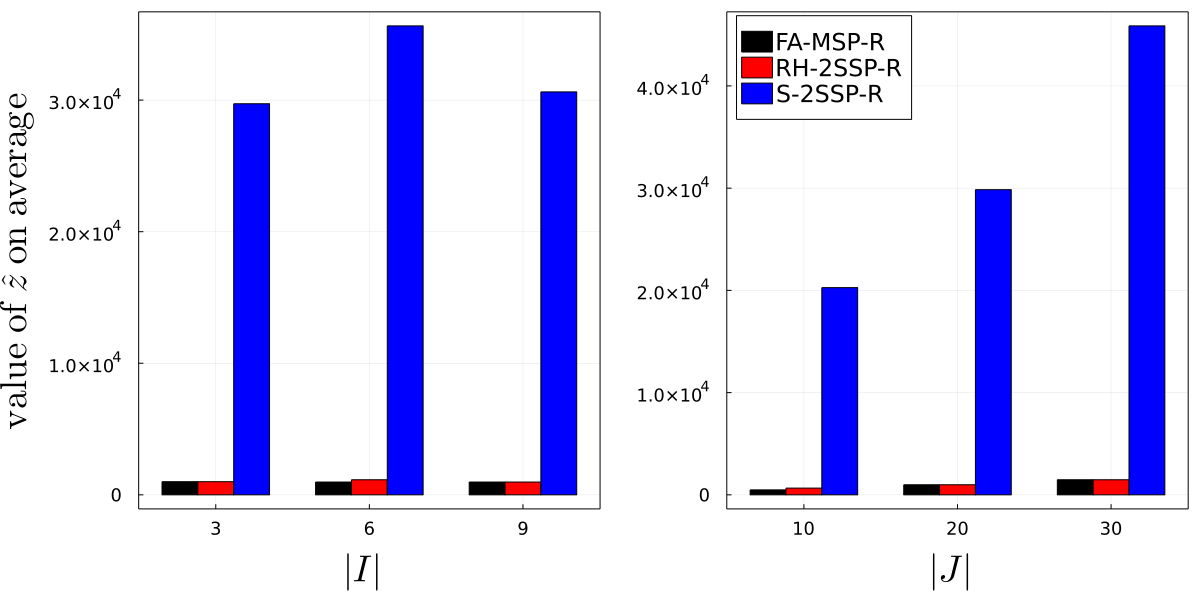}
    \label{fig:randT_avg_0.001}}
    \subfigure[$\nu = 0.60$.]{
    \includegraphics[scale=0.25]{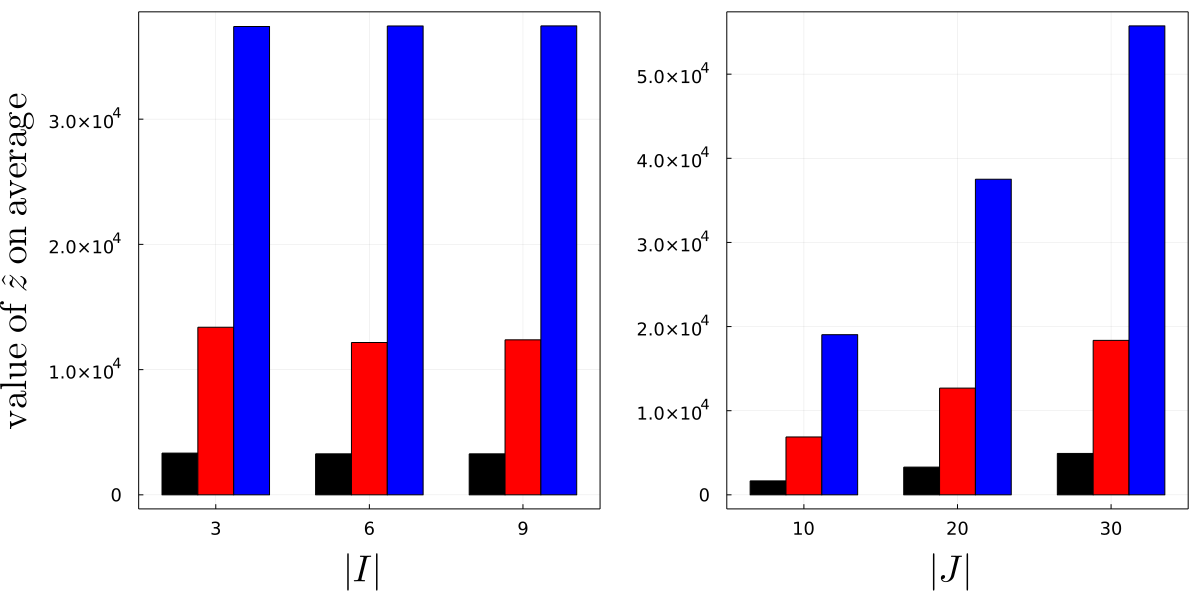}
    \label{fig:randT_avg_0.60}}
    \subfigure[$\nu = 5.00$.]{
    \includegraphics[scale=0.25]{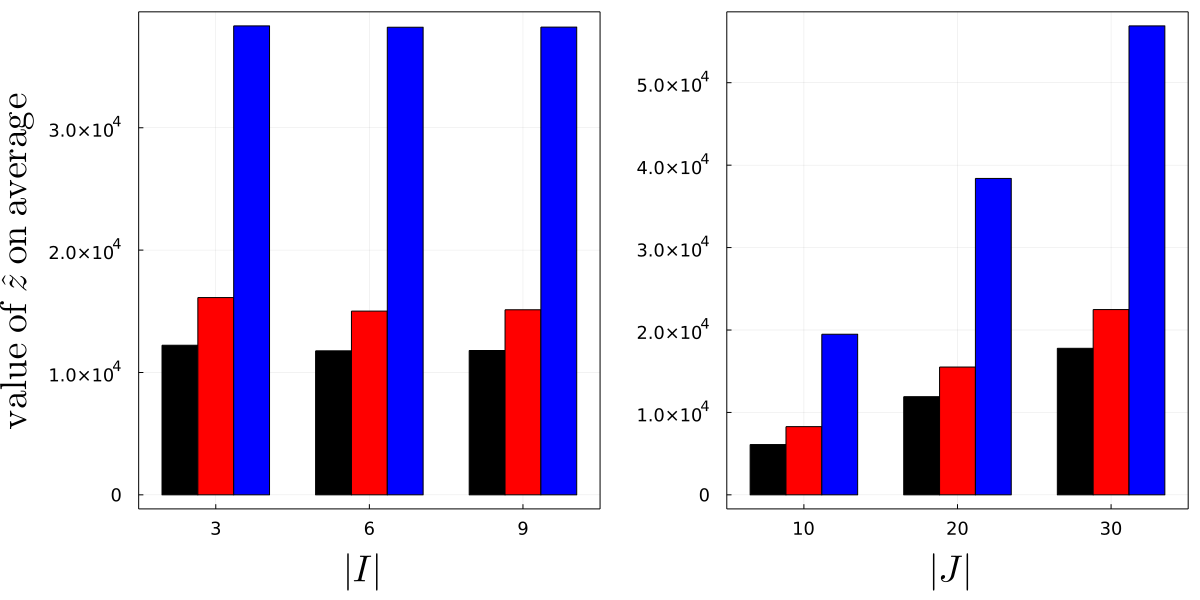}
    \label{fig:randT_avg_5.00}}
    \caption{ $\hat{z}$ values averaged across the different SPs and DPs for $\nu \in \{0.001, 0.6, 5.00\}$ when $T$ is \textit{random}.}
    \label{fig:randT_avg}
\end{figure}
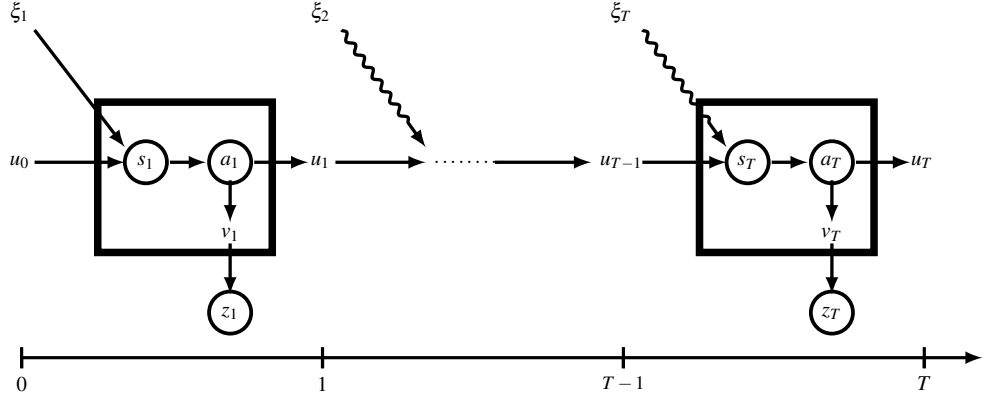
\begin{figure}[htbp]
\label{fig:MSP_models}
\centering
\begin{tikzpicture}[scale=0.8, transform shape]
\node [plain_circle_node] (u_0) at (0,0) {$u_0$};
\node [plain_circle_node] (xi_1) at (0,2.5) {$\xi_1$};
\node [solid_circle_node] (s_1) at (2.1,0) {$s_1$};
\node [solid_circle_node] (a_1) at (3.5,0) {$ a_1$};
\node [plain_circle_node] (v_1) at (3.5,-1.2) {$ v_1$};
\node [solid_circle_node] (z_1) at (3.5,-2.5) {$z_1$};
\node [plain_circle_node] (u_1) at (5,0) {$ u_1$};
\node [plain_circle_node] (xi_2) at (5,2.5) {$\xi_2$};
\node [plain_circle_node] (s_2) at (6.975,0) {};
\node [plain_circle_node] (u_T-2) at (7.9825,0) {};
\node [plain_circle_node] (u_T-1) at (10,0) {$ u_{T-1}$};
\node [solid_circle_node] (s_T) at (12.1,0) {$s_T$};
\node [plain_circle_node] (xi_T) at (10,2.5) {$\xi_T$};
\node [solid_circle_node] (a_T) at (13.5,0) {$ a_T$};
\node [plain_circle_node] (v_T) at (13.5,-1.2) {$v_T$};
\node [solid_circle_node] (z_T) at (13.5,-2.5) {$z_T$};
\node [plain_circle_node] (u_T) at (15,0) {$ u_T$};
\node [solid_square_node] (t1) at (2.75,-0.25) {};
\node [plain_square_node] (t2) at (7.5,0) {$\dots \dots \dots$};
\node [solid_square_node] (t_T) at (12.75,-0.25) {};

\draw [-latex,line width=1.25pt,color=black] (0.25,0) -- (1.75,0);
\draw [-latex,line width=1.25pt,color=black] (0.25,2.2) -- (1.75,0.25);
\draw [-latex,line width=1.25pt,color=black] (5.25,0) -- (6.75,0);
\draw [-latex,line width=1.25pt,color=black,decorate, decoration={snake,amplitude=.4mm,segment length=2mm,post length=3mm}] (5.35,2.2) -- (6.75,0.25);
\draw [-latex,line width=1.25pt,color=black] (10.35,0) -- (11.75,0);
\draw [-latex,line width=1.25pt,color=black,decorate, decoration={snake,amplitude=.4mm,segment length=2mm,post length=3mm}] (10.35,2.2) -- (11.75,0.25);
\draw [-latex,line width=1.25pt,color=black] (2.5,0) -- (3.1,0);
\draw [-latex,line width=1.25pt,color=black] (12.5,0) -- (13.1,0);
\draw [-latex,line width=1.25pt,color=black] (3.9,0) -- (4.8,0);
\draw [-latex,line width=1.25pt,color=black] (13.9,0) -- (14.8,0);
\draw [-latex,line width=1.25pt,color=black] (7.9,0) -- (9.5,0);

\draw [-latex,line width=1.25pt,color=black] (3.5,-0.3) -- (3.5,-1);
\draw [-latex,line width=1.25pt,color=black] (13.5,-0.3) -- (13.5,-1);

\draw [-latex,line width=1.25pt,color=black] (3.5,-1.35) -- (3.5,-2.2);
\draw [-latex,line width=1.25pt,color=black] (13.5,-1.35) -- (13.5,-2.2);
 
\draw[-latex, very thick] (0.01,-3.25) -- (16,-3.25);
\foreach \x in  {0,5,10,15}
\draw[shift={(\x,-3.25)},color=black, very thick] (1pt,5pt) -- (1pt,-5pt);
\draw[shift={(0,-3.25)},color=black, very thick] (1pt,0pt) -- (1pt,-5pt) node[below] {$0$};
\draw[shift={(5,-3.25)},color=black, very thick] (1pt,0pt) -- (1pt,-5pt) node[below] {$1$};
\draw[shift={(10,-3.25)},color=black, very thick] (1pt,0pt) -- (1pt,-5pt) node[below] {\footnotesize$T-1$};
\draw[shift={(15,-3.25)},color=black, very thick] (1pt,0pt) -- (1pt,-5pt) node[below] {$T$};
\end{tikzpicture}
\caption{A symbolic representation of the sequential decision-making process in MSP models.}
\end{figure}
\FloatBarrier
\section{Tables}
\label{sec:append_tables}
\begin{table}[htbp]
\centering
\begin{adjustbox}{width=\textwidth}
\begin{tabular}{@{}lllccccccccccccc@{}}
\toprule
 & & \multicolumn{4}{c}{$\alpha_1 = 1$} & \phantom{abc}& \multicolumn{4}{c}{$\alpha_1 = 3$} & \phantom{abc} & \multicolumn{4}{c}{$\alpha_1 = 5$}\\
\cmidrule{3-6} \cmidrule{8-11} \cmidrule{13-16}
  $\nu$ && $t=1$ & $t=2$ & $t=3$ & $t=4$ && $t=1$ & $t=2$ & $t=3$ & $t=4$ && $t=1$ & $t=2$ & $t=3$ & $t=4$\\ \midrule
0.1 &  & -   & -  & -  & 249.63 &  & -    & -   & -  & 1180.35 &  & -    & - & - & 1297.37 \\
0.2 &  & -   & 14.76 & 19.83 & 205.58 &  & 680.83  & 129.38 & 90.38 & 250.75  &  & 1173.99 & 3.63 & 8.44 & 69.70   \\
0.3 &  & -   & 20.77 & 26.14 & 178.30 &  & 1169.41 & 4.54   & 2.14  & 20.44   &  & 1210.36 & 0.09 & 0.82 & 52.14   \\
0.4 &  & -   & 22.66 & 26.74 & 167.85 &  & 1181.06 & 3.26   & 0.69  & 14.23   &  & 1217.58 & 0.01 & 0.19 & 48.59   \\
0.5 &  & 163.07 & 10.91 & 21.58 & 68.75  &  & 1197.35 & 0.37   & 0.55  & 8.61    &  & 1218.24 & - & 0.14 & 48.09   \\
0.6 &  & 175.08 & 28.68 & 16.46 & 54.86  &  & 1199.49 & 0.12   & 0.50  & 7.98    &  & 1218.24 & - & 0.14 & 47.41   \\
0.7 &  & 190.76 & 31.12 & 22.18 & 40.58  &  & 1199.49 & 1.30   & 0.35  & 7.45    &  & 1218.24 & - & - & 46.58   \\
0.8 &  & 195.77 & 32.90 & 22.26 & 33.08  &  & 1210.36 & 0.08   & 0.04  & 4.37    &  & 1218.24 & - & - & 46.58   \\
0.9 &  & 200.02 & 33.69 & 21.90 & 31.70  &  & 1210.57 & 0.08   & 0.03  & 4.32    &  & 1220.43 & - & - & 45.82   \\
1.0 &  & 204.91 & 34.15 & 21.74 & 30.66  &  & 1210.57 & -   & 0.02  & 3.21    &  & 1220.43 & - & - & 45.82    \\ 
\bottomrule
\end{tabular}
\end{adjustbox}
\caption{The average amount of procured relief items $\bar{f}_t$ (see~\eqref{eq:avg_procurement}) obtained from policies associated with the \ac{famspd} model as a function of the cost-scaling factor $\nu$ for different time periods $t=1, \dots, T$, and for different initial intensity levels $\alpha_1 \in \{1, 3, 5\}$.}
\label{tab:sens_results_det}
\end{table}
\begin{table}[htbp]
\centering
\begin{adjustbox}{width=\textwidth}
\begin{tabular}{@{}lllcccccccccccccccccccccc@{}}
\toprule
 & & \multicolumn{7}{c}{$\alpha_1 = 1$} & \phantom{abc}& \multicolumn{7}{c}{$\alpha_1 = 3$} & \phantom{abc} & \multicolumn{7}{c}{$\alpha_1 = 5$} \\ \cmidrule{3-9} \cmidrule{11-17} \cmidrule{19-25}
  $\nu$ && $t=1$ & $t=2$ & $t=3$ & $t=4$ & $t=5$ & $t=6$ & $t=7$ && $t=1$ & $t=2$ & $t=3$ & $t=4$ & $t=5$ & $t=6$ & $t=7$ && $t=1$ & $t=2$ & $t=3$ & $t=4$ & $t=5$ & $t=6$ & $t=7$ \\ \midrule
0.1 &  & -  & - & 0.74 & 9.02  & 23.60 & 29.47 & 7.38 &  & -   & - & 10.38 & 106.26 & 241.05 & 199.38 & 50.99 &  & -   & - & 13.66 & 120.03 & 325.65 & 212.89 & 57.82 \\
0.2 &  & -  & - & 0.96 & 9.75  & 26.33 & 27.25 & 6.00 &  & 61.41  & - & 14.21 & 117.85 & 220.31 & 162.25 & 36.83 &  & 113.12 & - & 17.29 & 127.04 & 273.65 & 164.05 & 43.04 \\
0.3 &  & -  & 2.08 & 3.46 & 10.07 & 25.64 & 25.87 & 4.54 &  & 331.89 & 1.86 & 7.22  & 60.69  & 129.83 & 100.87 & 24.96 &  & 441.09 & 0.37 & 8.92  & 64.30  & 156.63 & 93.97  & 23.87 \\
0.4 &  & -  & 3.24 & 5.14 & 10.63 & 24.80 & 24.79 & 4.22 &  & 492.66 & 0.61 & 4.64  & 38.97  & 82.51  & 69.42  & 18.73 &  & 586.68 & 0.26 & 5.87  & 42.45  & 112.55 & 67.19  & 16.69 \\
0.5 &  & 29.21 & 2.16 & 4.55 & 7.78  & 18.57 & 19.59 & 3.78 &  & 572.88 & 0.20 & 3.89  & 30.05  & 62.92  & 54.81  & 15.81 &  & 653.57 & 0.18 & 5.00  & 35.49  & 94.71  & 55.71  & 13.91 \\
0.6 &  & 40.77 & 1.95 & 4.55 & 7.20  & 16.50 & 17.47 & 3.28 &  & 619.16 & 0.13 & 3.16  & 25.85  & 53.47  & 47.73  & 13.93 &  & 719.02 & 0.21 & 4.22  & 29.40  & 79.18  & 46.53  & 11.36 \\
0.7 &  & 48.94 & 4.30 & 3.89 & 6.54  & 14.88 & 16.13 & 2.95 &  & 653.56 & 0.03 & 2.92  & 23.12  & 48.11  & 42.54  & 12.54 &  & 807.35 & 0.13 & 3.34  & 21.34  & 60.55  & 35.66  & 8.25  \\
0.8 &  & 56.21 & 5.86 & 3.39 & 6.05  & 13.55 & 15.12 & 2.70 &  & 674.71 & - & 3.15  & 21.58  & 44.99  & 39.78  & 11.83 &  & 899.40 & 0.06 & 2.33  & 14.09  & 43.81  & 25.20  & 5.37  \\
0.9 &  & 59.72 & 7.05 & 3.67 & 6.00  & 12.57 & 14.35 & 2.65 &  & 700.32 & - & 2.68  & 19.47  & 41.97  & 37.02  & 11.01 &  & 943.86 & 0.05 & 1.87  & 11.29  & 36.18  & 20.59  & 4.15  \\
1.0 &  & 62.87 & 7.73 & 4.14 & 5.71  & 11.99 & 14.01 & 2.49 &  & 763.03 & - & 2.00  & 15.82  & 34.51  & 30.22  & 9.11  &  & 984.88 & - & 1.53  & 8.93   & 29.71  & 16.78  & 3.33 \\
\bottomrule
\end{tabular}
\end{adjustbox}
\caption{The average amount of procured relief items $\bar{f}_t$ (see~\eqref{eq:avg_procurement}) obtained from policies associated with the \ac{famspr} model as a function of the cost-scaling factor $\nu$ for different time periods $t=1, \dots, T$, and for different initial intensity levels $\alpha_1 \in \{1, 3, 5\}$.}
\label{tab:sens_results_random}
\end{table}
\begin{table}[htbp]
    \centering
    \begin{adjustbox}{width=\textwidth}
        \begin{tabular}{@{} lll @{}} 
            \toprule
            Parameter         & Description & Value/Formula \\ \midrule
            $|J|$      			& number of \ac{dps}   		& $\in \{10, 20, 30\}$.\\
            $|I|$      			& number of \ac{sps}   		& $\in \{3, 6, 9\}$.\\
            $j_x$      			& $x$-coordinates of the \ac{dps}   		& $\in U(0,700)$.\\
            $j_y$      			& $y$-coordinates of the \ac{dps}   		& $\in U(0,100)$.\\
            $i_x$      			& $x$-coordinates of the \ac{sps}   		& $\in U(0,700)$.\\
            $i_y$      			& $y$-coordinates of the \ac{sps}   		& $ \in U(100,200)$.\\
            ($0_x,0_y$)      		& the $x$ and $y$ coordinates of the \ac{mdc}                                     		& $(350,450)$.\\
            $\mathrm{A}$	& state space for $\alpha_t$ (hurricane intensity) 		& $\{0, 1, \dots, 5\}$.\\
            $L_{x}$	& state space for $\ell_{x,t}$ (landfall spatial dimension) 		& $\{[0,100), \dots, [600,700)\}$.\\
            $L_{y}$	& state space for $\ell_{y,t}$ (landfall temporal dimension) 	& $\{[-350,-300),\dots, [-50,0), [0,\infty)\}$.\\
            $\mathbf{P}^{\alpha}$    & one-step transition probability matrix for $\alpha_t$	& see~\eqref{eq:intensity_MC}.\\
            $\mathbf{P}^{\ell_x}$    & one-step transition probability matrix for $\ell_{x,t}$	& see~\eqref{eq:hlocation_MC}.\\
            $\mathbf{P}^{\ell_y}$    & one-step transition probability matrix for $\ell_{y,t}$	& see~\eqref{eq:vlocation_MC}.\\
            $\mathcal{T}$ & a set of transient states when $T$ is random & $\{\xi_t = (\alpha_t, \ell_{x,t}, \ell_{y,t}) \; | \; \alpha_t \neq 0, \ell^{-}_{y,t} < 0\}$. \\
            $\mathcal{A}$ & a set of absorbing states when $T$ is random & $\{\xi_t = (\alpha_t, \ell_{x,t}, \ell_{y,t}) \; | \; \alpha_t = 0$ or $ \ell^{-}_{y,t} \geq 0\}$. \\

            $T_{\max}$          & maximum possible number of stages	& 8. \\
            $\floor{\hat T}$    & floor of the expected number of steps before reaching the absorbing state $\ell_{y} = [0,\infty)$ & 5.\\
            $M$                 & number of discrete points (exact coordinates) for a potential landfall location ($\ell_{x,t}$) & 10. \\
            $\ell^m_{x,t}$ &  exact coordinates for a potential hurricane spatial landfall location at a given state $\ell_{x,t}$ & see~\eqref{eq:exact_coor}.\\
            $d^{\xi_t}_{j,t}$ & demand at DP $j\in J$ for given intensity $\alpha_t$, a horizontal location $\ell^{m}_{x,t}$, and vertical location $\ell_{y,t}$ & see~\eqref{eq:demand_function}.\\
            $\mathbb{P}(\ell^m_{x,t} | \ell_{x,t} = \ell)$ & probability of hurricane landfall in coordinates $\ell^m$ given state $\ell_{x,t}=\ell$ & $1/M$.\\
            $\bar d$            & maximum possible demand for relief items at a given DP & 400. \\
            $\bar \delta$       & maximum possible distance to a DP which generates a request for relief items & 300. \\ 
            $\omega$ & fuel cost -- used in calculating the transportation costs & 0.0038. \\
            $\nu$ & logistic costs scaling factor used in the main results & $\in \{0.001, 0.60, 5.00\}$.\\
            $c^b_{ii',t}$ & unit cost for prepositioning/rerouting the relief items from an SP $i \in \{0\} \cup I$ to an an SP $i' \in I$ & see~\eqref{eq:ship_cost}.\\
            $c_{ij,t}^a$ & unit cost for delivering the relief items from an SP $i \in I$ to a DP $j \in J$ & see~\eqref{eq:deliver_cost}.\\
            $\beta$ & base cost used in the calculation of $h, c_t^h, p$ and $q$ & $5$.\\
            $c^h_{i,t}$ & unit cost for holding a relief item at an SP $i\in I$ in period $t$ & $0.2\beta$.\\
            $h_t$ & unit cost for procuring a relief item from the \ac{mdc} in period $t$ & $\beta (1+\nu(t-1))$.\\
            $p$ & unit cost for demand shortage (penalty for unsatisfied demand per item) & $80\beta$.\\
            $q$ &  salvage value per item & -$0.05\beta$.\\
            ${u_i}$ & inventory capacity at SP $i \in I$ & $\in U(0.05\bar{d}\frac{|J|}{|I|}, 0.5\bar{d}\frac{|J|}{|I|})$.\\
            \midrule
            $\hat{z}$ & sample mean for the policy performance & $\sum_{n=1}^N z(\hat{\xi}^n)/N$.\\
            $\hat{\sigma}$ & sample standard deviation for the policy performance & $\sqrt{\sum_{n=1}^N \left(z(\hat{\xi}^n)-\hat{z}\right)^2/(N-1)}$.\\
            $\hat{z}_{\pm}$ & $95\%$ confidence level for the mean policy performance & $\hat{z}\pm 1.96 \hat{\sigma}/\sqrt{{N}}$. \\
            ${x}_{0}$ & initial inventory levels at different \ac{sps} &  $\mathbf{0}$. \\
            $\xi_1$ & initial state in the initial experiments & $(1, [100,200), [-350,-300))$\\
            $N$ & number of out-sample paths used to evaluate different policies & 1000.\\
            $R$ & number of scenario sample paths used in the training of \ac{2ssp} models & 100.\\
            $\bar n $ & maximum number of iterations used during the training of different models  & $10^5$. \\
            $\epsilon$ & tolerance parameter & $10^{-5}$.\\  
            $\bar{m}$ & a stalling parameter used in the termination criterion & $500$.\\
            time limit & the time limit used during the training of different models & $3\times 60^2$ seconds.\\
            \bottomrule
        \end{tabular}
    \end{adjustbox}
    \caption{Summary of the parameters and values used in the problem data and implementation.}
    \label{tab:data}
\end{table}
\end{appendix}

\end{document}